\documentclass[onefignum,onetabnum]{siamart220329}

\DeclareSymbolFont{CMlargesymbols}{OMX}{cmex}{m}{n} 
\DeclareMathDelimiter{(}{\mathopen} {operators}{"28}{CMlargesymbols}{"00}
\DeclareMathDelimiter{)}{\mathclose}{operators}{"29}{CMlargesymbols}{"01}
\DeclareMathAlphabet\mathcal{OMS}{cmsy}{m}{n} 
\SetMathAlphabet\mathcal{bold}{OMS}{cmsy}{b}{n} 

\usepackage{listings}
\usepackage{lipsum}
\usepackage{amsfonts}
\usepackage{graphicx}
\usepackage{epstopdf}
\usepackage{algorithmic}
\ifpdf
  \DeclareGraphicsExtensions{.eps,.pdf,.png,.jpg}
\else
  \DeclareGraphicsExtensions{.eps}
\fi

\usepackage[numbers,sort&compress]{natbib}

\newsiamremark{remark}{Remark}
\newsiamremark{hypothesis}{Hypothesis}
\crefname{hypothesis}{Hypothesis}{Hypotheses}
\newsiamthm{claim}{Claim}

\headers{HPS on triangular and deforming surfaces}{ G.~Zavalani}

\title{Fast High-Order Spectral Solvers for PDEs on Triangulated Surfaces
with Applications to Deforming Surfaces
\thanks{Submitted to the editors DATE.
}
}

\author{Gentian Zavalani  \thanks{Technische Universit{\"a}t Dresden, %
  Institute of Numerical Mathematics, %
  01062 Dresden, Germany. \\
Email: \url{gentian.zavalani@tu-dresden.de}
}
}

\usepackage{amsopn}

\ifpdf
\hypersetup{
  pdftitle={Fast spectral methods on triangular and deforming surfaces},
  pdfauthor={G. Zavalani}
}
\fi

\usepackage{multirow} 
\usepackage{mathptmx}
\usepackage[mathcal]{eucal}
\usepackage{nccmath, amssymb, mathtools}
\usepackage[utf8]{inputenc}
\usepackage[T1]{fontenc}
\usepackage{float}
\usepackage{overpic}
\usepackage{mathtools}
\usepackage{tgpagella}
\usepackage{multirow}
\usepackage{float}
\usepackage[parfill]{parskip}

\usepackage{graphicx}
\usepackage{lettrine}
\usepackage{enumitem}
\usepackage{mathrsfs}
\usepackage{tikz}
\usetikzlibrary{matrix}
\usepackage{tikz-cd}
\usepackage{booktabs}
\usepackage{caption}
\usepackage{subcaption}
 \definecolor{highlightcolor}{rgb}{0.7, 0.85, 1.0} 

\usepackage{microtype}
\usepackage{bm}
\usepackage{enumitem}
\usepackage{url}
\usepackage{xcolor}
\usepackage{hyperref}
\usepackage{bookmark} 

\definecolor{OurRed}{rgb}{0.64, 0.30, 0.30}

\graphicspath{{images/}}

\makeatletter
\def\input@path{{images/}}
\makeatother

\newcommand{\Cheb}{\mathrm{Cheb}}

\newcommand{\R}{\mathbb{R}}
\newcommand{\surf}{\Gamma}

\newcommand{\dt}{\Delta t}

\newcommand{\p}{\partial}

{\bf}{\it}

\makeatletter
\let\cref@override@label@type\@gobbletwo
\makeatother

\begin{document}

\maketitle


\begin{abstract}
In this paper, we extend the classical quadrilateral based hierarchical
Poincaré--Steklov (HPS) framework to triangulated geometries. Traditionally, the
HPS method takes as input an unstructured, high-order quadrilateral mesh and
relies on tensor-product spectral discretizations on each element. To overcome
this restriction, we introduce two complementary high-order strategies for
triangular elements: a
reduced quadrilateralization approach which is
straightforward to implement, and triangle based spectral element method based on Dubiner
polynomials. We show numerically that
these extensions preserve the spectral accuracy, efficiency, and fast
direct-solver structure of the HPS framework. The method is further extended to
time dependent and evolving surfaces, and its performance
is demonstrated through numerical experiments on reaction--diffusion systems, and geometry driven surface evolution.
\end{abstract}
\begin{keywords}
spectral methods; HPS solver; domain decomposition; Dubiner
polynomials; triangulated meshes; Turing patterns; evolving surfaces.
\end{keywords}

\section{Introduction}
Partial differential equations (PDEs) on static and evolving surfaces are commonly discretized using finite difference, finite element, and spectral methods. Finite difference schemes are simple to implement but offer only algebraic convergence, requiring many grid points and substantial memory. Surface finite element methods (SFEMs)~\cite{dziuk2013finite} provide geometric flexibility but rely on low-order polynomial approximations—typically cubic or quartic~\cite{jonsson2017cut}—which can lead to ill-conditioned systems and make high-order accuracy difficult to achieve despite the availability of efficient solvers such as multigrid~\cite{briggs2000multigrid} and UMFPACK~\cite{davis2004algorithm}. Moreover, when the underlying surface mesh has limited smoothness, attaining high-order accuracy remains particularly challenging. Spectral methods, by contrast, attain very high accuracy with relatively few degrees of freedom by employing high-degree polynomial approximations, and spectral element methods extend this accuracy to flexible geometries~\cite{fornberg1998practical}, though at the cost of potentially dense linear systems. When combined with domain-decomposition techniques such as Schwarz iterations~\cite{canuto1988schwarz} or hierarchical Poincaré–Steklov (HPS) solvers~\cite{gillman2014direct, martinsson2013direct}, the overall complexity can be reduced to $\mathcal{O}(N^{3/2})$ or even $\mathcal{O}(N)$. In this work, we adopt such a domain-decomposition strategy with spectral collocation on each element,\footnote{Other element-wise discretizations, such as finite elements, could also be employed.} together with a hierarchical direct solver based on nested dissection~\cite{george1973nested}, yielding the HPS method~\cite{martinsson2019fast}. This framework eliminates local degrees of freedom through Dirichlet-to-Neumann or Schur-complement operators and recursively merges them in a binary tree, producing an efficient and robust direct solver without assembling a global stiffness matrix.

To date, HPS has been most effective on quadrilateral surface meshes because spectral collocation is remarkably efficient on tensor-product domains due to the inherent structure of the expansions they employ. The main drawback, however, is
the inability to handle complex geometries. In order to handle highly complex geometries the use of triangular elements
is generally preferred. Therefore, it is tempting to try to marry the efficiency of tensor products with the flexibility of triangular geometries.  Although prior work has remarked that the HPS framework could, in principle, be adapted to triangular elements, to the best of our knowledge, no high-order formulation or implementation has been presented. In this work, we propose an extension of the HPS framework to triangulated surface discretizations. 
Building on its formulation for quadrilateral meshes, we introduce two complementary high-order strategies tailored to triangular geometries: a quadrilateralization approach that embeds each triangular element into a quadrilateral patch compatible with the existing HPS formulation, and a basis-construction approach that employs Dubiner polynomials to enable direct spectral collocation on triangular elements. A guiding objective of the paper is to provide a detailed and self-contained description of the
quadrilateral based HPS framework and to identify the modifications required for a
stable and high-order extension to triangular elements.

The remainder of the paper is organized as follows. Section~\ref{sec:IP} introduces the multivariate interpolation framework on reference elements, including tensor-product interpolation on the hypercube and total degree interpolation on simplices. In Sections~\ref{sub:local_discret} and~\ref{sub:local_op} we revisit the quadrilateral based HPS domain decomposition framework~\cite{fortunato2022highorder, martinsson2013direct}. Unlike previous treatments, we provide a detailed construction of high-order surface parametrizations and collocation operators, as these components are essential for the triangular extensions introduced later. Section~\ref{domain_compos} summarizes the domain decomposition procedure. Section~\ref{sec:extend_rohmbus} presents two strategies for extending the HPS
method to triangular surface meshes: a reduced quadrilateralization procedure
and a triangle based spectral element method based on Dubiner
polynomials.  Section~\ref{time_depdn} extends the framework to time-dependent problems, and Section~\ref{turing_model} illustrates its performance through numerical experiments involving interacting Turing systems.

Section~\ref{sec:evol_surface} discusses how the framework can be applied to
evolving surfaces by introducing a time-dependent surface representation and
extending the projection method to account for geometric evolution. Section~\ref{sec:surf_app} introduces triangulated surface evolution and motivates the use of an arbitrary Lagrangian–Eulerian (ALE) framework to control mesh distortion. Finally, Sections~\ref{sec:surface_growth} and~\ref{patter_evolve} examine two geometric growth models and study their influence on the formation, distribution, and symmetry of patterns over time.

\section{Multivariate polynomial interpolation on reference elements}\label{sec:IP}

We consider polynomial interpolation on the tensor-product hypercube
 $\Omega_d=[-1,1]^d$ and on the $d$-simplex $\Delta_d \subset
\mathbb{R}^d$, using Lagrange bases defined on appropriate nodal sets.

\paragraph{Polynomial spaces and multi-index notation}

For $x = (x_1,\dots,x_d)\in \mathbb{R}^d$ and
$\alpha = (\alpha_1,\dots,\alpha_d)\in \mathbb{N}_0^d$, we write
\[
  x^\alpha = \prod_{i=1}^d x_i^{\alpha_i}, \qquad
  \|\alpha\|_\infty = \max_{1\le i\le d} \alpha_i, \qquad
  \|\alpha\|_1 = \alpha_1 + \cdots + \alpha_d.
\]

On the hypercube $\Omega_d$, we consider the tensor product polynomial space
\[
  \Pi_{d,n}(\Omega_d)
  = \mathrm{span}\{x^\alpha : \alpha \in A_{d,n}\},
  \qquad 
  A_{d,n} = \{\alpha \in \mathbb{N}_0^d : \|\alpha\|_\infty \le n\},
\]
consisting of polynomials of degree at most $n$ in each coordinate.

On the simplex $\Delta_d$, approximation is based instead on the total degree polynomial space

\[
  \Pi_{d,n}(\Delta_d)
  = \mathrm{span}\{x^\alpha : \alpha \in \mathbb{N}_0^d : \|\alpha\|_1 \le n\},
\]
whose dimension is 
\[
  N_n = \binom{n+d}{d}.
\]

\paragraph{Tensor-product interpolation on the hypercube}

Let $P_i = \{p_{0,i},\dots,p_{n,i}\} \subset [-1,1]$ denote a univariate 
interpolation grid with $|P_i| = n+1$. The tensor-product grid
\[
  G_{d,n} = P_1 \times \cdots \times P_d
\]
contains $(n+1)^d$ nodes. Each multi-index
$\alpha = (\alpha_1,\dots,\alpha_d)\in A_{d,n}$ corresponds to a grid point
\[
  p_\alpha 
  = (p_{\alpha_1,1},\dots,p_{\alpha_d,d}) \in G_{d,n}.
\]

For $i=1,\dots,d$, let $\ell_{j,i}$ denote the univariate Lagrange basis
functions associated with $P_i$. The resulting tensor-product Lagrange
basis on $G_{d,n}$ is
\begin{equation*}\label{eq:L}
L_{\boldsymbol{\alpha}}(x_1,\dots,x_d)
=
\prod_{i=1}^d l_{\alpha_i,i}(x_i),
\qquad
l_{j,i}(x)
=
\prod_{\substack{k=0\\k\ne j}}^{n}
\frac{x-p_{k,i}}{p_{j,i}-p_{k,i}}.
\end{equation*}
which satisfies $L_\alpha(p_\beta) = \delta_{\alpha\beta}$.

The interpolation operator
$Q_{G_{d,n}} : C^0(\Omega_d)\to \Pi_{d,n}(\Omega_d)$
is defined by
\[
  (Q_{G_{d,n}}f)(p_\alpha) = f(p_\alpha),
  \qquad \alpha \in A_{d,n},
\]
and admits the representation
\begin{equation}\label{main_poly}
  Q_{G_{d,n}} f(x)
  = \sum_{\alpha\in A_{d,n}} f(p_\alpha)\, L_\alpha(x).
\end{equation}

In this work we use the Chebyshev--Lobatto grid
\[
  \mathrm{Cheb}_n = \bigl\{\cos(k\pi/n): 0\le k\le n \bigr\},
  \qquad
  \mathrm{Cheb}_{d,n} = \bigotimes_{i=1}^d \mathrm{Cheb}_n,
\]
which is standard in spectral collocation methods.

\paragraph{Total-degree interpolation on the simplex}

For elements lacking tensor-product structure, such as simplices,
interpolation is based on the total degree space $\Pi_{d,n}(\Delta_d)$.$\qquad$  
Let
  $\widehat X_n = \{\hat x_m\}_{m=0}^{N_n} \subset \Delta_d$,
be any unisolvent nodal set~\cite{isaac2020recursive}, and let 
$\{\phi_k\}_{k=0}^{N_n}$ denote a basis for $\Pi_{d,n}(\Delta_d)$,
typically Dubiner polynomials~\cite{dubiner1991spectral}. The corresponding Vandermonde matrix is
\[
  K_{mk} = \phi_k(\hat x_m),
  \qquad 0 \le m,k \le N_n.
\]
Since the nodes are unisolvent, $K$ is invertible. The Lagrange basis
functions are then
\[
  \ell_m(x) = \sum_{k=0}^{N_n} (K^{-1})_{km}\,\phi_k(x),
  \qquad m = 0,\dots,N_n,
\]
and satisfy $\ell_m(\hat x_{m'}) = \delta_{mm'}$.

The associated interpolation operator
$Q_{\Delta_d,n}: C^0(\Delta_d)\to \Pi_{d,n}(\Delta_d)$
is given by
\begin{equation}  \label{eq:TDinterp}
  Q_{\Delta_d,n} f(x)
  = \sum_{m=0}^{N_n} f(\hat x_m)\,\ell_m(x).
\end{equation}

Together, these approximation spaces form the basis for the discretization
of quadrilateral and triangular surface elements, respectively.

\section{High-order parametric surface approximation}\label{sub:local_discret}
We consider a general elliptic surface PDE on $\Gamma$,

\begin{equation}\label{eq:problem}
\mathcal{L}_{\Gamma} u(\mathbf{x}) = f(\mathbf{x}), \quad \mathbf{x} \in \Gamma,
\end{equation}

where $f(\mathbf{x})$ is a smooth function on $\Gamma$ and $\mathcal{L}_{\Gamma}$ is a variable-coefficient linear second-order elliptic surface operator. If $\Gamma$ is not a closed surface, Eq.~\eqref{eq:problem} may also be subject to boundary conditions, e.g., $u(\mathbf{x}) = h(\mathbf{x})$ for $\mathbf{x} \in \partial \Gamma$ and some function $h$ .

To numerically solve Eq.~\eqref{eq:problem}, we need to discretize $\mathcal{L}_{\Gamma}, f$, and $h$, and represent these
objects in some finite-dimensional basis. The numerical error involved in a discretization of PDEs on curved surfaces depends on two properties of the discretization, the representation of the objective function and the representation of the geometry. Thus, a higher-order scheme is only possible with also a higher-order description of the surface approximation.
Let $\Gamma^{\text{quad}}_h \subset \mathcal{N}_{\delta}$\footnote{
The open tubular neighborhood of the normal bundle is defined as $\mathscr{N}_\delta := \left\{ \mathbf{x} \in \mathbb{R}^{d+1} \,\middle|\, d(\mathbf{x}, \Gamma) < \delta \right\},$
where \( d : \mathscr{N}_\delta \rightarrow \Gamma \) denotes the signed distance function.
} be a shape-regular,  composed of finitely many regular and quasi-uniform $d$-dimensional  hypercubes with diameter $h$,  topologically equivalent to the smooth surface $\Gamma$. The collection of these hypercubes is denoted by $\widehat{\mathcal{K}}_h$, which provides a representation of $\Gamma^{\text{quad}}_h$:

\begin{equation}
\Gamma^{\text{quad}}_h = \bigcup_{\widehat{\mathcal{E}} \in \widehat{\mathcal{K}}_h} \widehat{\mathcal{E}},\; \text{with}\; K=|\widehat{\mathcal{K}}_h|.
\end{equation}

 We assume that the elements do not overlap, i.e., for \(\widehat{\mathcal{E}}_i, \widehat{\mathcal{E}}_j \in \widehat{\mathcal{K}}_h\), we have that \(\operatorname{int}(\widehat{\mathcal{E}}_i) \cap \operatorname{int}(\widehat{\mathcal{E}}_j) = \emptyset\) and if \(\widehat{\mathcal{E}}_i \cap \widehat{\mathcal{E}}_j = I \neq \emptyset\) and \(\dim(I) = d-1\), it is called an intersection of \(\widehat{\mathcal{E}}_i\) and \(\widehat{\mathcal{E}}_j\) and is assumed to be a subset of an \((d-1)\)-dimensional facet of \(\widehat{\mathcal{E}}_i\) and \(\widehat{\mathcal{E}}_j\), respectively. Each element $\widehat{\mathcal{E}} \in \widehat{\mathcal{K}}_h$ is parametrized over a reference element $\Omega_{d} \subset \mathbb{R}^d$ by an invertible and differentiable mapping:

\begin{equation}
\mu_{\widehat{\mathcal{E}}} : \Omega_{d} \to \widehat{\mathcal{E}},
\end{equation}

referred to as the geometry mapping of $\widehat{\mathcal{E}}$.

Additionally, we assume the existence of a bijective mapping $\mathbf{X}: \Gamma^{\text{quad}}_h \to \Gamma$, such that the smooth surface $\Gamma$ can be represented as a union of non-overlapping mapped elements:

\begin{equation}
\Gamma = \bigcup_{\widehat{\mathcal{E}}  \in \widehat{\mathcal{K}}_h} \mathbf{X}(\widehat{\mathcal{E}} ) = \bigcup_{\widehat{\mathcal{E}}  \in \widehat{\mathcal{K}}_h} \mathbf{X}(\mu_{\widehat{\mathcal{E}}} (\Omega_{d})) =: \bigcup_{\widehat{\mathcal{E} } \in \widehat{\mathcal{K}}_h} \mathbf{X}_{\widehat{\mathcal{E}}}(\Omega_{d}).
\end{equation}

With this property, $\Gamma^{\text{quad}}_h$ is referred to as the reference surface (or reference domain) of $\Gamma$, and the collection $\{\mathbf{X}_{\widehat{\mathcal{E}}}\}_{\widehat{\mathcal{E}} \in \widehat{\mathcal{K}}_h}$ as its reference parametrization.

 \begin{figure}[t!]
    \centering
    \begin{tikzpicture}
        \node[inner sep=0pt] at (0,0) {\includegraphics[clip,width=1.0\columnwidth]{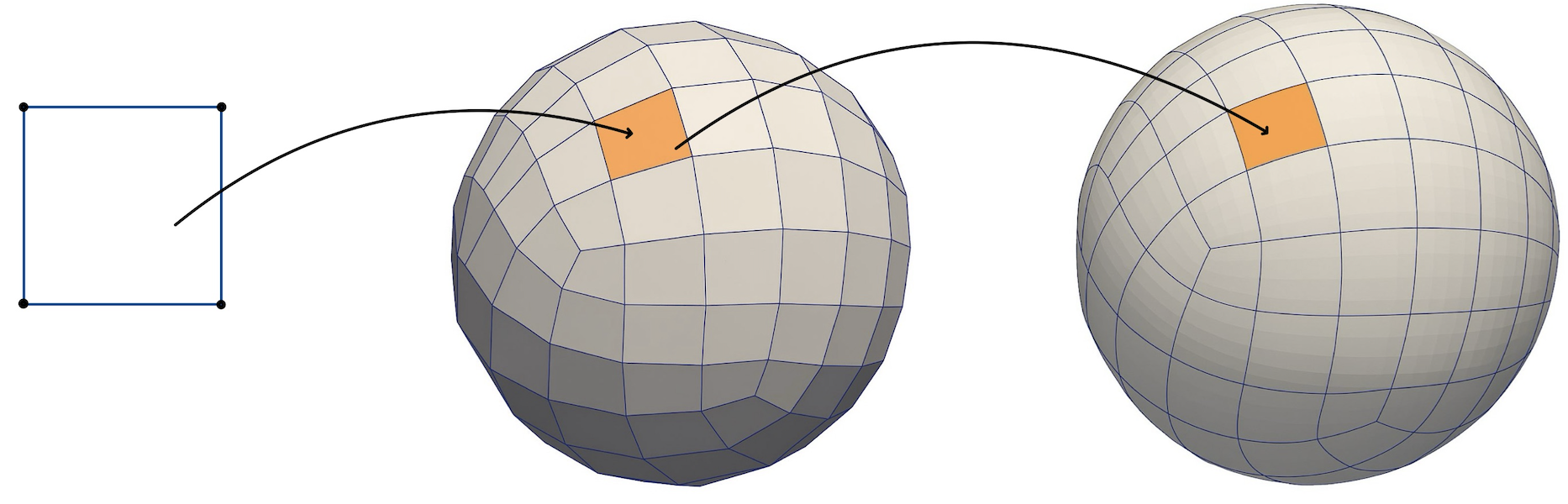}};
        
    
        \node[anchor=north west] at (1.5,1.6) {$\pi_k$};

        \node[anchor=north west] at (4.7,1.1) {$\mathcal{E}_k\,$};

          \node[anchor=north west] at (4.7,3.3) {$\Gamma^{\text{quad}}_{h,n}$};

        \node[anchor=north west] at (1.32,3.0) {$\psi_{k}:=\pi_{k}\circ\mu_{k}$};
        \node[anchor=north west] at (-2.,1.1) {$\widehat{\mathcal{E}}_k$};
        \node[anchor=north west] at (-1.8,3.1) {$\Gamma^{\text{quad}}_h$};
        \node[anchor=north west] at (-3.6,1.7) {$\mu_{k}$};
        \node[anchor=north west] at (-6.3,0.3) {$[-1,1]^2$};
    \end{tikzpicture}
  \vspace{-20pt} 
  \caption{Construction of a surface parametrization over the reference square 
    $[-1,1]^2$ via closest-point projection from the piecewise affine approximation 
    $\Gamma^{\mathrm{quad}}_{h}$.}
  \label{Figure_quad_app_frame}
\end{figure}
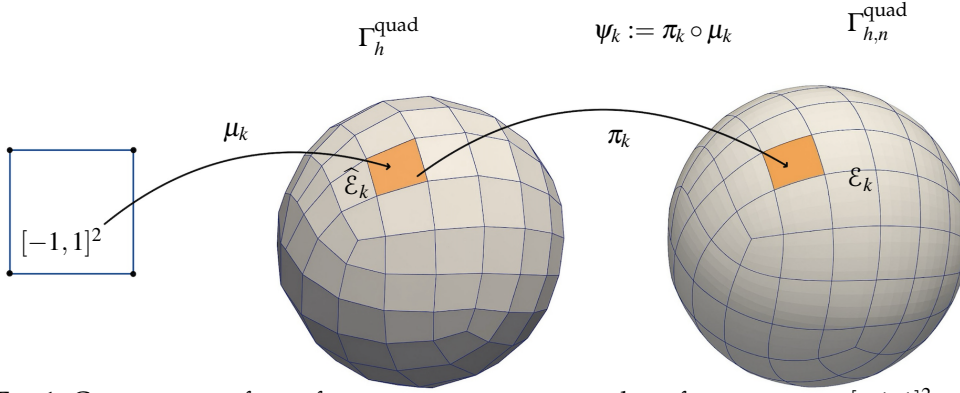 

The Jacobian of the parametrization $\mathbf{X}_{\widehat{\mathcal{E}}}$ at $\mathbf{x}$ is  
\[
D\mathbf{X}_{\widehat{\mathcal{E}}}(\mathbf{x}) : \mathbb{R}^d \to \mathbb{R}^{d+1}.
\]  

The first fundamental form is the symmetric and positive definite matrix $\mathbf{g} \in \mathbb{R}^{d \times d}$ defined as  
\[
\mathbf{g}(\mathbf{x}) := D\mathbf{X}_{\widehat{\mathcal{E}}}(\mathbf{x})^T D\mathbf{X}_{\widehat{\mathcal{E}}}(\mathbf{x}) \quad \forall \mathbf{x} \in \Omega_{d}.
\]  

The corresponding volume element $g(\mathbf{x})$ is given by  
\[
g(\mathbf{x}) = \sqrt{\det \mathbf{g}(\mathbf{x})}.
\]

Let  $Q_{G_{d,n}}\mathbf{X} \in \Pi_{d,n}$ be the $n$th-order Lagrange polynomial interpolation of the mapping $\mathbf{X}$ on the reference element $\Omega_{d}$. In the following, we consider the restriction of the closest-point projection $\pi$ to $\Gamma^{\text{quad}}_h$ as a mapping $\pi : \Gamma^{\text{quad}}_h \to \mathbb{R}^{d+1}$, and define $n^{\text{th}}$ -order cubic parameterization of the surface. 

\begin{definition}[$n^{\text{th}}$-order parametrization]\label{def:46}
Given an  $C^r$-regular cubical parametrization $\psi_{k}:=\pi_{k}\circ\mu_{k} : \Omega_{d} \to \widehat{\mathcal{E}}_k,\;
k=1,\dots,K,\, r\geq2
\,,$ where $\pi_{k}$ is the closest point projection onto the target surface $\Gamma$ and $\mu_{k}$ is a multi-linear map (see Figure~\ref{Figure_quad_app_frame}).
We say that the mesh \( \Gamma^{\mathrm{quad}}_h \) is of order \( n \) if each element is obtained by polynomial interpolation of the values \( \{ \psi_k(p_{\boldsymbol{\alpha}}) \}_{\boldsymbol{\alpha} \in A_{d,n}} \) sampled at the points \( \{ p_{\boldsymbol{\alpha}} \}_{\boldsymbol{\alpha} \in A_{d,n}} \).
\end{definition}

This implies that on each element, we can approximate the $C^r$-regular parametrization maps through interpolation using the points $\{\psi_{k}(p_{\boldsymbol{\alpha}})\}_{{\boldsymbol{\alpha}}\in A_{d,n}}$. This involves computing a $n^{\text{th}}$-order vector-valued polynomial approximation:
\begin{equation}\label{main.poly}
    Q_{G_{d,n}}\psi_{k}\left(\mathrm{x}\right) =\sum_{{\boldsymbol{\alpha}} \in A_{d,n}}\psi_{k}(p_{\boldsymbol{\alpha}})L_{{\boldsymbol{\alpha}}} ,\quad k=1,\ldots, K .
\end{equation}
Partial derivatives of the elemental coordinate maps, $\partial_{\xi_j} Q_{G_{d,n}} \psi_k$, for $k = 1,\ldots,K$ and $j = 1,\ldots,d$, are used to construct an $n^{\text{th}}$-order approximation of the volume element $g(\mathbf{x})$

\begin{equation}\label{main.poly1}
g^n(\mathbf{x}) = \sqrt{\det \left(\mathbf{D}Q_{G_{d,n}}\psi_{k}(\mathbf{x})^T \mathbf{D}Q_{G_{d,n}}\psi_{k}(\mathbf{x})\right)}.
\end{equation}

Mapping the piecewise flat surface results in a higher-order approximation $\Gamma^{\text{quad}}_{h,n} := Q_{G_{d,n}}\psi_{k}(\Gamma^{\text{quad}}_h)$ of the surface $\Gamma$ that is contained in $\mathcal{N}_{\delta}$ if the mesh size is small enough. In the following, we assume that $h$ is chosen correspondingly. Associated with the discrete surface $\Gamma^{\text{quad}}_{h,n}$ is a set of surface elements
\[
\mathcal{K}_{h,n} := \{ Q_{G_{d,n}}\psi_{k}(\widehat{\mathcal{E}}) \mid \widehat{\mathcal{E}} \in \widehat{\mathcal{E}}_h \},
\]
such that $\Gamma^{\text{quad}}_{h,n} = \bigcup_{\mathcal{E} \in \mathcal{K}_{h,n}} \mathcal{E}$.

\section{High-order spectral
collocation on a single element}\label{sub:local_op} 
We now describe a spectral collocation method for discretizing Eq.~\eqref{eq:problem} on a single surface element $\mathcal{E}_{k}\subset\Gamma$. The discretization is derived directly from the strong formulation of the problem. Let $p_{{\boldsymbol{\alpha}}}\in\Cheb_{d,n},\; {\boldsymbol{\alpha}}\in A_{d,n}$ denote points in a tensor product Chebyshev-Lobatto grid. Let $\mathbf{u}\in \R^{(n+1)^d}$  denote a vector holding approximations to $u$ at $p_{{\boldsymbol{\alpha}}}\in \Cheb_{d,n}$. Then just like Eq.~\eqref{main.poly}, we have:
\begin{equation}
        Q_{G_{d,n}}\mathbf{u}\left(\mathrm{x}\right) =\sum_{{\boldsymbol{\alpha}} \in A_{d,n}}\mathbf{u}(p_{\boldsymbol{\alpha}})L_{{\boldsymbol{\alpha}}} ,\quad p_{{\boldsymbol{\alpha}}}\in\Cheb_{d,n}.
\end{equation}

To discretize $\mathcal{L}_{\Gamma}$ on the element $\mathcal{E}_k$, we first compute the discrete operators on the reference cube $\Omega_d$ and then map them to $\mathcal{E}_k$ using the numerical coordinate mapping $\mathbf{\psi}_k$.  Let \( \mathcal{D} \in \R^{(n+1) \times (n+1) }\) be the one-dimensional spectral differentiation matrix \cite{trefethen2000spectral} associated with Chebyshev-Lobatto points on the interval \([-1,1]\) and let \( I \in \R^{(n+1)\times (n+1)} \) be the identity matrix. The $d$-dimensional differentiation matrices on the reference hypercube \(\Omega_{d}\) in the \(\xi_1,\; \xi_1,\ldots, \xi_d\)-directions can be constructed through the Kronecker products

\begin{align*}
\mathcal{D}_{\xi_1} &= I \otimes \cdots \otimes I \otimes \mathcal{D} \\
\mathcal{D}_{\xi_2} &= I \otimes \cdots \otimes  \mathcal{D} \otimes I \\
\mathcal{D}_{\xi_i} &= I \otimes \cdots \otimes  \mathcal{D} \otimes  \cdots \otimes I \\
& \quad \text{(with } \mathcal{D} \text{ in the } (d-i+1) \text{ position)} \\
\vdots & \\
\mathcal{D}_{\xi_d} &= \mathcal{D} \otimes I \otimes  \cdots \otimes I
\end{align*}
where $\mathcal{D}_{\xi_i}\in \R^{(n+1)^d \times (n+1)^d},\; i=1,2,\ldots,d. $ Let \(\mathbf{M}[\mathbf{u}] \in \R^{(n+1)^d \times (n+1)^d}\) denote the diagonal multiplication matrix formed by placing the entries of \(\mathbf{u}\) along the diagonal. Differentiation matrices corresponding to the components of the surface gradient on \(\mathcal{E}_k\) are given by\footnote{For instance, the discrete Laplace–Beltrami operator may be discretized as $\Delta_{\Gamma} \approx \sum_{i=1}^{d+1}(\mathcal{D}_{x_{i}}^{\Gamma})^{2}$.}
\begin{equation}
        \mathcal{D}_{x_{j}}^{\Gamma} = \sum_{i=1}^{d}\mathbf{M}[\p_{x_{j}}\xi_i] \mathcal{D}_{\xi_i},\; j=1,2,\ldots,d+1\,.
\end{equation}
The operator \(\mathcal{L}_{\Gamma}\) given in Eq.~\eqref{eq:problem} is then discretized on the element \(\mathcal{E}_k\) as the \((n+1)^d \times (n+1)^d\) matrix
\begin{equation}
    \mathbf{L}_{\mathcal{E}_k} = \sum_{i=1}^{d+1} \sum_{j=i}^{d+1} \mathbf{M}[a_{ij}] \mathcal{D}_i^{\Gamma} \mathcal{D}_j^{\Gamma} + \sum_{i=1}^{d+1} \mathbf{M}[b_i] \mathcal{D}_i^{\Gamma} + \mathbf{M}[c],
\end{equation}
with all variable coefficients sampled on the grid points $p_{{\boldsymbol{\alpha}}}\in \Cheb_{d,n}$. In other words, we can think of the matrix $\mathbf{L}_{\mathcal{E}_k}$ as a discrete approximation to the differential operator $\mathcal{L}_{\mathcal{E}_k}\;$.
Using this discretization scheme, Eq.~\eqref{eq:problem} can be written as \((n+1)^d \times (n+1)^d\) linear system.
\begin{equation}\label{main:disc1}
  \mathbf{L}_{\mathcal{E}_k}\mathbf{u}= \mathbf{f}\;.
\end{equation}
To prepare for the imposition of boundary conditions, we partition the index set \( \{1, \dots, (n+1)^d \} \) for a given element \( \mathcal{E}_k \) into interior (\( I_{\operatorname{int}} \)) and boundary (\( I_{\partial} \)) subsets, so that
\[
\{1, 2, \ldots, (n+1)^d \} = I_{\operatorname{int}} \cup I_{\partial},
\]
with \( (n - 1)^d \) points in the interior and \( 2nd^{d - 1} \) on the boundary. This partition is used to divide both vectors and matrices into blocks. For instance, given a matrix \( \mathbf{A} \), the submatrix \( \mathbf{A}^{i,b} \) consists of the rows indexed by \( I_{\operatorname{int}} \) and columns indexed by \( I_{\partial} \).~The solution vector is written as
\begin{equation*}
\mathbf{u}_i = \mathbf{u}(I_{\operatorname{int}}) \quad \text{and} \quad \mathbf{u}_b = \mathbf{u}(I_{\partial}).
\end{equation*}
Partitioning the $\mathbf{L}_{\mathcal{E}_{k}}$ and reordering the degrees of freedom in Eq.~\eqref{main:disc1} in the order $\{I_{\operatorname{int}}, I_{\partial}\}$ (i.e. interior then boundary) gives a block linear
system,

\begin{equation}\label{main:disc2}
\begin{bmatrix}
\mathbf{L}_{\epsilon_k}^{i,i} & \mathbf{L}_{\epsilon_k}^{i,b} \\
\mathbf{L}_{\epsilon_k}^{b,i} & \mathbf{L}_{\epsilon_k}^{b,b}
\end{bmatrix}
\begin{bmatrix}
\mathbf{u}^i \\
\mathbf{u}^b
\end{bmatrix}
=
\begin{bmatrix}
\mathbf{f}^i \\
\mathbf{f}^b
\end{bmatrix}.
\end{equation}

To impose Dirichlet boundary conditions, we set $\mathbf{u} = \mathbf{h}^b$ on $\partial \mathcal{E}_k$, where $\mathbf{h}^b \in \mathbb{R}^{2dn^{d-1} \times 1}$ denotes the vector of boundary values of $h(x)$.  Substituting these into the system and eliminating the boundary unknowns \( \mathbf{u}_b \) via a Schur complement \cite{mathew2008domain} yields the reduced system:
\[
\mathbf{L}_{\mathcal{E}_k}^{i,i} \mathbf{u}^i = \mathbf{f}^i - \mathbf{L}_{\mathcal{E}_k}^{i,b} \mathbf{h}^b.
\]
The interior solution is then expressed as
\[
\mathbf{u}^i = \left(\mathbf{L}_{\mathcal{E}_k}^{i,i}\right)^{-1} \mathbf{f}^i - \mathbf{S}_{\mathcal{E}_k} \mathbf{h}^b,
\]
where $\mathbf{S}_{\mathcal{E}_k} := -\left(\mathbf{L}_{\mathcal{E}_k}^{i,i}\right)^{-1} \mathbf{L}_{\mathcal{E}_k}^{i,b}$ is known as the solution operator and is of size $(n+1)^d \times 2d n^{d-1}$.



\section{Domain decomposition methods}\label{domain_compos}

The spectral collocation method described in the previous section converges very quickly as the number of points \( n \) increases provided the solution \( u \) of the Eq.~\eqref{eq:problem} is smooth. The matrix \( \mathbf{L}_{\mathcal{E}_k} \) that arises from the discretization has some structure and contains many zeros, but it is still considerably denser than the matrices produced by finite difference or finite element methods. One way to reduce this density is to use domain decomposition methods, which also lend themselves well to implementation on parallel computing architectures.

In a domain decomposition approach, the computational domain \( \Gamma \) is partitioned into smaller elements \( \mathcal{E}_k \), for \( k = 1, \ldots, K \), which may touch or overlap. The original problem \eqref{eq:problem} is then reformulated on each element, resulting in a family of smaller subproblems, where  each subproblem can be solved independently using spectral collocation. However, to ensure that the local solutions $u_k(\mathbf{x})$, each defined solely on a element $\mathcal{E}_k$ for $k = 1, \ldots, K$, fit together and form a
smooth solution of the PDE \eqref{eq:problem} on the entire computational domain $\Gamma$, they have to satisfy matching conditions. For two adjacent (non-overlapping) elements $\mathcal{E}_1$ and $\mathcal{E}_2$, the solution is required to be $C^1$-continuous across their common interface. In particular, both the solution and its binormal derivative must be continuous along the shared boundary. This requirement can be expressed as
    
\begin{equation}\label{eq:contin}
     u_1 (\mathbf{x}) = u_2 (\mathbf{x}), \quad \mathbf{x} \in \partial \mathcal{E}_{1}\cap\partial \mathcal{E}_{2}
\end{equation}
    \begin{equation}\label{eq:binormal}
     \partial_{\mathbf{n}_b}u_1(\mathbf{x}) =- \partial_{\mathbf{n}_b} u_2(\mathbf{x}), \quad \mathbf{x} \in \partial \mathcal{E}_{1}\cap\partial \mathcal{E}_{2}
\end{equation}
where $\mathbf{n}_b$ denotes the outward-pointing binormal vector. The minus sign arises from the opposing orientations of the binormals on the two adjacent elements.

The use of Chebyshev basis functions naturally leads to shared boundary degrees of freedom between neighboring quadrilateral elements, making the enforcement of continuity of the solution straightforward. Enforcing continuity of the binormal derivative is more delicate; among the available approaches, an effective strategy is to employ the Poincar\'e--Steklov operator~\cite{quarteroni2008numerical}, also known as the Dirichlet-to-Neumann (DtN) map. Originally introduced by V.~A.~Steklov, this operator maps prescribed Dirichlet data (solution values) to the corresponding Neumann data (binormal derivatives) on the boundary.



 For a given element \( \mathcal{E}_k \), the Dirichlet-to-Neumann operator, denoted by \( \mathrm{DtN}_{\mathcal{E}_k} \), computes the outward fluxes corresponding to prescribed Dirichlet boundary data.  Having constructed the solution operators \( \mathbf{S}_{\mathcal{E}_k} \) that solve the PDE locally on each element \( \mathcal{E}_k \), for \( k = 1, \ldots, K \), the Dirichlet-to-Neumann operator is given by the product of the binormal derivative operator \( \mathcal{D}_{\mathcal{E}_k} \) and the solution operator \( \mathbf{S}_{\mathcal{E}_k} \):
\[
\mathrm{DtN}_{\mathcal{E}_k} = \mathcal{D}_{\mathcal{E}_k} \, \mathbf{S}_{\mathcal{E}_k}.
\]
\subsection{Merging Dirichlet-to-Neumann maps}\label{merge_section}

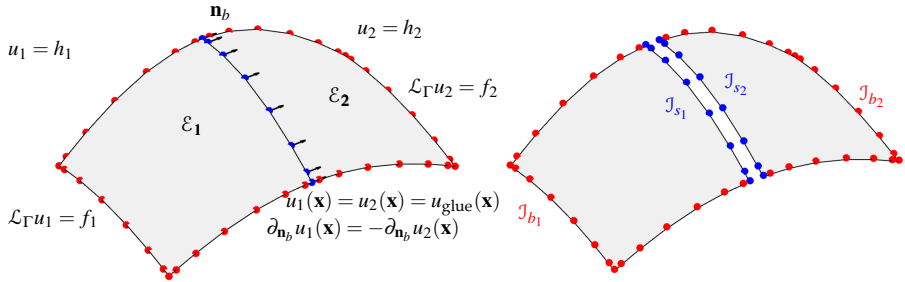
\begin{figure}[htb]
    \centering
    \begin{overpic}[width=0.43\textwidth,clip,trim=0cm 0cm 0.2cm 0cm]{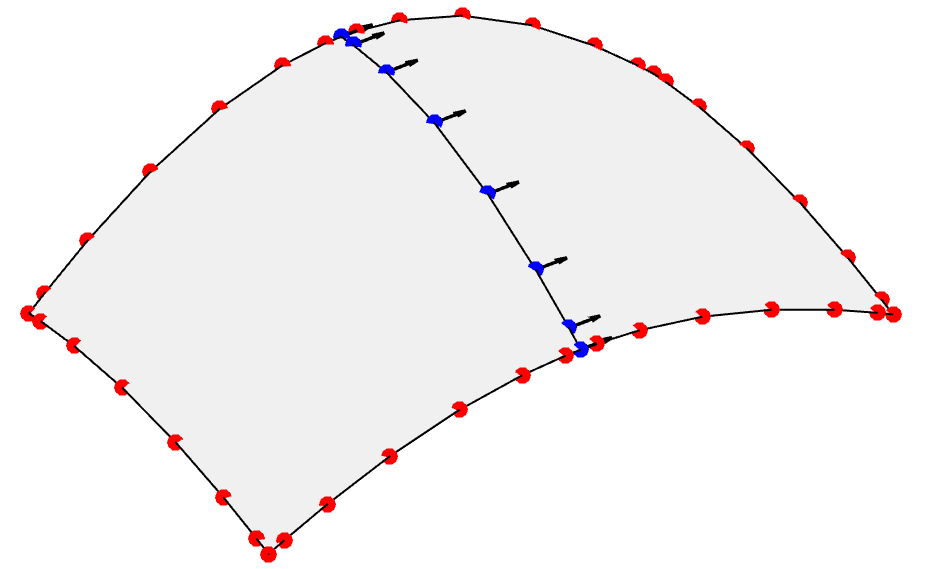}
            \put(32,36){\footnotesize$\mathcal{E}_{\bm{1}}$}
        \put(67,43){\footnotesize$\mathcal{E}_{\bm{2}}$}
        \put(-9,54){\footnotesize$u_1 = h_1$}
        \put(74,59){\footnotesize$u_2 = h_2$}
        \put(-9,14){\footnotesize$\mathcal{L}_{\Gamma} u_1 = f_1$}
        \put(86,45){\footnotesize$\mathcal{L}_{\Gamma}u_2 = f_2$}
        \put(39,63.3){\footnotesize$\mathbf{n}_{b}$}
        \put(57,18){\footnotesize$u_1(\mathbf{x}) = u_2(\mathbf{x})=u_{\text{glue}}(\mathbf{x})$}
        \put(52.4,11.5){\footnotesize$\partial_{\mathbf{n}_b} u_1(\mathbf{x}) =- \partial_{\mathbf{n}_b} u_2(\mathbf{x})$}
    \end{overpic}
    \hspace{1em}
    \begin{overpic}[width=0.405\textwidth,clip,trim=0cm 0cm 0.2cm 0cm]{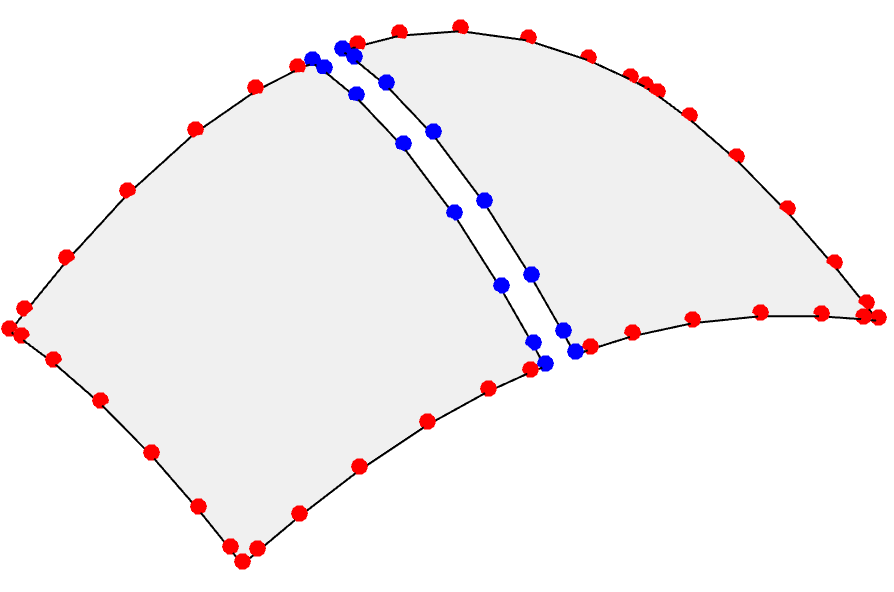}
        \put(3,16) {\footnotesize\color{red}$\mathcal{I}_{b_1}$}
        \put(89,46){\footnotesize\color{red}$\mathcal{I}_{b_2}$}
        \put(40,43){\footnotesize\color{blue}$\mathcal{I}_{s_1}$}
        \put(55,49){\footnotesize\color{blue}$\mathcal{I}_{s_2}$}
    \end{overpic}

    \caption{Interface coupling of two adjacent surface elements $\mathcal{E}_1$ and $\mathcal{E}_2$ via continuity of the solution and its binormal derivative.
Red points denote boundary collocation points $\mathcal{I}_{b_1}$ and $\mathcal{I}_{b_2}$, while blue points indicate aligned interface collocation points $\mathcal{I}_{s_1}$ and $\mathcal{I}_{s_2}$.}

    \label{fig:two_glued_patches}
\end{figure}

Assume now that $\Gamma$ is partitioned into adjacent elements. For simplicity, we consider the case of two elements, denoted by $\mathcal{E}_1$ and $\mathcal{E}_2$, as illustrated in Figure~\ref{fig:two_glued_patches}. We denote their common interface by~$
\mathcal{I}_{12} = \partial \mathcal{E}_1 \cap \partial \mathcal{E}_2$ .

\begin{equation}\label{eq:glue_1}
        \mathcal{L}_{\Gamma} u_1 (\mathbf{x}) = f_1 (\mathbf{x}), \quad \mathbf{x} \in \mathcal{E}_1,\quad \text{with}\quad     u_1 (\mathbf{x}) = h_1 (\mathbf{x}), \quad \mathbf{x} \in \partial \mathcal{E}_1 \setminus \mathcal{I}_{12},
\end{equation}
\begin{equation}\label{eq:glue_2}
   \mathcal{L}_{\Gamma} u_2 (\mathbf{x}) = f_2 (\mathbf{x}), \quad \mathbf{x} \in \mathcal{E}_2, \quad \text{with}\quad     u_2 (\mathbf{x}) = h_2 (\mathbf{x}), \quad \mathbf{x} \in \partial \mathcal{E}_2 \setminus \mathcal{I}_{12},
 \end{equation}
with continuity conditions:
\begin{align*}
    u_1 (\mathbf{x}) &= u_2 (\mathbf{x}), \quad \mathbf{x} \in \mathcal{I}_{12}, \nonumber\\
    \partial_{\mathbf{n}_b} u_1(\mathbf{x}) &=- \partial_{\mathbf{n}_b} u_2(\mathbf{x}), \quad \mathbf{x} \in \mathcal{I}_{12}.\nonumber
\end{align*}

 The patching problem \eqref{eq:glue_1} and \eqref{eq:glue_2} can be regarded as two decoupled, four-sided Dirichlet problems when given a suitable piece of Dirichlet data along  $\mathcal{I}_{12}$. That is, there exists an interface function $u_{\text{glue}}$ such that \eqref{eq:glue_1} and \eqref{eq:glue_2}  is equivalent to

\begin{equation}\label{eq:glue_1_1}
\begin{aligned}
\mathcal{L}_{\Gamma} u_1(\mathbf{x}) &= f_1(\mathbf{x}), && \mathbf{x}\in\mathcal{E}_1,\\
u_1(\mathbf{x}) &= h_1(\mathbf{x}), && \mathbf{x}\in \partial\mathcal{E}_1\setminus \mathcal{I}_{12},\\
u_1(\mathbf{x}) &= u_{\text{glue}}(\mathbf{x}), && \mathbf{x}\in\mathcal{I}_{12}.
\end{aligned}
\end{equation}

\begin{equation}\label{eq:glue_2_1}
\begin{aligned}
\mathcal{L}_{\Gamma} u_2(\mathbf{x}) &= f_2(\mathbf{x}), && \mathbf{x}\in\mathcal{E}_2,\\
u_2(\mathbf{x}) &= h_2(\mathbf{x}), && \mathbf{x}\in \partial\mathcal{E}_2\setminus \mathcal{I}_{12},\\
u_2(\mathbf{x}) &= u_{\text{glue}}(\mathbf{x}), && \mathbf{x}\in\mathcal{I}_{12}.
\end{aligned}
\end{equation}

To determine the unknown interface function \( u_{\text{glue}} \), we construct a solver referred to as the interface solution operator, denoted by \( \mathbf{S}_{\text{glue}} \), such that
 $u_{\text{glue}}(\mathbf{x}) = \mathbf{S}_{\text{glue}} \begin{bmatrix}
h_1(\mathbf{x}) \\
h_2(\mathbf{x})
\end{bmatrix}$. This solver is constructed using local operators from each element. Specifically, we first build local solvers for the elements $\mathcal{E}_1$ and $\mathcal{E}_2$, and then combine parts of these operators to create the interface solution operator $\mathbf{S}_{\text{glue}}$. 

Once the interface function $u_{\text{glue}}$ is determined, the two subproblems from equations \eqref{eq:glue_1_1} and \eqref{eq:glue_2_1} become independent and can be solved separately by applying spectral collocation on $\mathcal{E}_1$ and $\mathcal{E}_2$.

Let $\mathcal{I}_{s_1}, \mathcal{I}_{s_2} \subset \{1, \dots, 2dn^{d-1} \}$ be index sets corresponding to the shared interface points (i.e., the points which are interior to the merged domain $\mathcal{E}_\text{glue}:=\mathcal{E}_1 \cup \mathcal{E}_2$) with respect to elements $\mathcal{E}_1$ and $\mathcal{E}_2$, and $\mathcal{I}_1$ and $\mathcal{I}_2$ the remaining indices corresponding to the boundary points (see Figure~\ref{fig:two_glued_patches}, right panel). As the points of $\mathcal{E}_1$ and $\mathcal{E}_2$ are identical along the shared interface $\mathcal{I}_{12}$, continuity of the solution across $\mathcal{I}_{12}$ simply means that\footnote{The superscripts \( s_1 \), \( s_2 \) denote restriction to interface indices} $\mathbf{u}^{s_1}_1 = \mathbf{u}^{s_2}_2$, so let us denote these solution values by $\mathbf{u}_{\text{glue}}$.  

Continuity conditions across element interfaces are imposed locally using discrete Dirichlet-to-Neumann (DtN) operators.
\[
\mathrm{DtN}_{\mathcal{E}_k} = \mathcal{D}_{\mathcal{E}_k} \mathbf{S}_{\mathcal{E}_k}, \quad k=1,2.
\]

Since $\mathrm{DtN}_{\mathcal{E}_k}$ maps from the quadrilateral boundary to itself, it is a square matrix of size $2dn^{d-1} \times 2dn^{d-1}$. We define the boundary data vectors \( \mathbf{h}_1^b, \mathbf{h}_2^b \in \mathbb{R}^{2d n^{d-1}} \) for the two neighboring elements as
\[
\mathbf{h}_1^b(\mathbf{x}) :=
\begin{bmatrix}
\mathbf{h}_1(\mathbf{x}) \\
\mathbf{u}_{\text{glue}}(\mathbf{x})
\end{bmatrix}, \quad
\mathbf{h}_2^b(\mathbf{x}) :=
\begin{bmatrix}
\mathbf{h}_2(\mathbf{x}) \\
\mathbf{u}_{\text{glue}}(\mathbf{x})
\end{bmatrix},
\]
where \( \mathbf{h}_1(\mathbf{x}) \) and \( \mathbf{h}_2(\mathbf{x}) \) correspond to the non-shared boundary degrees of freedom, and \( \mathbf{u}_{\text{glue}}(\mathbf{x}) \) represents the shared interface values between the two elements.

The corresponding homogeneous solutions on each element are then given by
\[
\mathbf{w}_k = \mathbf{S}_{\mathcal{E}_k} \mathbf{h}_k^b, \quad k = 1,2.
\]

The continuity condition for the binormal derivative across the shared interface requires that

\begin{equation}\label{eq:interface_sol}
\left(\mathrm{DtN}_{\mathcal{E}_1} \mathbf{h}_1^b + \mathcal{D}_{\mathcal{E}_1}\bm{v}_1\right)^{s_1} 
+ \left(\mathrm{DtN}_{\mathcal{E}_2}\mathbf{h}_2^b + \mathcal{D}_{\mathcal{E}_2}\bm{v}_2\right)^{s_2} = 0,
\end{equation}

where \( \bm{v}'_k := \mathcal{D}_{\mathcal{E}_k}\bm{v}_k \) denotes particular fluxes associated with particular solutions on each element. Rewriting the system \eqref{eq:interface_sol} in terms of the interface unknown \( \mathbf{u}_{\text{glue}} \), we obtain separate linear systems for the homogeneous and particular interface unknowns:


\begin{align}
\left( \mathrm{DtN}_{\mathcal{E}_1}^{s_1 s_1} + \mathrm{DtN}_{\mathcal{E}_2}^{s_2 s_2} \right) \mathbf{w}_{\text{glue}} &= \mathrm{DtN}_{\mathcal{E}_1}^{s_1 b_1} \mathbf{h}_1 + \mathrm{DtN}_{\mathcal{E}_2}^{s_2 b_2} \mathbf{h}_2, \label{eq:interface_homo}\\
\left( \mathrm{DtN}_{\mathcal{E}_1}^{s_1 s_1} + \mathrm{DtN}_{\mathcal{E}_2}^{s_2 s_2} \right) \bm{v}_{\text{glue}} &= \bm{v}'^{s_1}_1 + \bm{v}'^{s_2}_2 \label{eq:interface_partic}.
\end{align}
with local operators constructed on each element.
By solving Eq.~\eqref{eq:interface_homo}, we obtain the interface solution operator \( \mathbf{S}_{\text{glue}} \), which is defined as follows:
\begin{equation}
\mathbf{S}_{\text{glue}} := -\left( \mathrm{DtN}_{\mathcal{E}_1}^{s_1 s_1} + \mathrm{DtN}_{\mathcal{E}_2}^{s_2 s_2} \right)^{-1} \left[\mathrm{DtN}_{\mathcal{E}_1}^{s_1 b_1} \quad\mathrm{DtN}_{\mathcal{E}_2}^{s_2 b_2}\right]\,.
\end{equation}
The Schur complement also allows us to write down the Dirichlet-to-Neumann operator for the merged domain. Using the new interface solution operator $\mathbf{S}_{\text{glue}}$, we construct the new Dirichlet-to-Neumann operator

\begin{equation}
\mathrm{DtN}_{\text{glue}} := 
\begin{bmatrix}
    \mathrm{DtN}_{\mathcal{E}_1}^{b_1 b_1} & 0 \\
    0 & \mathrm{DtN}_{\mathcal{E}_2}^{b_2 b_2}
\end{bmatrix}
+
\begin{bmatrix}
    \mathrm{DtN}_{\mathcal{E}_1}^{b_1 s_1} \\
    \mathrm{DtN}_{\mathcal{E}_2}^{b_2 s_2}
\end{bmatrix}
\mathbf{S}_{\text{glue}}.
\end{equation}
In a similar manner, the specific flux for the merged domain is given by
\begin{equation}
\bm{v}'_{\text{glue}} =
\begin{bmatrix}
    \bm{v}'^{b_1}_1 \\
    \bm{v}'^{b_2}_2
\end{bmatrix}
+
\begin{bmatrix}
    \mathrm{DtN}_{\mathcal{E}_1}^{b_1 s_1} \\
    \mathrm{DtN}_{\mathcal{E}_2}^{b_2 s_2}
\end{bmatrix}
\bm{v}_{\text{glue}}.
\end{equation}

\subsection{The hierarchical scheme}\label{se:hps_scheme}
The term hierarchical refers to the organization of the computational mesh into a hierarchical binary tree structure \cite{martinsson2019fast}. Each element of the mesh, before any merging or gluing occurs, corresponds to a leaf box in this tree.
At each stage of the merging process, two neighboring leafs are combined to form a larger domain. The result of merging two such elements is a pair of operators defined on the merged domain:
 \( \mathcal{E}_{\text{glue}}  \): (1) a solution operator, \( \mathbf{S}_{\text{glue}} \), that solves for the unknown interface values \( \mathbf{u}_{\text{glue}} \) inside \( \mathcal{E}_{\text{glue}}  \); and (2) a Dirichlet-to-Neumann (DtN) operator, \( \mathrm{DtN}_{\text{glue}} \), that maps boundary data to outward fluxes on \( \mathcal{E}_{\text{glue}}  \). These operators encapsulate all necessary information to solve the PDE on \( \mathcal{E}_{\text{glue}}  \). After merging, the new parent element\( \mathcal{E}_{\text{glue}}  \)is functionally equivalent to the original elements \( \mathcal{E}_1 \) or \( \mathcal{E}_2 \), and thus can be treated as a single element, ready to be merged again with additional domains. 
\begin{figure}[htb]
    \centering

    \begin{subfigure}[b]{0.45\textwidth}
        \centering
        \begin{overpic}[width=\textwidth]{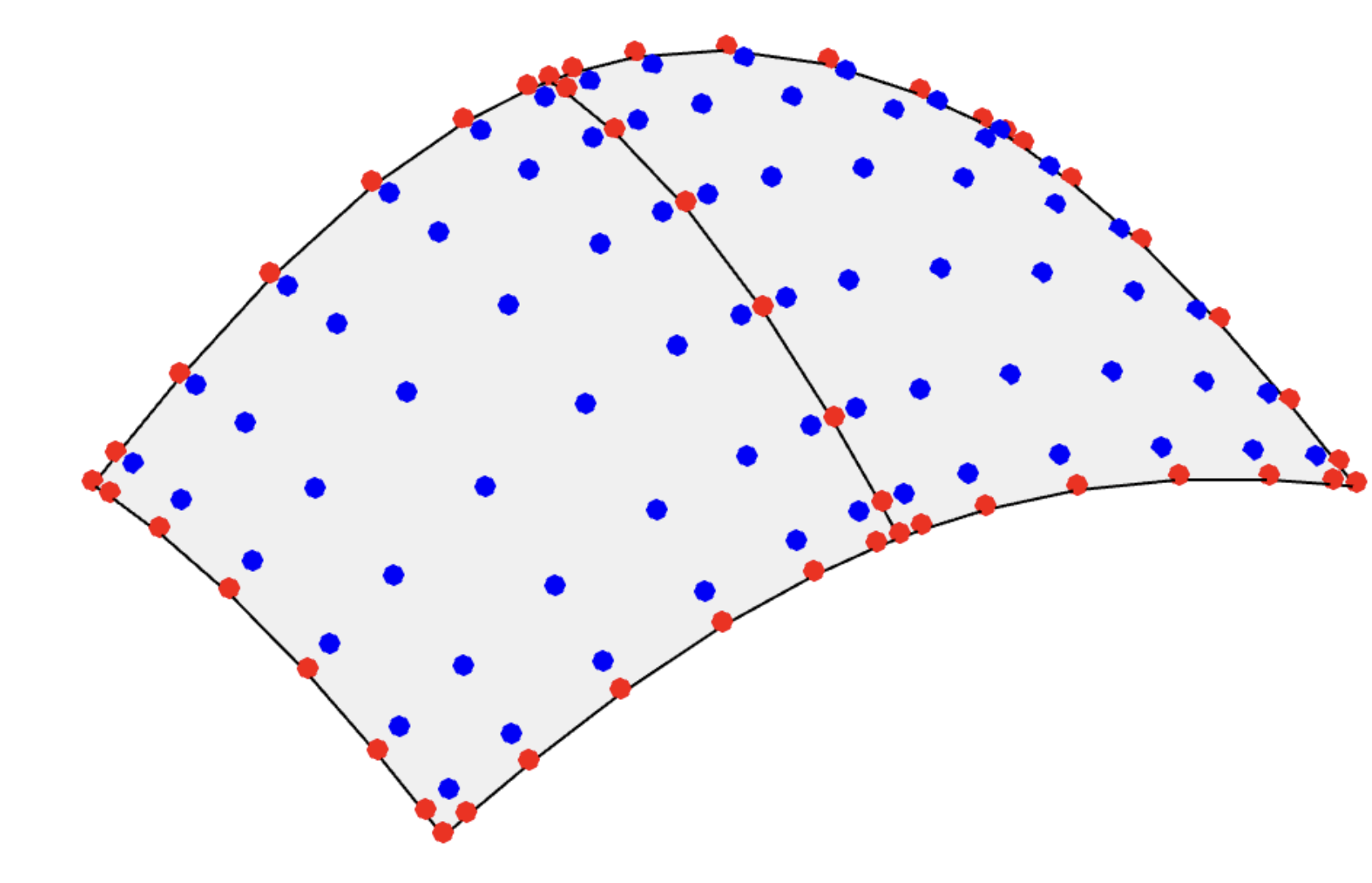}
            \put(32,36){\footnotesize$\mathcal{E}_{1}$}
            \put(70,43){\footnotesize$\mathcal{E}_{2}$}
                      \put(95,56){\footnotesize$\textbf{Step 1}$}
        \end{overpic}
    \end{subfigure}
    \hspace{1em}
    \begin{subfigure}[b]{0.45\textwidth}
        \centering
        \begin{overpic}[width=\textwidth]{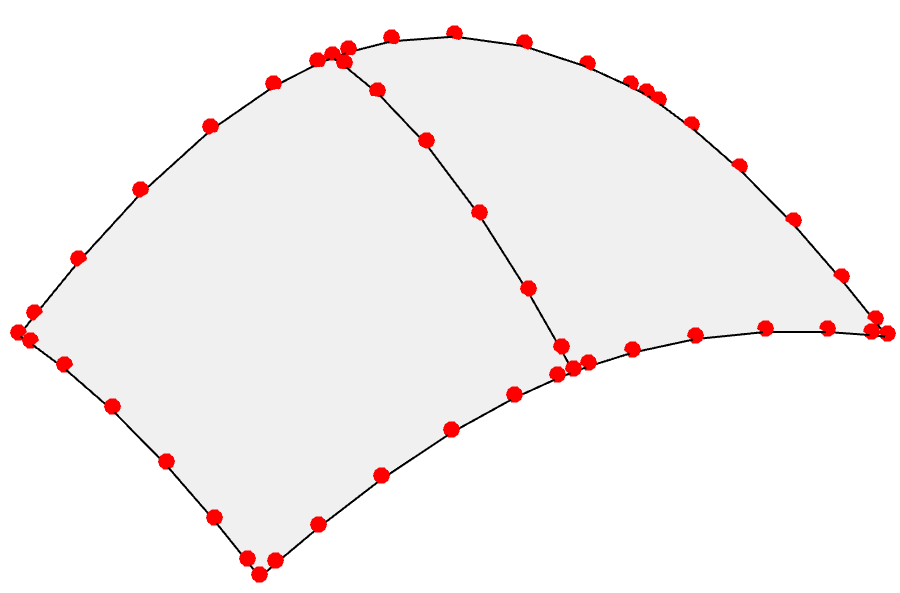}
            \put(32,36){\footnotesize$\mathcal{E}_{1}$}
            \put(15,28){\footnotesize$\textcolor{blue}{u_{1}(\mathbf{x})} = \mathbf{S}_{\mathcal{E}_{1}} \begin{bmatrix}
            \textcolor{red}{h_1(\mathbf{x})} \\
            \textcolor{red}{u_{\text{glue}}(\mathbf{x})}
            \end{bmatrix}$}
            \put(40,7)
            {\scriptsize$\textcolor{red}{\frac{\partial u_1}{\partial n_b}}=\text{DtN}_{\mathcal{E}_{\bm{1}}}  
\begin{bmatrix}
\textcolor{red}{h_1(\mathbf{x}}) \\
\textcolor{red}{u_{\text{glue}}(\mathbf{x})}
\end{bmatrix}$}
\put(75,18)
            {\scriptsize$\textcolor{red}{\frac{\partial u_2}{\partial n_b}}=\text{DtN}_{\mathcal{E}_{\bm{1}}}  
\begin{bmatrix}
\textcolor{red}{h_2(\mathbf{x}}) \\
\textcolor{red}{u_{\text{glue}}(\mathbf{x})}
\end{bmatrix}$}
            \put(67,43){\footnotesize$\mathcal{E}_{2}$}
            \put(67,62){\footnotesize$\textcolor{blue}{u_{2}(\mathbf{x})} = \mathbf{S}_{\mathcal{E}_{2}} \begin{bmatrix}
            \textcolor{red}{h_2(\mathbf{x})} \\
            \textcolor{red}{u_{\text{glue}}(\mathbf{x})}
            \end{bmatrix}$}
        \end{overpic}
    \end{subfigure}

    \vspace{1em}

    \begin{subfigure}[b]{0.45\textwidth}
        \centering
        \begin{overpic}[width=\textwidth]{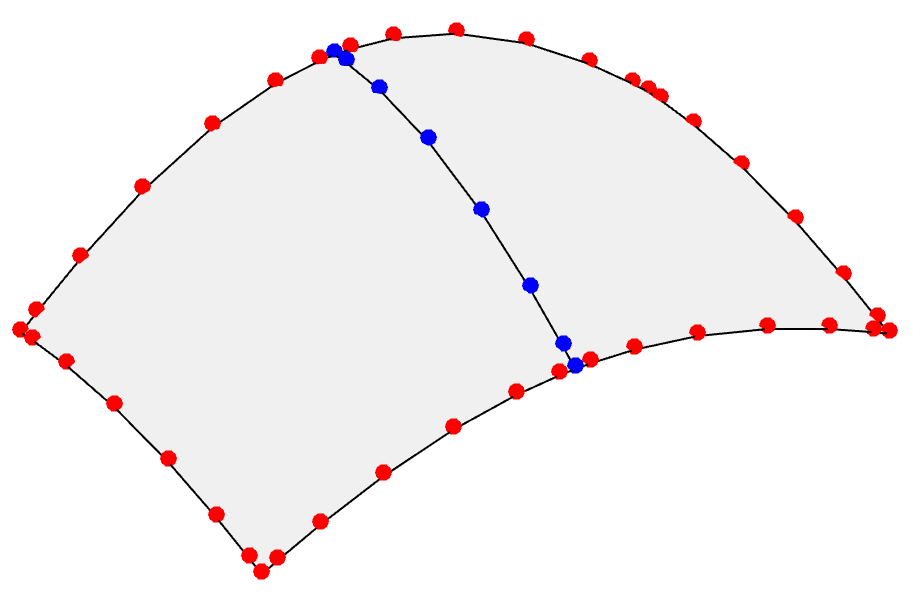}
            \put(32,36){\footnotesize$\mathcal{E}_{1}$}
            \put(55,13){\footnotesize$\textcolor{blue}{u_{\text{glue}}(\mathbf{x})} = \mathbf{S}_{\text{glue}} \begin{bmatrix}
            \textcolor{red}{h_1(\mathbf{x})} \\
            \textcolor{red}{h_{2}(\mathbf{x})}
            \end{bmatrix}$}
            \put(67,43){\footnotesize$\mathcal{E}_{2}$}
                        \put(95,56){\footnotesize$\textbf{Step 2}$}
        \end{overpic}
    \end{subfigure}
    \hspace{1em}
    \begin{subfigure}[b]{0.45\textwidth}
        \centering
        \begin{overpic}[width=\textwidth]{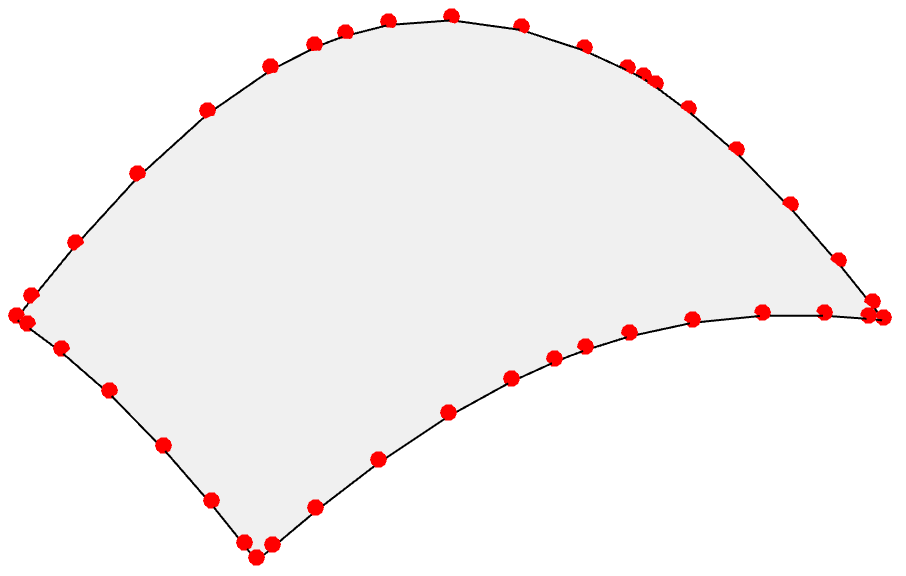}
            \put(32,36){\footnotesize$\mathcal{E}_{\text{glue}}:=\mathcal{E}_1 \cup \mathcal{E}_2$}
        \end{overpic}
    \end{subfigure}

    \caption{
   An overview of the hierarchical merge process.
    \textbf{Step 1:} For each patch \( \mathcal{E}_i \), remove the internal points (blue) and retain only the boundary points (red). Use the remaining data to construct the solution operator \( \mathbf{S}_{\mathcal{E}_i} \) and the Dirichlet-to-Neumann operator \( \mathrm{DtN}_{\mathcal{E}_i},\; i=1,2 \).
    \textbf{Step 2:} Merge patches in pairs as described in Subsection~\ref{merge_section}. Form the equilibrium equation using \( \mathrm{DtN}_{\mathcal{E}_1} \) and \( \mathrm{DtN}_{\mathcal{E}_2} \) and eliminate the interior points of the new larger patches. Construct the new solution operator \( \mathbf{S}_{\text{glue}} \) and corresponding DtN operator.
    }
    \label{fig:gluing_steps}
\end{figure}
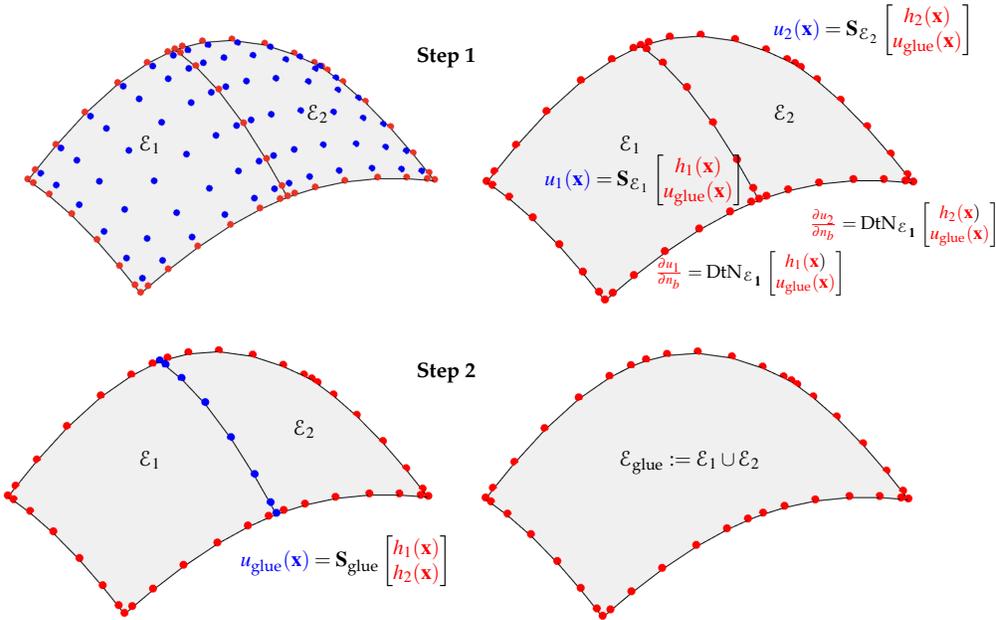
This hierarchical merging procedure, enabled by local operators, constitutes the
hierarchical Poincaré--Steklov scheme. For a detailed description of the
algorithmic construction and implementation, we refer the reader to
\cite{fortunato2022highorder}. An illustrated overview of the merge process is
described in Figure~\ref{fig:gluing_steps}.

\section{Triangular surface discretizations}\label{sec:extend_rohmbus}
The HPS domain decomposition framework described in the previous sections is
formulated under the assumption that the surface
\( \Gamma \subset \mathbb{R}^{d+1} \) is discretized by a high-order quadrilateral
mesh. This setting allows the use of tensor-product spectral discretizations and
is particularly well suited to structured surface representations. In many practical applications, however, the surface \( \Gamma \) is given as a triangulated mesh \(   \bigcup_{\widehat{T} \in \mathcal{\widehat{T}}_h} \widehat{T}\). Simplicial domains, such as triangles and tetrahedra, which are not tensor product domains offer far greater flexibility than Cartesian products of intervals (e.g., squares or cubes) when it comes to handling complex geometries through partitioning methods \cite{ciarlet2002finite}. The procedure described in \ref{merge_section} and \ref{se:hps_scheme} remains unchanged, provided that suitable high-order points and basis functions are constructed on the reference triangle, allowing the PDE to be discretized with high-order accuracy on each element, similarly to the quadrilateral case.  We propose two complementary
approaches for extending the HPS framework to triangular surface meshes.  
The first approach is a rhombus-based remeshing technique, which quadrilateralizes the mesh by converting triangular elements into quadrilaterals while minimizing the number of additional elements introduced by the centroid-based subdivision scheme~\cite{fortunato2021ultraspherical}. The second approach is a direct high-order discretization on triangles, relying on the use of Dubiner polynomials, which enable spectral collocation on triangular elements. Details of both approaches are presented in the following subsections.

\subsection{Rhombus-based remeshing}
\begin{figure*}[!t]
  \centering

  \begin{subfigure}{0.45\textwidth}
    \centering
    \includegraphics[width=0.7\linewidth]{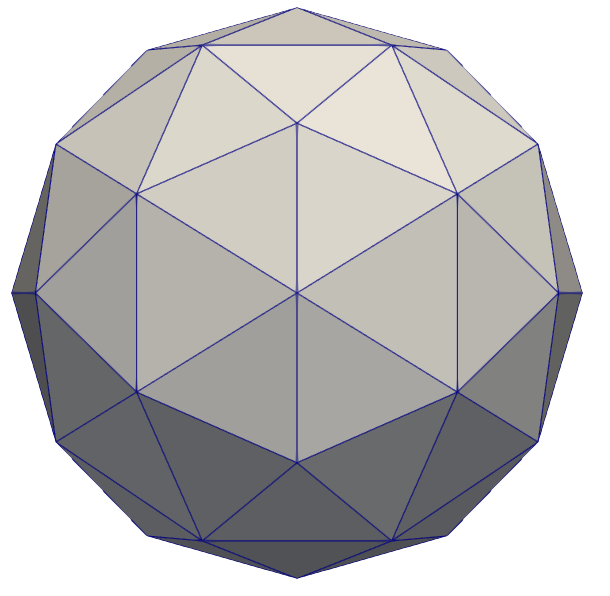}
    \label{fig:tri_quad_a}
  \end{subfigure}
  \hfill
  \begin{subfigure}{0.45\textwidth}
    \centering
    \includegraphics[width=0.7\linewidth]{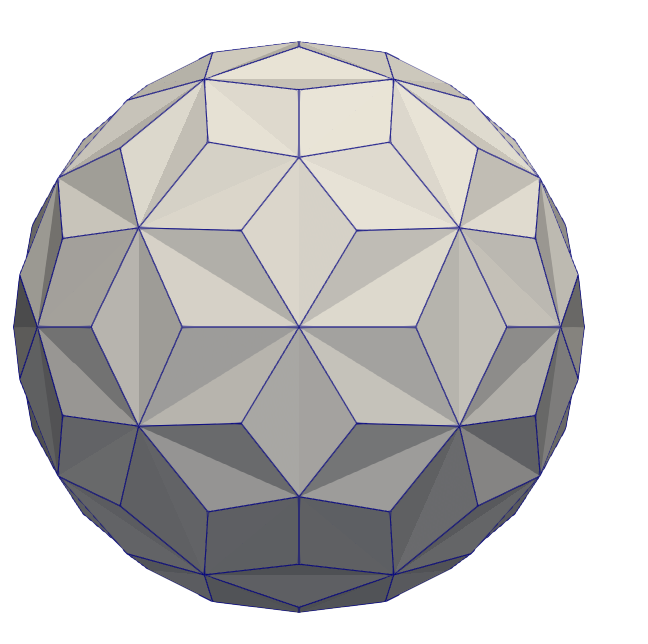}
    \label{fig:tri_splitting3}
  \end{subfigure}

 \caption{Reference triangular mesh $\Gamma_h^{\mathrm{tri}}$ (80 triangles) and the corresponding quadrilateral mesh generated using the proposed rhombus-based remeshing method (120 quads).}

  \label{fig:quad_methods_2x2}
\end{figure*}

The method is straightforward: we replace the original triangulated mesh with a conforming collection of quadrilateral patches 
$\{\widehat{\mathcal{E}}_k\}_{k=1}^K$, each equipped with a smooth mapping \( \psi_k : \Omega_d \rightarrow \widehat{\mathcal{E}}_k \), and a point set \( \{ \psi_k(p_{\boldsymbol{\alpha}}) \}\) with \({\boldsymbol{\alpha} \in A_{d,n}} \) for use in high-order interpolation.
Let \( \{ \widehat{T}_k \}_{k=1}^{N} \) denote the set of triangles in \( \Gamma_h^{\text{tri}} \). For any pair of adjacent triangles \( \widehat{T}_i, \widehat{T}_j \in \widehat{\mathcal{T}}_{h} \) sharing an edge \( e_{ij} \), let \( c_i \) and \( c_j \) denote the centroids of the triangles, and let \( \{ v_{a}, v_{b} \} \in e_{ij} \) be the vertices of the shared edge. A quadrilateral patch is then formed as
\[
\widehat{\mathcal{E}}_k = \text{Quad}(v_a, c_j, v_b, c_i).
\]
Now, the elements $\widehat{\mathcal{E}}_k$ (see Figure~\ref{fig:quad_methods_2x2}) can be naturally parameterized using cubical maps \( \psi_k : \Omega_d \rightarrow \widehat{\mathcal{E}}_k \), and a corresponding set of mapped points \( \{ \psi_k(p_{\boldsymbol{\alpha}}) \}_{\boldsymbol{\alpha} \in A_{d,n}} \subset \widehat{\mathcal{E}}_k \cap \Gamma \). Then, the resulting mesh \( \Gamma_h^{\text{quad}} \) can be interpolated to obtain  an $n^{\text{th}}$-order cubical parametrization in the sense of Definition~\ref{def:46}, thereby enabling the direct application of the HPS scheme introduced in the preceding sections. Compared to other quadrilateralization strategies such as centroid-based subdivision~\cite{fortunato2021ultraspherical}, the proposed construction forming rhombus-shaped quadrilaterals from the triangular mesh yields the smallest number of quadrilateral elements.

\begin{figure}[htbp]
  \centering
  
  \begin{subfigure}[b]{0.48\textwidth}
    \includegraphics[width=\linewidth]{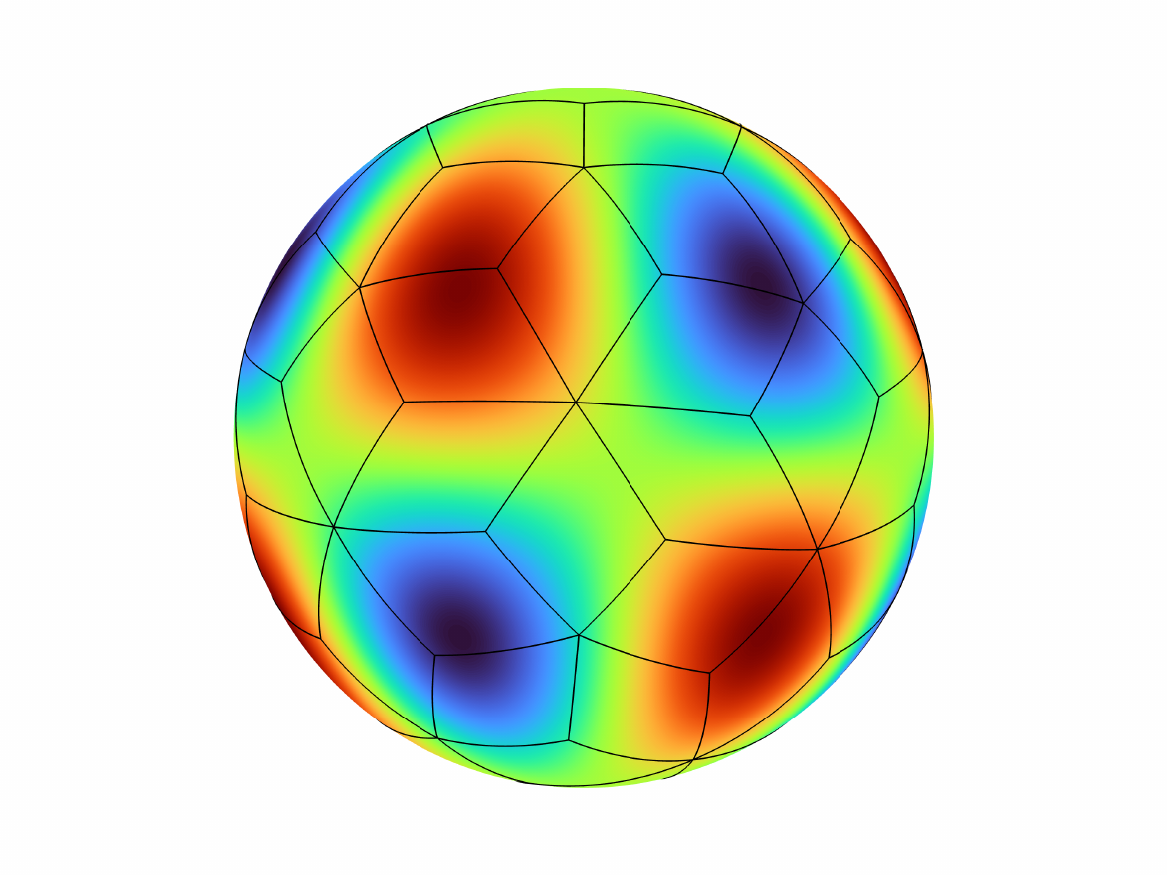}
    \caption{Computed solution}
    \label{fig:sphere_velocity_newn}
  \end{subfigure}
  \hfill
  \begin{subfigure}[b]{0.48\textwidth}
    \includegraphics[width=\textwidth]{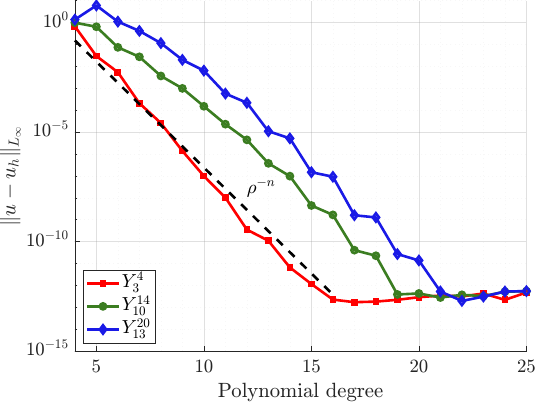}
    \caption{$L_\infty$-error vs.\ degree} 
    \label{fig:sphere_error_l2}
  \end{subfigure}

  \caption{The surface mesh of the sphere geometry and the solution of~\eqref{eq:exp22}.  
  (Right) $L_\infty$-error versus polynomial degree for three harmonic test functions, compared with an exponential fit rate $\rho^{-n}$, where $\rho \approx 9.3$.}
  \label{fig:sphere_velocity_new}
\end{figure}
\paragraph{HPS on the sphere}
To illustrate the remeshing procedure, we consider a model problem on the unit sphere.  
Specifically, we study the Laplace--Beltrami equation
\begin{equation}\label{eq:exp21}
- \Delta_\Gamma u + u = f,
\qquad \text{on } \Gamma = \{ \bm{x} \in \mathbb{R}^3 : \|\bm{x}\|_2 = 1 \},
\end{equation}
where the exact solution is chosen to be the spherical harmonic 
$u(\bm{x}) = Y_\ell^m(\bm{x})$ with $(\ell,m) = (4,3)$.
The corresponding right-hand side is therefore
\[
f(\bm{x}) = -\ell(\ell+1) Y_\ell^m(\bm{x}) + Y_\ell^m(\bm{x}),
\]
so that $u$ satisfies~\eqref{eq:exp21} exactly.

The initial triangular mesh is transformed into rhombus-shaped quadrilateral elements, after which the vertices are projected onto the sphere using the closest-point map~$\pi(\bm{x}) := \frac{\bm{x}}{\|\bm{x}\|}.$
A high-order curved geometry is then obtained by interpolating $\pi$ piecewise over the flat quadrilaterals using $\Cheb_{d,n}$ grid.

Figure~\ref{fig:sphere_velocity_newn} shows the harmonic function evaluated on a rhombus-based sphere mesh.
The right plot in Figure~\ref{fig:sphere_velocity_new} displays the decay of the $L_\infty$-error as a function of polynomial degree, using three different spherical harmonics $Y_\ell^m$ as exact solutions. The results confirm exponential convergence consistent with spectral accuracy. The best-exponential fit rate $\rho^{-n}$, with $\rho \approx 9.3$, is included to indicate the expected decay rate. 

\subsection{Dubiner polynomials on the reference triangle}

As in quadrilateral based spectral element methods, triangular spectral element
discretizations require an orthogonal polynomial basis together with a suitable
set of interpolation nodes. On quadrilateral reference elements, such as the square
\( [-1,1]^2 \), these ingredients arise naturally from tensor-product constructions
based on Chebyshev polynomials and the associated Chebyshev points.

For triangular elements, however, this tensor-product structure is no longer
available. We therefore employ Dubiner polynomials, which form an
orthogonal basis on the reference triangle $\Delta_2 = \{(\xi,\eta)\in\mathbb{R}^2 : \xi \ge 0,\ \eta \ge 0,\ \xi+\eta \le 1\}.$ 
The Dubiner basis is obtained by applying the classical collapsed-coordinate mapping to
tensor-product Jacobi polynomials on the reference square and provides a natural
high-order spectral discretization on triangular domains.

More precisely, for total polynomial degree $n$, the Dubiner basis functions are
\[
\varphi_{ij}(\xi,\eta)
  := C_{ij}\, 2^j (1-\eta)^i\,
      J_i^{0,0}\!\left( \frac{2\xi}{1-\eta} - 1 \right)
      J_j^{2i+1,0}(2\eta - 1),
\qquad 0 \le i,\, j \le n,\; i+j \le n,
\]
where the normalization constant is
\[
C_{ij} := \sqrt{ \frac{2(2i+1)(i+j+1)}{4^i} }.
\]
Here $J_m^{\alpha,\beta}$ denotes the Jacobi polynomial orthogonal with 
respect to the weight
\[
w(x) = (1-x)^\alpha (1+x)^\beta,
\]
that is,
\[
\int_{-1}^1 (1-x)^\alpha (1+x)^\beta\,
  J_m^{\alpha,\beta}(x)\, J_q^{\alpha,\beta}(x)\, dx
  = \frac{2}{2m+1}\, \delta_{mq}.
\]

Let
\[
\Pi_{2,n}(\Delta_2)= \operatorname{span}\{\varphi_{ij}\}_{i+j\le n},
\qquad 
N_n = \frac{(n+1)(n+2)}{2},
\]
and let
\[
\widehat{X}_n = \{\, (\xi_i,\eta_j) \mid 0\le i\le n,\; 0\le j\le n-i \,\}
\]
denote a set of spectral nodes on the reference simplex $\Delta_2$, generated by
the recursive, parameter-free construction~\cite{isaac2020recursive}.
In our implementation, the underlying one-dimensional seed grid consists of
second-kind Chebyshev (Chebyshev--Lobatto) points, so that the induced triangular
nodes coincide with Chebyshev boundary nodes on each edge. For notational convenience, we enumerate the nodes $(\xi_i,\eta_j)$ by a single
index $m=0,\dots,N_n$, writing $(\xi_m,\eta_m)$, and likewise enumerate the basis
functions $\varphi_{ij}$ by $k=0,\dots,N_n$, writing $\varphi_k$.

For any function $u$ defined on $\Delta_2$, we denote by $Q_{\Delta_2,n} u \in \Pi_{2,n}(\Delta_2)$
its polynomial interpolant satisfying
\[
Q_{\Delta_2,n}u(\xi_m,\eta_m) = u(\xi_m,\eta_m), \qquad m = 0,\dots,N_n.
\]

We expand $Q_{\Delta_2,n}u$ in the Dubiner basis by enumerating
$\varphi_{ij}$ with a single index $k=1,\dots,N_n$:
\[
Q_{\Delta_2,n}u(\xi,\eta) = \sum_{k=0}^{N_n} c_k \, \varphi_k(\xi,\eta),
\]
where $c_k$ are the modal coefficients.
Evaluating this representation at the collocation nodes yields
\[
u_m = u(\xi_m,\eta_m)
    = \sum_{k=0}^{N_n} c_k \, \varphi_k(\xi_m,\eta_m).
\]

Introducing the Dubiner Vandermonde matrix
\[
K_{mk} = \varphi_k(\xi_m,\eta_m),
\]
this relation is written compactly as
\[
\mathbf{u} = K \boldsymbol{c}.
\]
Since the Dubiner nodes are good interpolation points, the matrix $K$ is well conditioned,
and hence the modal coefficients follow from
\[
\boldsymbol{c} = K^{-1}\mathbf{u}.
\]

To compute derivatives of $Q_{\Delta_2,n} u$ at the collocation points, we differentiate
the Dubiner expansion:
\[
\partial_\xi Q_{\Delta_2,n} u(\xi_m,\eta_m)
  = \sum_{k=0}^{N_n} c_k\, \partial_\xi \varphi_k(\xi_m,\eta_m),
\qquad
\partial_\eta Q_{\Delta_2,n} u(\xi_m,\eta_m)
  = \sum_{k=0}^{N_n} c_k\, \partial_\eta \varphi_k(\xi_m,\eta_m).
\]
Define the derivative Vandermonde matrices
\[
(K_\xi)_{mk} = \partial_\xi \varphi_k(\xi_m,\eta_m),
\qquad
(K_\eta)_{mk} = \partial_\eta \varphi_k(\xi_m,\eta_m).
\]
Using $\boldsymbol{c}=K^{-1}\mathbf{u}$, we obtain the nodal derivative vectors
\[
\mathbf{u}_\xi = K_\xi K^{-1}\mathbf{u},
\qquad
\mathbf{u}_\eta = K_\eta K^{-1}\mathbf{u}.
\]
Thus, the differentiation matrices on the reference triangle are
\begin{equation}\label{eq:diff_tri}
D_\xi = K_\xi K^{-1}, \qquad
D_\eta = K_\eta K^{-1},
\end{equation}
which map nodal values to directional derivatives in reference coordinates:
\[
\mathbf{u}_\xi = D_\xi \mathbf{u}, \qquad
\mathbf{u}_\eta = D_\eta \mathbf{u}.
\]

Here we remark that the differentiation matrices on triangular elements
have dimension $(N_n+1)\times (N_n+1)$, whereas in quadrilateral case they had dimension $(n+1)^2\times(n+1)^2$, owing to the Cartesian
tensor-product structure of the underlying mesh.

Alternatively, one may adopt a nodal Lagrange basis on each reference triangle, yielding a physical–space representation of the solution, rather than using the Dubiner polynomials as in~\cite{eskilsson2004triangular}, which leads to a modal formulation.

We now seek an explicit formula for the Lagrange basis by representing
them in terms of the reference basis, i.e.
\[
\ell_m(\xi,\eta)
    = \sum_{k=0}^{N_n} A_{mk}\, \varphi_k(\xi,\eta).
\]
We then use the cardinal property of the Lagrange polynomials
\[
\delta_{mj}
    = \sum_{k=0}^{N_n} A_{mk}\, \varphi_k(\xi_j,\eta_j),
\]
where $\delta$ is the Kronecker delta function, to determine that $A = K^{-T}$.
Since $K_{jk} = \varphi_k(\xi_j,\eta_j)$
is the generalized Vandermonde matrix, we construct
the Lagrange polynomials as follows:
\[
\ell_m(\xi,\eta)
    = \sum_{k=0}^{N_n} (K^{-1})_{km}\,
      \varphi_k(\xi,\eta),
\]
where $\ell_m$ are the multivariate Lagrange polynomials associated with
the node set $\widehat{X}_n$.

Building on this foundation, let us now consider a triangulation of the surface $\Gamma$, denoted by
\begin{equation*}
    \Gamma^{\text{tri}}_h = \bigcup_{\widehat{T} \in \mathcal{\widehat{T}}_h} \widehat{T},\quad |\mathcal{\widehat{T}}_h|=K,\;
\end{equation*}
which is not necessarily a high-order approximation of $\Gamma$.
In the following, we  define $n^{\text{th}}$-order simplex parametrization of the surface for each $\widehat{T} \in \mathcal{\widehat{T}}_h $ 
\begin{equation}\label{eq:re_para}
\phi_k : \Delta_2 \subset \mathbb{R}^2 \to \widehat{T}_k\subset\mathbb{R}^3,
\;
\phi_k = \pi_k \circ \tau_k,
\quad k=1,\dots,K 
\,.
\end{equation}

 \begin{definition}[Order-$n$ simplex parametrization]
Let
$\phi_k : \Delta_2 \subset \mathbb{R}^2 \to \widehat{T}_k\subset\mathbb{R}^3,
\;
\phi_k = \pi_k \circ \tau_k,
\; k=1,\dots,K,$ where $\pi_k$ is the closest point projection onto the surface $\Gamma$ 
and $\tau_k$ is affine.  
We say that the mesh $\Gamma^{\mathrm{tri}}_h$ is of order $n$ if each 
element is obtained by polynomial interpolation (in total degree) of$\{\phi_k(\xi_m, \eta_m)\}_{m=0}^{N_n},$ 
sampled at the triangular spectral nodes 
$\{(\xi_m,\eta_m)\}_{m=0}^{N_n}$.
\end{definition}

For every triangle \(\widehat{T} \in \mathcal{\widehat{T}}_{h}\), we compute  $\{\phi_k(\xi_m, \eta_m)\}_{m=0}^{N_n},$  and define an \(n\)-th order \(d\)-dimensional triangle \(T\) by applying polynomial interpolation of order \(n\) to the coordinates of the projected points $\{\phi_k(\xi_m, \eta_m)\}_{m=0}^{N_n},$ namely:
\begin{equation}
    Q_{\Delta_2,n}\phi_{k}(\xi,\eta)
    =
    \sum_{i=0}^{n}\sum_{j=0}^{n-i}
      \phi_{k}(\xi_{i}, \eta_{j})\, \varphi_{ij}(\xi,\eta),
    \qquad (\xi,\eta)\in \Delta_2,\quad k=1,\ldots, K.
\end{equation}\textbf{}

Mapping the piecewise flat surface results in a higher-order approximation $\Gamma^{\text{tri}}_{h,n} := Q_{\Delta_2,n}\phi_{k}(\Gamma^{\text{tri}}_h)$ of the surface $\Gamma$ that is contained in $\mathcal{N}_{\delta}$ if $h$ is small enough. Associated with the discrete surface $\Gamma^{\text{tri}}_{h,n}$ is a set of surface elements
\[
\mathcal{T}_{h,n} := \{ Q_{\Delta_2,n}\phi_{k}(\widehat{T}) \mid \widehat{T} \in \mathcal{T}_{h} \},
\]
such that
\begin{equation*}
    \Gamma^{\text{tri}}_{h,n} = \bigcup_{T \in \mathcal{T}_{h,n}} T\,.
\end{equation*}

Given $\Gamma^{\text{tri}}_{h,n}$, we can compute an $n^{\text{th}}$-order approximation of the volume element $g(\mathbf{x})$,

\[
g^n(\mathbf{x}) = \sqrt{\det \left(\mathbf{D}Q_{\Delta_2,n}\phi_{k}(\mathbf{x})^T \mathbf{D}Q_{\Delta_2,n}\phi_{k}(\mathbf{x})\right)}.
\]

\subsection{High-order spectral collocation on a triangular element}
\label{sub:local_op_tri}

We now describe the spectral collocation discretization of a surface
operator on a single triangular surface element $T_k$.  
Functions defined on $T_k$ are represented at the same interpolation nodes used for the geometry.

Let \(u\) be a function defined on \(T_k\), and define $u_{ij} := u\!\left(\phi_k(\xi_i,\eta_j)\right)$. The corresponding Dubiner interpolant is then
\begin{equation}\label{eq:u_interp_tri}
u(\xi,\eta)
  = \sum_{i=0}^{n}\sum_{j=0}^{\,n-i}
      u_{ij}\,\varphi_{ij}(\xi,\eta),
  \qquad (\xi,\eta)\in \Delta_2 .
\end{equation}

To discretize differential operators on \(T_k\), we follow the same strategy as in the quadrilateral case (see Section~\ref{sub:local_op}). We first construct the differentiation matrices on the reference triangle $\Delta_2$ and then map them to the physical element via the geometric transformation \(\phi_k\).

Consider the local parametrization $\phi_k(\xi,\eta)
  = \bigl(x_1(\xi,\eta),\, x_2(\xi,\eta),\, x_3(\xi,\eta)\bigr)$.
Using \eqref{eq:diff_tri}, the components of the surface gradient are given by
\begin{equation}\label{eq:surfgrad_tri}
D^{\Gamma}_{x_j}
  = M[\partial_{x_j}\xi]\, \mathcal{D}_\xi
    + M[\partial_{x_j}\eta]\, \mathcal{D}_\eta,
\qquad j = 1,2,3.
\end{equation}
For the second-order surface operator
\[
\mathcal{L}_\Gamma u
=
\sum_{i,j=1}^3 a_{ij}\,\partial_{x_i}^\Gamma\partial_{x_j}^\Gamma u
+
\sum_{i=1}^3 b_i\,\partial_{x_i}^\Gamma u
+
c\,u,
\]
the discrete operator on the element $\mathcal{E}_k$ is
\begin{equation}\label{eq:L_tri}
\mathbf{L}_{\mathcal{E}_k}
=
\sum_{i=1}^{3}\sum_{j=i}^{3}
  M[a_{ij}]\, D^\Gamma_{x_i}\, D^\Gamma_{x_j}
\;+\;
\sum_{i=1}^{3} M[b_i]\, D^\Gamma_{x_i}
\;+\;
M[c].
\end{equation}

As in the quadrilateral setting, it is convenient to reorder the degrees of
freedom. For the triangular element, we partition the Dubiner nodes
according to whether the Chebyshev points $(\xi_m,\eta_m)$ lie in the interior of
the reference triangle or on its edges. The interior index set has size $|I_{\operatorname{int}}| = \frac{(n-2)(n-3)}{2}$,
while the boundary index set has size $|I_{\partial}| = 3(n-1)$.

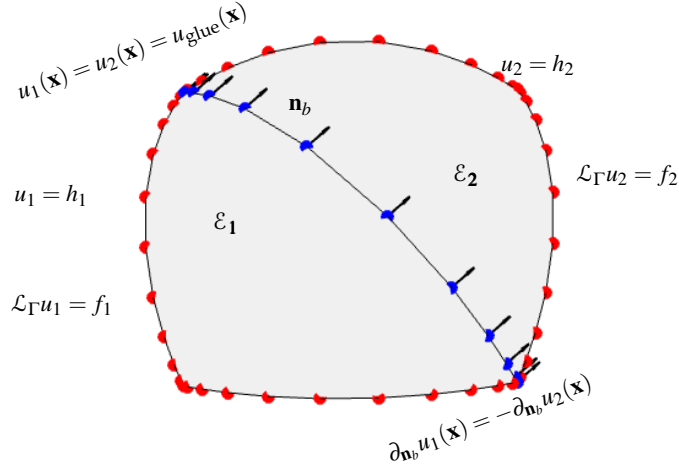
\begin{figure}[htb]
    \centering
    \begin{overpic}[width=0.7\textwidth]{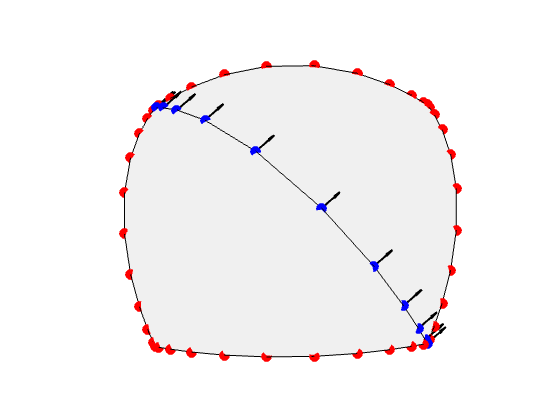}
        \put(32,36){\small$\mathcal{E}_{\bm{1}}$}
        \put(67,43){\small$\mathcal{E}_{\bm{2}}$}
        \put(3,40){\small$u_1 = h_1$}
        \put(74,59){\small$u_2 = h_2$}
        \put(2.5,24){\small$\mathcal{L}_{\Gamma} u_1 = f_1$}
        \put(85,43){\small$\mathcal{L}_{\Gamma}u_2 = f_2$}
        \put(43,53.3){\small$\mathbf{n}_{b}$}

        \put(3,55){%
          \rotatebox{20}{\small
            $u_1(\mathbf{x}) = u_2(\mathbf{x}) = u_{\text{glue}}(\mathbf{x})$
          }%
        }

        \put(57,2){%
          \rotatebox{20}{\small
            $\partial_{\mathbf{n}_b} u_1(\mathbf{x}) = -\partial_{\mathbf{n}_b} u_2(\mathbf{x})$
          }%
        }
    \end{overpic}

    \caption{
    Interface coupling of two surface elements, enforcing continuity of the solution and its binormal derivative. 
    Red points mark boundary collocation nodes, and blue points show aligned interface nodes used for coupling.
    }
    \label{fig:tri_two_glued_patches}
\end{figure}

Using the standard Schur--complement construction, directly analogous to the
quadrilateral case, we obtain the triangular solution operators
$\mathbf{S}_{T_k}$, whose dimensions are
$\frac{n(n+1)}{2} \times 3(n-2)$.
Having constructed $\mathbf{S}_{T_k}$ to solve the PDE locally on each
triangular element $T_k$, for $k = 1,\ldots,K$, the corresponding
Dirichlet-to-Neumann operator is defined by
$
\mathrm{DtN}_{T_k} = \mathcal{D}_{T_k}\,\mathbf{S}_{T_k}$,
which is a square matrix of size $3(n-2)\times 3(n-2)$. The merging strategy
described in~\ref{merge_section} then applies without modification.
Figure~\ref{fig:tri_two_glued_patches} illustrates the interface coupling
between two triangular surface elements, enforcing continuity of both the
solution and its binormal derivative.

\subsection{Computational asymptotic complexity}\label{sec:performance}
For an order-$n$ discretization using triangular spectral elements, the total
number of degrees of freedom associated with the mesh
$\{\widehat{T}_k\}_{k=1}^{K}$ scales as
$
N \sim K\,\frac{(n+1)(n+2)}{2}.
$
As in the quadrilateral based HPS method \cite{martinsson2013direct, fortunato2022highorder} , the computational complexity consists
of three main stages and is governed by the number of local degrees of freedom.
On each triangular patch, constructing the local solution operator
$\bm{S}_{\widehat{T}_k}$ and the corresponding Dirichlet-to-Neumann operator
$\mathrm{DtN}_{\widehat{T}_k}$ requires
$O(N n^{6})$ operations, since the local collocation matrices have size
$O(n^{2})$. The merging process over the hierarchical tree incurs a cost of
$O(N^{3/2} n^{3})$, while the overall complexity of the solve stage is
$O(n^{2} N \log N + N n^{3})$. The memory requirements scale similarly to the
solve stage, as each level of the hierarchy stores dense representations of the
solution and Dirichlet-to-Neumann operators.

    \begin{figure}[htb]
    \centering

    \begin{subfigure}[t]{0.48\textwidth}
        \centering
        \includegraphics[width=\textwidth]{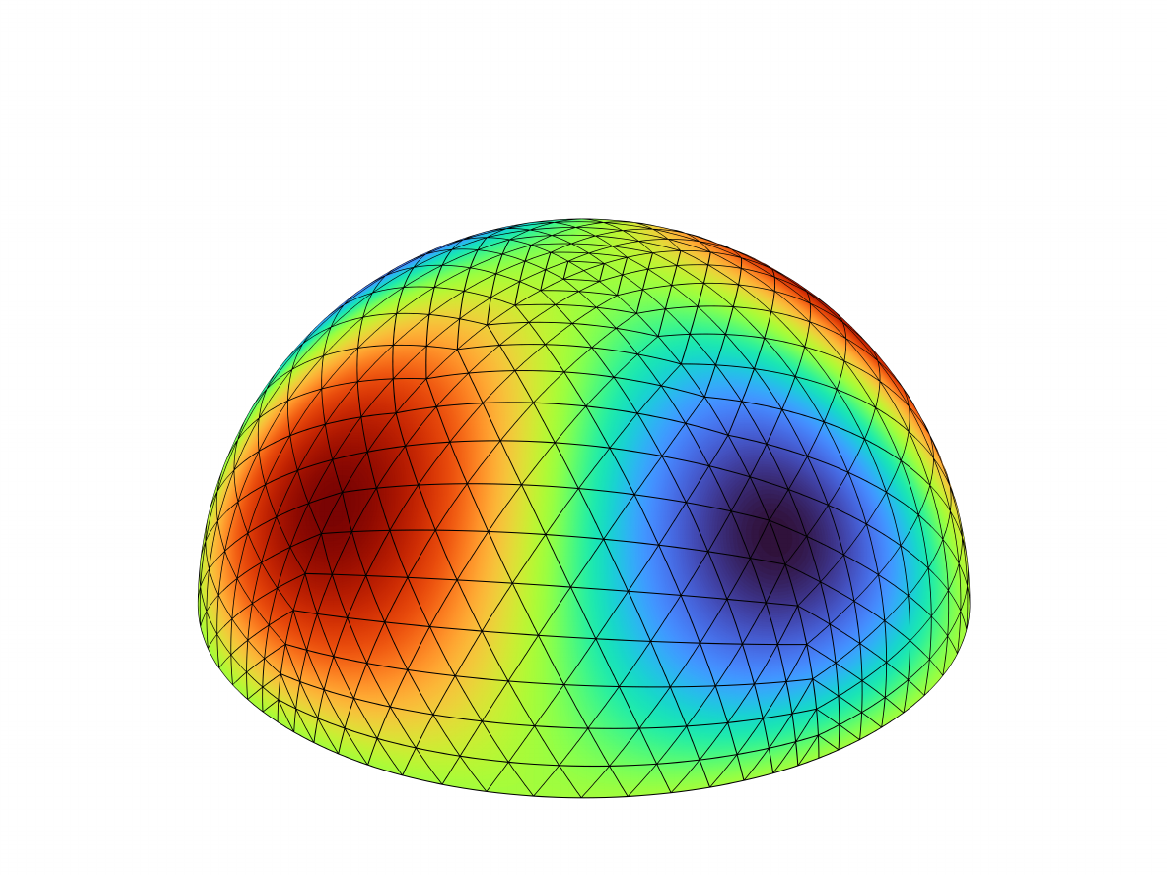}
    \end{subfigure}
    \hfill
    \begin{subfigure}[t]{0.48\textwidth}
        \centering
        \includegraphics[width=\textwidth]{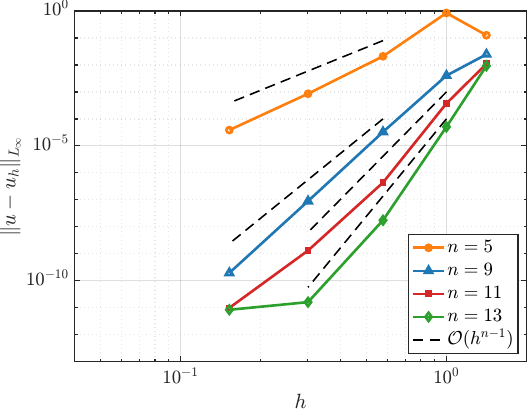}
    \end{subfigure}

    \caption{%
    A high-order triangular patch mesh is used to discretize the Laplace--Beltrami equation on the upper hemisphere. 
    The spectral discretization exhibits high-order convergence consistent with the rate \(\mathcal{O}(h^{n-1})\).%
    }
    \label{fig:convergence_sphere1}
\end{figure}
\paragraph{Laplace--Beltrami problem on the sphere}

Following the convergence experiments in \cite{fortunato2022highorder}, we study the Laplace--Beltrami problem on the unit sphere $\Gamma$
seeking \(u\) such that
\begin{equation}\label{eq:exp22}
    -\Delta_{\Gamma} u = f \qquad \text{on } \Gamma = \{\mathbf{x} \in \mathbb{R}^3 \mid \|\mathbf{x}\|_2 = 1\} .
\end{equation}
Exact solutions are chosen from the family of spherical harmonics \(Y_\ell^m\), which satisfy
\[
    -\Delta_{\Gamma} Y_\ell^m = \ell(\ell+1)\, Y_\ell^m .
\]

Two geometrical configurations are considered.  
First, the problem is restricted to the upper hemisphere, where we take $
    u(\mathbf{x}) = Y_3^2(\mathbf{x}),$
and impose Dirichlet boundary conditions along the equator, given by
\[
    h(\bm{x}) = 0.25 \sqrt{\frac{105}{\pi}} \, (x_1^2 - x_2^2)\, x_3 .
\]
Second, we consider the closed unit sphere and choose~$u(\mathbf{x}) = Y_{20}^{10}(\mathbf{x})$.
In both cases, the forcing term is defined as
\[f(\mathbf{x}) = -\ell(\ell+1)\, Y_\ell^m(\mathbf{x}),\] 
so that the exact solution satisfies \eqref{eq:exp22}.

\begin{figure}[htb]
    \centering

    \begin{subfigure}[t]{0.48\textwidth}
        \centering
        \includegraphics[width=\textwidth]{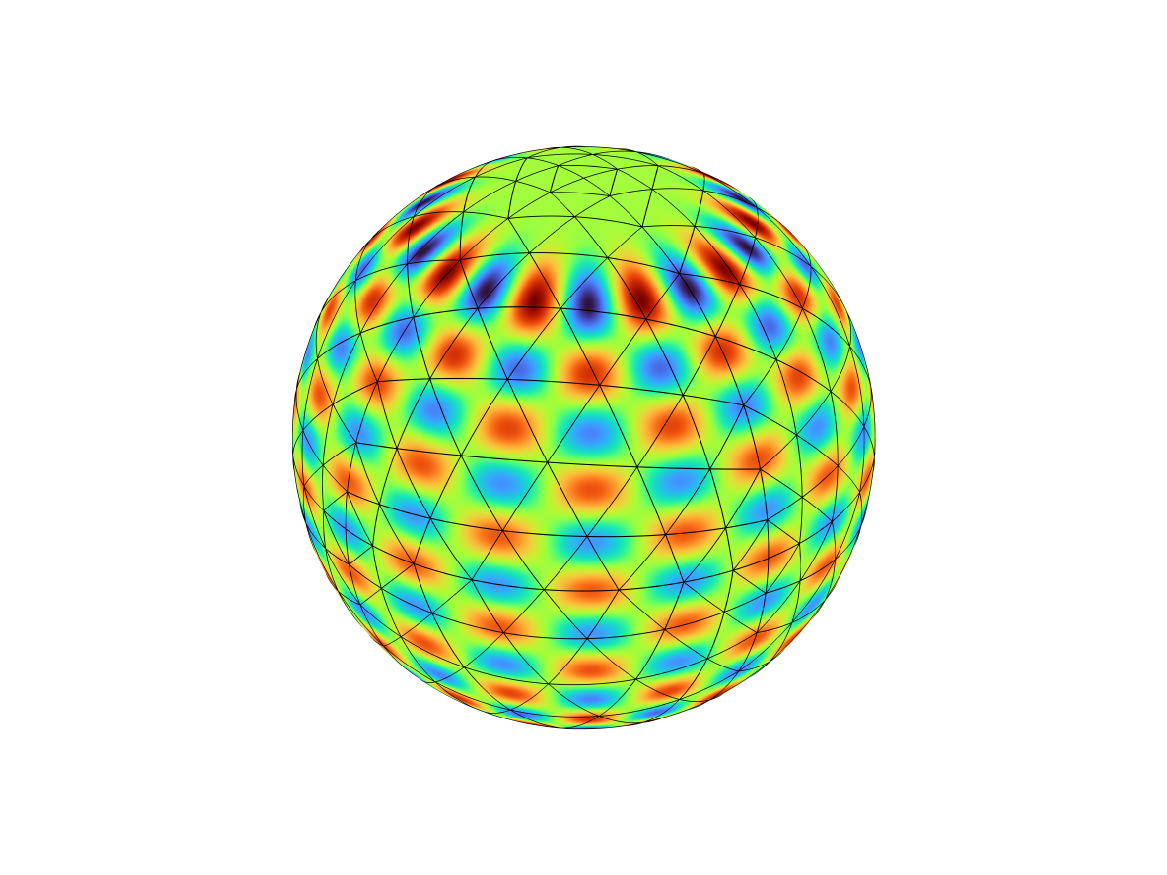}
    \end{subfigure}
    \hfill
    \begin{subfigure}[t]{0.48\textwidth}
        \centering
        \includegraphics[width=\textwidth]{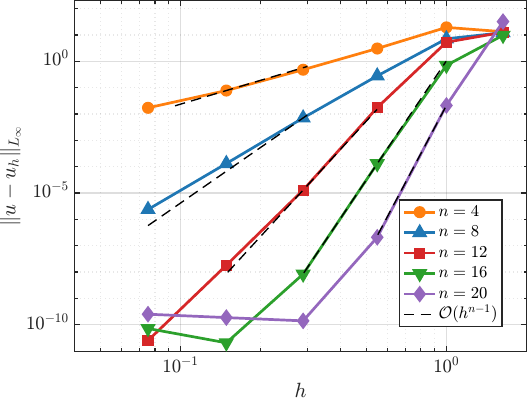}
    \end{subfigure}
\caption{
A high-order triangulated sphere mesh is used to approximate a spherical harmonic. The observed algebraic decay aligns with the fit rate \( h^{n-1} \).
}
  \label{fig:convergence_sphere}
\end{figure}

The solution is evaluated on a triangulated approximation of the sphere and used to assess the convergence of the triangle based hierarchical HPS discretization.  
Figure~\ref{fig:convergence_sphere1} shows the spatial convergence under \(h\)-refinement for polynomial degrees \(n = 5, 9, 11,\) and \(13\). The relative \(L_\infty\)-error decays algebraically with the mesh size \(h\), following the expected rate \(\mathcal{O}(h^{\,n-1})\); the corresponding reference slopes are included for comparison. The surface plots (left) depict the computed solutions on the triangulated geometry. Figure~\ref{fig:convergence_sphere} shows analogous convergence results for the closed unit sphere, with algebraic \(L_\infty\)-error decay at the expected rate \(\mathcal{O}(h^{\,n-1})\). These results demonstrate that the
convergence behavior of the triangle based HPS discretization is comparable to
that of the quadrilateral based HPS method~\cite{fortunato2022highorder}.

\section{Time-dependent equations}\label{time_depdn}
While the HPS scheme is designed to solve stationary problems modeled by linear elliptic
PDEs, it can also be useful for accelerating  time-dependent problems modeled by parabolic PDEs. 
We start by considering the reaction-diffusion systems, which are widely regarded as key mechanisms for pattern formation in a variety of contexts, including biological, chemical, physical, and even economic processes. A general reaction–diffusion system describing \(N\) interacting species defined on a closed, smooth surface \(\Gamma \subset \Omega \subset \mathbb{R}^{d+1}\) can be written in the form:

\begin{equation}\label{eq:main_diff}
\frac{\partial \bm{u}}{\partial t} =\nabla_\Gamma \cdot \left( \mathbf{D} \nabla_\Gamma \bm{u} \right)+ \bm{F}(\bm{u}) ,
\end{equation}

where $\bm{u}:=(u_1, u_2, \dots, u_N),\;$ $\bm{F}$ represents the reaction kinetics, also referred to as the source term, and $\mathbf{D}$ denotes
the diffusion tensor. The reason a fast direct solver for elliptic equations is useful is that each time step in \eqref{eq:main_diff} requires solving an elliptic problem, and when there are many timesteps it can become advantageous to use a direct solver. This is of course especially true if the complexity (as it is here) of the linear solver is good.
 
In the discretization process, we handle the spatial discretization, as described in the previous sections, using a domain decomposition approach with spectral collocation on each element. This yields a system of ordinary differential equations (ODEs) in time, given by:

\begin{equation}\label{discret_tdp}
\frac{d\bm{u}}{dt} = L_\Gamma \bm{u} + \bm{F}(\bm{u}),
\end{equation}

 The term \( L_\Gamma \bm{u} \) arises from the diffusion components, while \( \bm{F}(\bm{u}) \) originates from the reaction components.  Because the diffusive term is typically stiff \cite{strikwerda2004finite}, the use of explicit schemes usually necessitates excessively small time steps. This can result in computations which are prohibitively expensive in three
spatial dimensions. Fully implicit treatment, however, requires the implicit treatment of the nonlinear reaction term, \( \bm{F}(\bm{u})\), at every time step. This can be
particularly expensive and undesirable because the Jacobian of \( \bm{F}(\bm{u})\) could be dense and is typically non-definite,  or non-symmetric which makes fast iterative solution techniques \cite{varga1962matrix} less efficient and more difficult to implement.  
Moreover, explicitly handling the nonlinear reaction term 
\(\bm{F}(\bm{u})\)  is easy to implement and adds relatively little computational effort per time step. Additionally, many well-known time-stepping methods applied to \eqref{discret_tdp} are either first-order (such as backward Euler) or result in only a weak reduction of high-frequency error components (such as Crank-Nicolson). In this work,  we employ the Implicit-Explicit Backward Differentiation Formula (IMEX-BDF) family of schemes \cite{ascher1995implicit}, combining backward differentiation for implicit terms with Adams-Bashforth for explicit terms. 

Let $\Delta t > 0$ denote the time step, and $\bm{u}^n(\mathbf{x}) \approx \bm{u}(\mathbf{x}, n \Delta t)$ represent the approximate solution at step $n$. Time discretization with the $M^{\text{th}}$ order IMEX-BDF scheme yields a steady-state problem at each step:

\begin{equation}\label{eq:imex_bdf}
\left( I - \omega \Delta t L_\Gamma \right) \bm{u}^{n+1} = \sum_{i=0}^{M-1} a_i \bm{u}^{n-i} + \Delta t \sum_{i=0}^{M-1} b_i\,\bm{F}(\bm{u}^{n-i}),
\end{equation}
which can be written more compactly as
\[
A \bm{u}^{n+1} = \bm{f}^n,
\]
where
\[
A := I - \omega \Delta t L_\Gamma, \quad 
\bm{f}^n := \sum_{i=0}^{M-1} a_i \bm{u}^{n-i} + \Delta t \sum_{i=0}^{M-1} b_i \bm{F}(\bm{u}^{n-i}).
\]
For instance, using the IMEX-BDF1 scheme, where $\omega = 1, \; a_0 = 1, \;b_0 = 1$,
yields the equation
\begin{equation}\label{eq:imex_bdf1}
\left( I - \Delta t L_\Gamma \right) \bm{u}^{n+1} = \bm{u}^{n} + \Delta t \,\bm{F}(\bm{u}^{n}),
\end{equation}
which must be solved once per time step to compute $\bm{u}^{n+1}$ from $\bm{u}^{n}$, yielding a cost per time of $\mathcal{O}(N \log K + N(n+1))$

The values of \( \omega \), \( a_i \), and \( b_i \) for IMEX-BDF schemes of order 2 to 4 are given as follows:
\[
\begin{aligned}
&\text{IMEX-BDF2:} \; &&\omega = \tfrac{2}{3}, && a = \left(\tfrac43, -\tfrac13\right), && b = \left(\tfrac43, -\tfrac23\right), \\
&\text{IMEX-BDF3:} \; &&\omega = \tfrac{6}{11}, && a = \left(\tfrac{18}{11}, -\tfrac{9}{11}, \tfrac{2}{11}\right), && b = \left(\tfrac{18}{11}, -\tfrac{18}{11}, \tfrac{6}{11}\right), \\
&\text{IMEX-BDF4:} \; &&\omega = \tfrac{12}{25}, && a = \left(\tfrac{48}{25}, -\tfrac{36}{25}, \tfrac{16}{25}, -\tfrac{3}{25}\right), && b = \left(\tfrac{48}{25}, -\tfrac{72}{25}, \tfrac{48}{25}, -\tfrac{12}{25}\right).
\end{aligned}
\]

 Regarding the accuracy of IMEX-BDF schemes, the following theorem provides important insights:

\begin{theorem}[\cite{ascher1995implicit}, Theorem 2.1]\label{conv.theorem}
An \(s\)-step IMEX-BDF scheme, as defined in \eqref{eq:imex_bdf}, cannot achieve an order of accuracy higher than \(s\).
\end{theorem}
We conclude this section by evaluating the accuracy and convergence of the
triangle based HPS scheme for a time-dependent problem. As a benchmark with a
known analytic solution, we consider isotropic diffusion on the surface of the
unit sphere. The governing equation is the surface diffusion equation,

\begin{equation}
    \frac{\partial u}{\partial t} = \Delta_{\Gamma} u \quad \text{on } \Gamma = \{\mathbf{x} \in \mathbb{R}^3 \mid \|\mathbf{x}\|_2 = 1\},
\end{equation}
The exact solution in spherical coordinates \( (\vartheta, \varphi) \), where \( \vartheta \in [0,\pi] \) is the polar angle and \( \varphi \in [0, 2\pi) \) is the azimuthal angle, is obtained via expansion in spherical harmonics \( Y_\ell^m \):
\begin{equation*}
    u(t, \vartheta, \varphi) = \sum_{\ell = 0}^{\infty} \sum_{m = -\ell}^{\ell} c_{\ell m}(0) Y_\ell^m(\vartheta, \varphi) e^{- \ell(\ell + 1)t},
\end{equation*}
with modal coefficients determined by the projection of the initial data:
\begin{equation*}
    c_{\ell m}(0) = \int_{\Gamma} (-1)^m Y_\ell^{-m}(\vartheta, \varphi)\, u(0, \vartheta, \varphi)\, dS.
\end{equation*}

To avoid effects due to spectral truncation and to isolate the numerical error, we consider the initial condition
\begin{equation*}
    u(0, \vartheta, \varphi) = Y_1^0(\vartheta, \varphi) = \sqrt{\frac{3}{4\pi}} \cos \vartheta,
\end{equation*}
which leads to a closed-form analytic solution due to the orthogonality of the spherical harmonics:
\begin{equation*}
    u(t, \vartheta, \varphi) = Y_1^0(\vartheta, \varphi)\, e^{-2t}.
\end{equation*}

\begin{figure}[htbp]
  \centering
  \begin{subfigure}[b]{0.45\textwidth}
    \includegraphics[width=1.\textwidth]
    {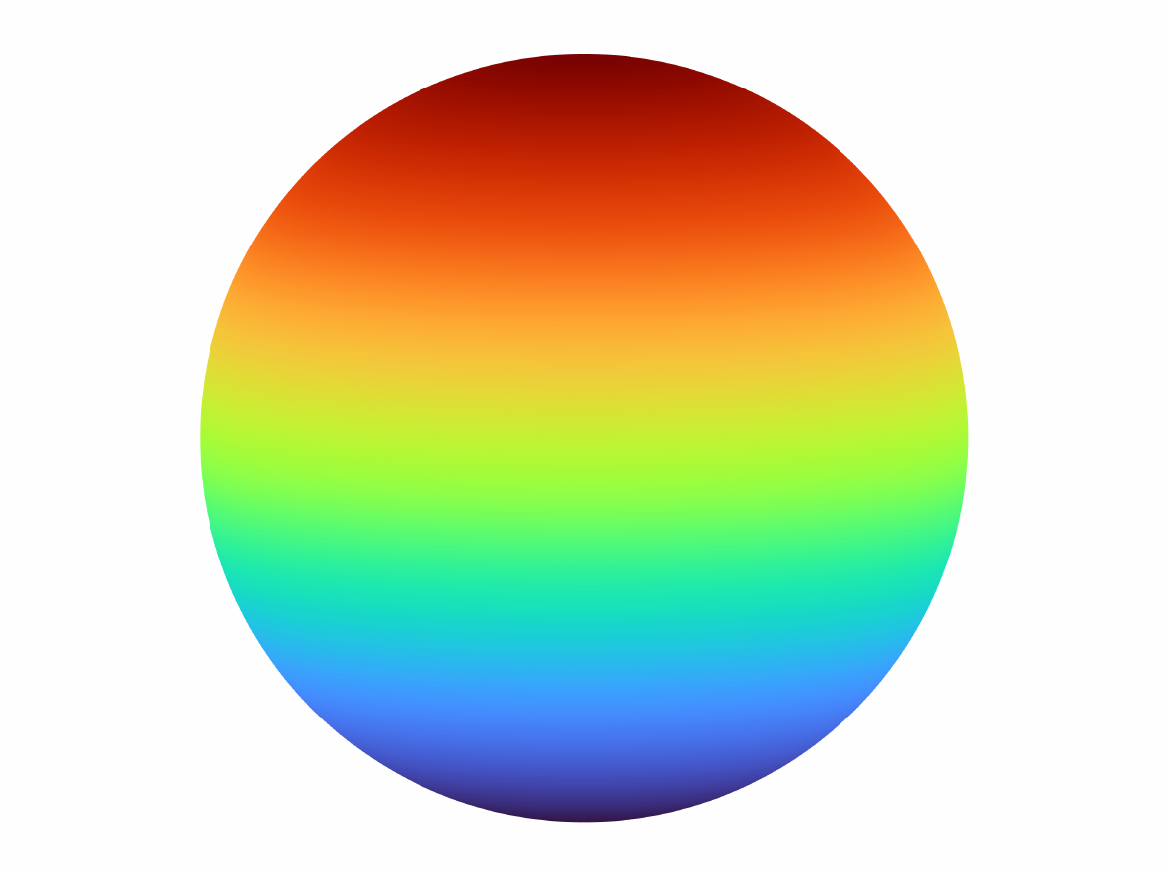}
  \end{subfigure}
  \hspace{0.5cm}
  \begin{subfigure}[b]{0.45\textwidth}
    \includegraphics[width=\textwidth]
    {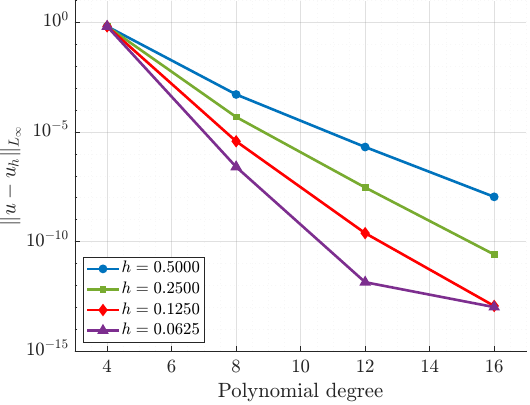}
    \label{fig:diffsphere_error_linf}
  \end{subfigure}
  \caption{%
 (Left) Computed solution at $t = 1$ using time step $\Delta t = 10^{-3}$. 
    (Right) $L_\infty$-error versus polynomial degree for different mesh sizes $h$. 
    The simulations were performed using the implicit--explicit backward differentiation formula (IMEX--BDF4) scheme.
  }
  \label{fig:diffsphere_error}
\end{figure}
The left panel of Figure~\ref{fig:diffsphere_error} shows the computed solution at the final time $t=1$, while the right panel illustrates the $L_\infty$-error as a function of polynomial degree for different mesh sizes. As expected, increasing the polynomial degree results in exponential convergence of the error, especially for finer meshes.

\section{Numerical study of spatial pattern
formation in Turing systems}\label{turing_model}

In Section~\ref{time_depdn}, we examined the pure surface diffusion equation as a benchmark problem to assess the accuracy and convergence of the proposed triangle based HPS scheme. Having established the method’s robustness for this simpler case, we now turn to a more complex and biologically motivated setting: the Turing model of spatial pattern formation, first proposed by Turing~\cite{turing1990chemical}. He showed that a system of two reacting and diffusing
chemicals could give rise to spatial patterns in chemical concentrations from initial near-homogeneity. A comprehensive overview of commonly used reaction kinetics and their underlying motivation can be found in~\cite{maini1997spatial}. In this section, we consider the Turing model on a surface $\Gamma$ involving an activator–inhibitor system ~\cite{barrio1999two}. Our aim is to investigate how the domain geometry, nonlinearities, and coupling of such systems influence the emergence of spatial patterns. The specific form of the system considered is:

\begin{subequations}\label{eq:reaction_diffusion}
\begin{align}
    \frac{\partial u_1}{\partial t} &= \delta_{u_1} \Delta_\Gamma u_1 + \alpha u_1 \left(1 - r_1 u_2^2\right) + u_2 \left(1 - r_2 u_1\right)\quad\quad\quad \text{with}\;   u_1(\mathbf{x},0)=u_{1}^{0}(\mathbf{x}),\label{eq:reaction_diffusion_u11}  \\
    \frac{\partial u_2}{\partial t} &= \delta_{u_2} \Delta_\Gamma u_2 + \beta u_2 \left(1 + \frac{\alpha r_1}{\beta} u_1 u_2\right) + u_1 \left(\gamma + r_2 u_2\right)\quad \text{with}\;   u_2(\mathbf{x},0)=u_{2}^{0}(\mathbf{x}), \label{eq:reaction_diffusion_u22} 
\end{align}
where \( \alpha,\beta,\gamma\), \( r_1 \),\;\( r_2 \),\; \( \delta_{u_1}\) and \(\delta_{u_2} \) are the parameters of the reaction-diffusion system. In the context of Eq.~\eqref{eq:main_diff}, we have the following: \[
\mathbf{D} = \begin{pmatrix}
\delta_{u_1} & 0 \\
0 & \delta_{u_2}
\end{pmatrix}, \quad
\mathbf{u} = \begin{pmatrix}
u_1 \\
u_2
\end{pmatrix}, \quad
\mathbf{F} = \begin{pmatrix}
\alpha u_1 \left(1 - r_1 u_2^2\right) + u_2 \left(1 - r_2 u_1\right) \\
\beta u_2 \left(1 + \frac{\alpha r_1}{\beta} u_1 u_2\right) + u_1 \left(\gamma + r_2 u_2\right)
\end{pmatrix}.
\]

\end{subequations}
In this system, \( u_1 \) and \( u_2 \) are morphogens with \( u_1 \) as the "activator" and \( u_2 \) as the "inhibitor". If \( \alpha = -\gamma \), then \( (u_1, u_2) = (0, 0) \) is a unique equilibrium point of this system. The reaction term contributes to the formation of concentration peaks of $u_1$ and $u_2$, whereas the diffusion term tends to smooth these peaks. The interplay between these opposing processes—reaction-driven peak formation and diffusion-induced peak smoothing—leads to the emergence of characteristic Turing patterns. 

To understand when such patterns arise, we consider the conditions under which the reaction–diffusion system \eqref{eq:reaction_diffusion} supports them. Specifically, model \eqref{eq:reaction_diffusion} exhibits Turing pattern formation when the following two conditions, known as the Turing criteria, are satisfied:
\begin{enumerate}[label=(\roman*)]
\item In the absence of diffusion, the system tends toward a spatially uniform, linearly stable steady state.
\item The steady state becomes unstable in the presence of diffusion due to the introduction of random perturbations.
\end{enumerate}

It is well known that patterns produced by the Turing model are influenced by domain geometry~\cite{bunow1980pattern}, and the nonlinear reaction term \( \mathbf{F}(\mathbf{u}) \), containing quadratic and cubic interactions, plays a key role in driving pattern development. To examine how geometry and nonlinearities affect pattern emergence, we simulate the reaction diffusion system~\eqref{eq:reaction_diffusion} on curved surfaces using the high-order triangle based HPS scheme.

\begin{figure}[h!]
    \centering
    \setlength{\tabcolsep}{2pt}
    \renewcommand{\arraystretch}{0.5}
    \begin{tabular}{ccc}
     \includegraphics[width=0.32\textwidth]{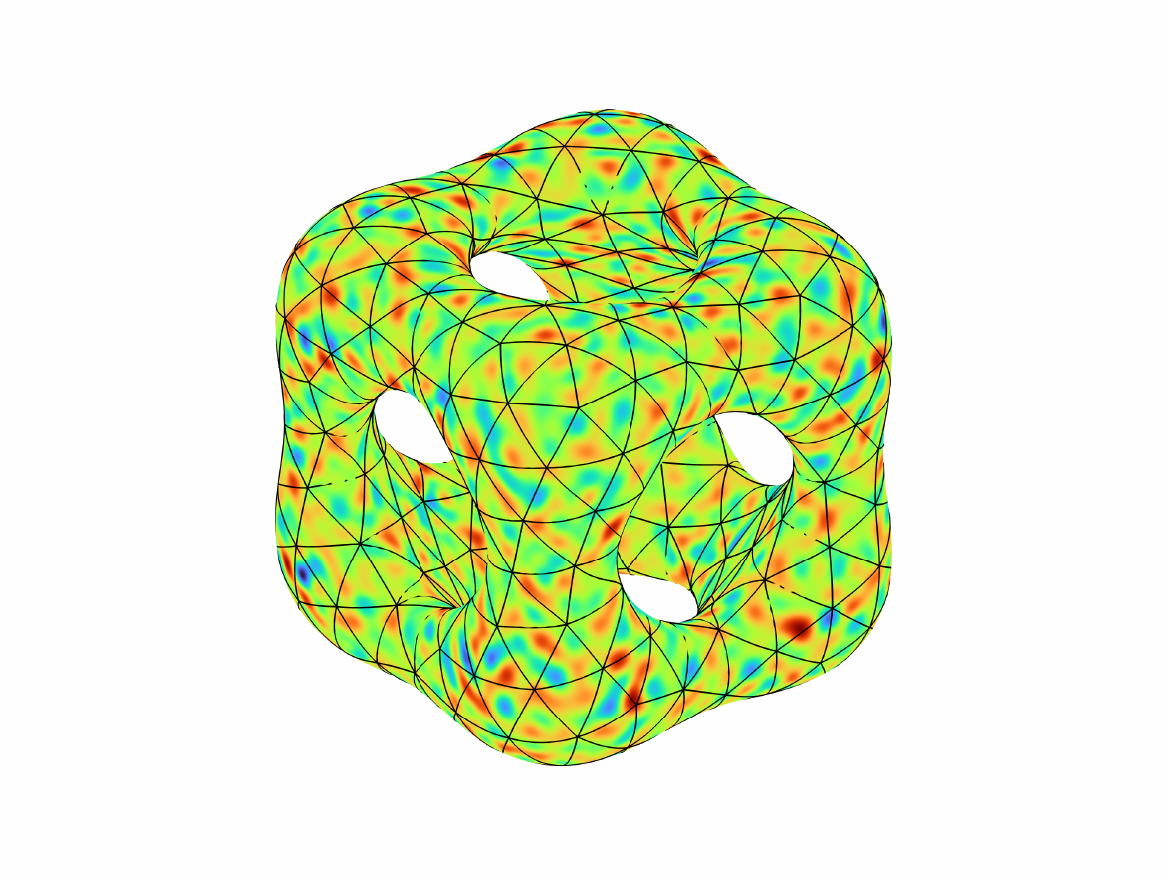} &
        \includegraphics[width=0.32\textwidth]{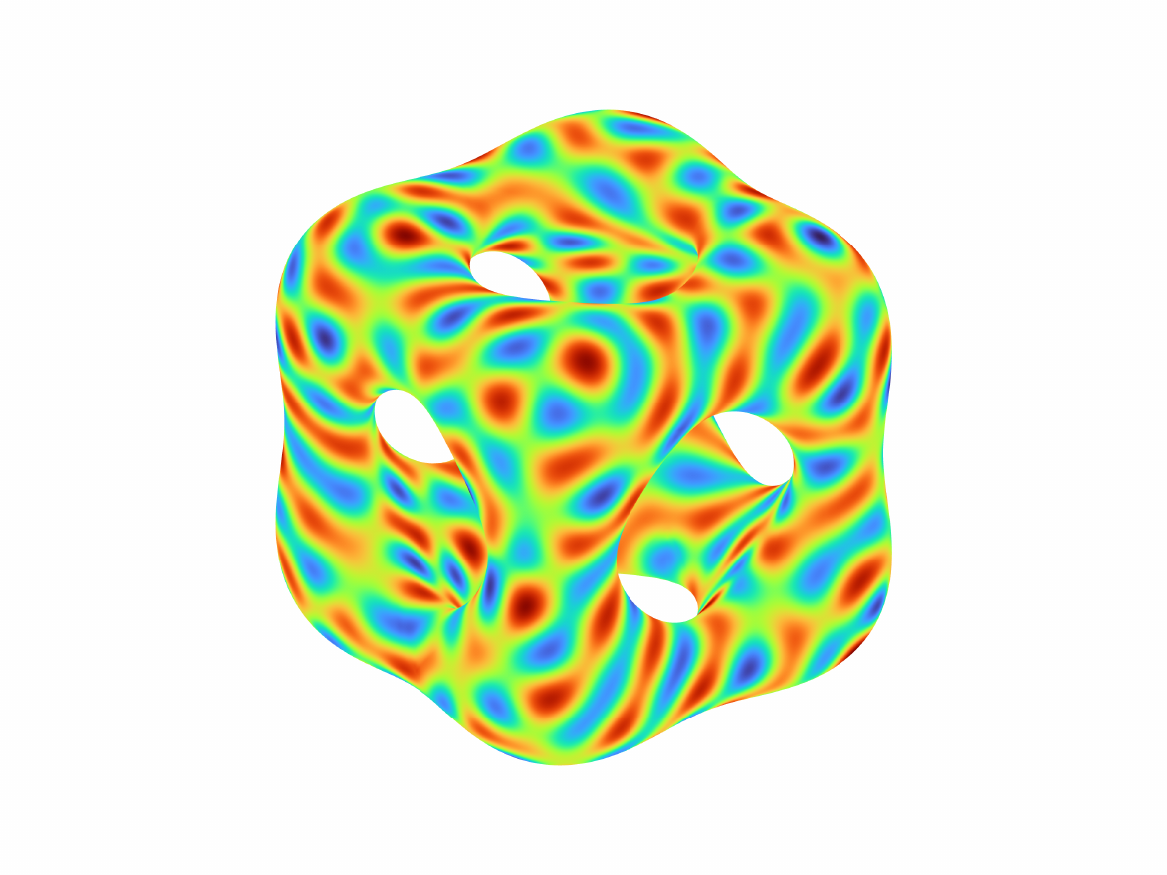} &
        \includegraphics[width=0.32\textwidth]{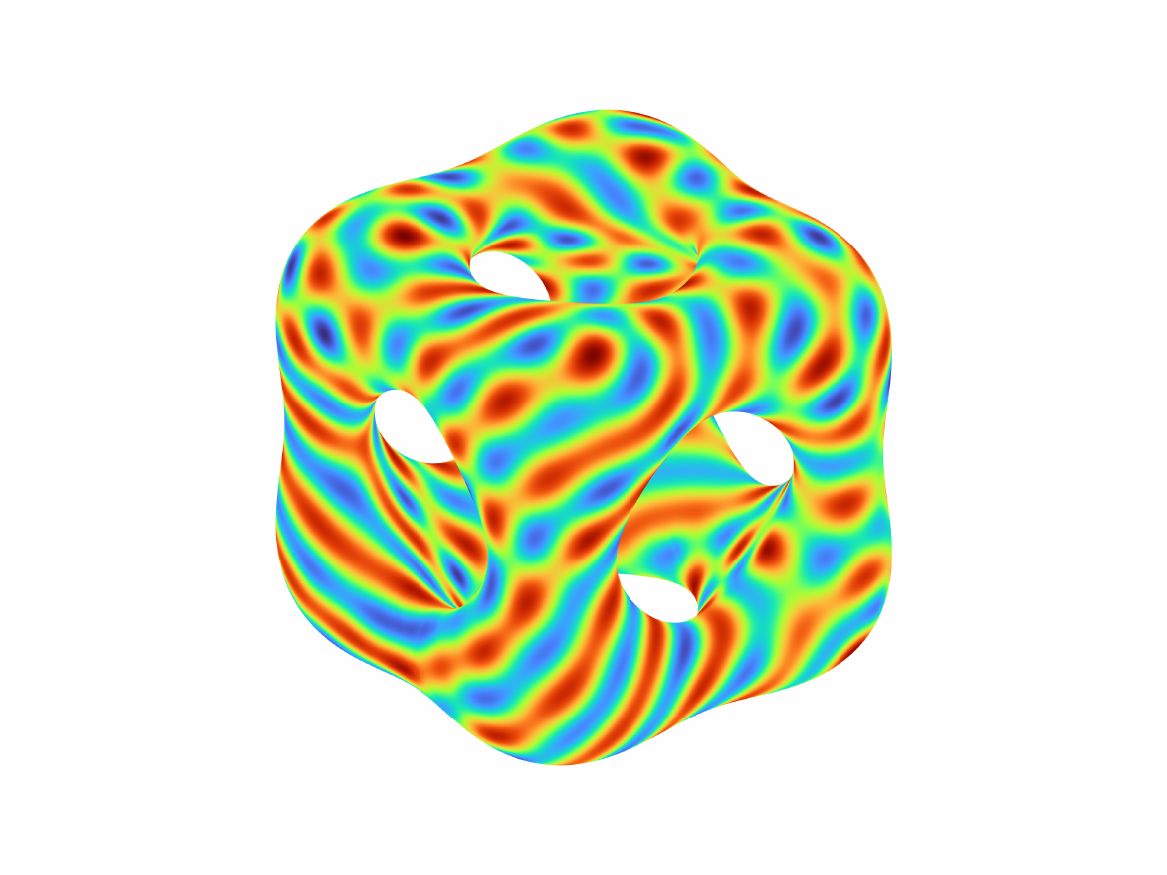} \\
        \includegraphics[width=0.32\textwidth]{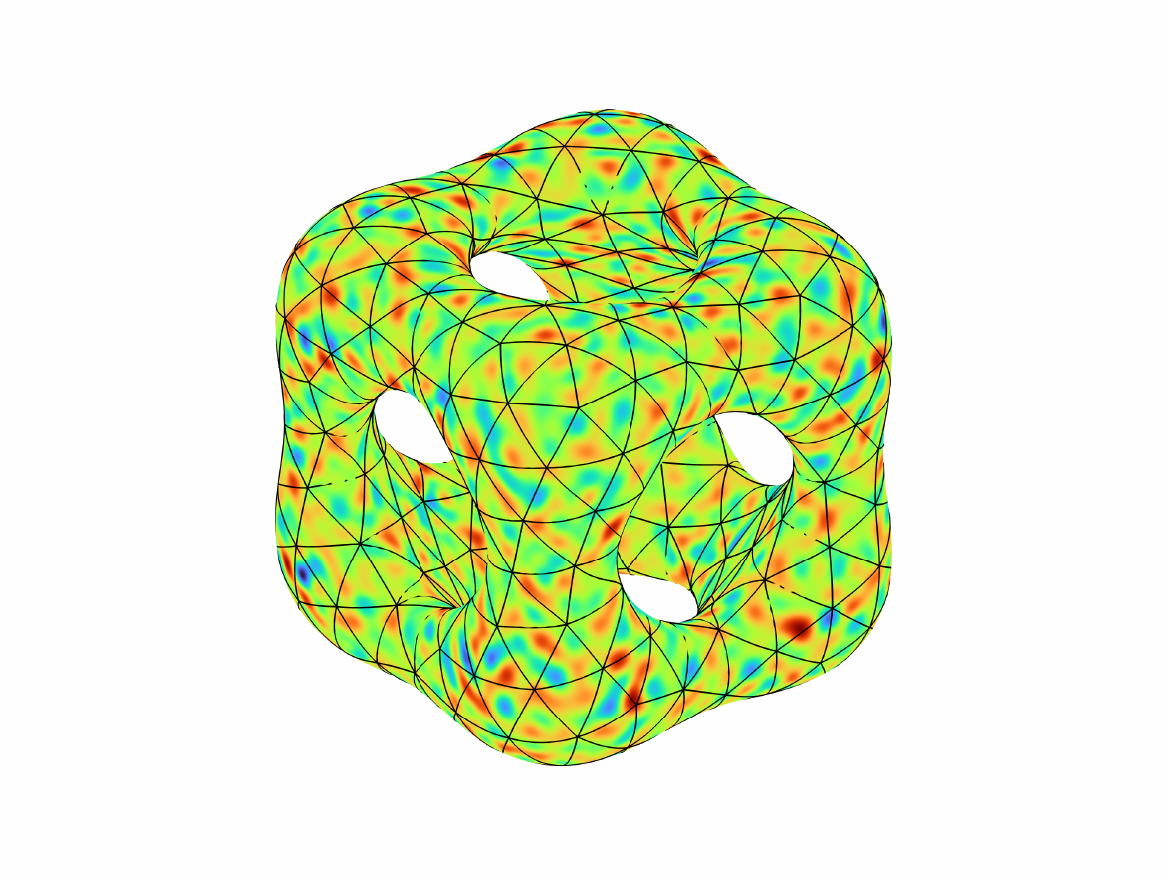} &
        \includegraphics[width=0.32\textwidth]{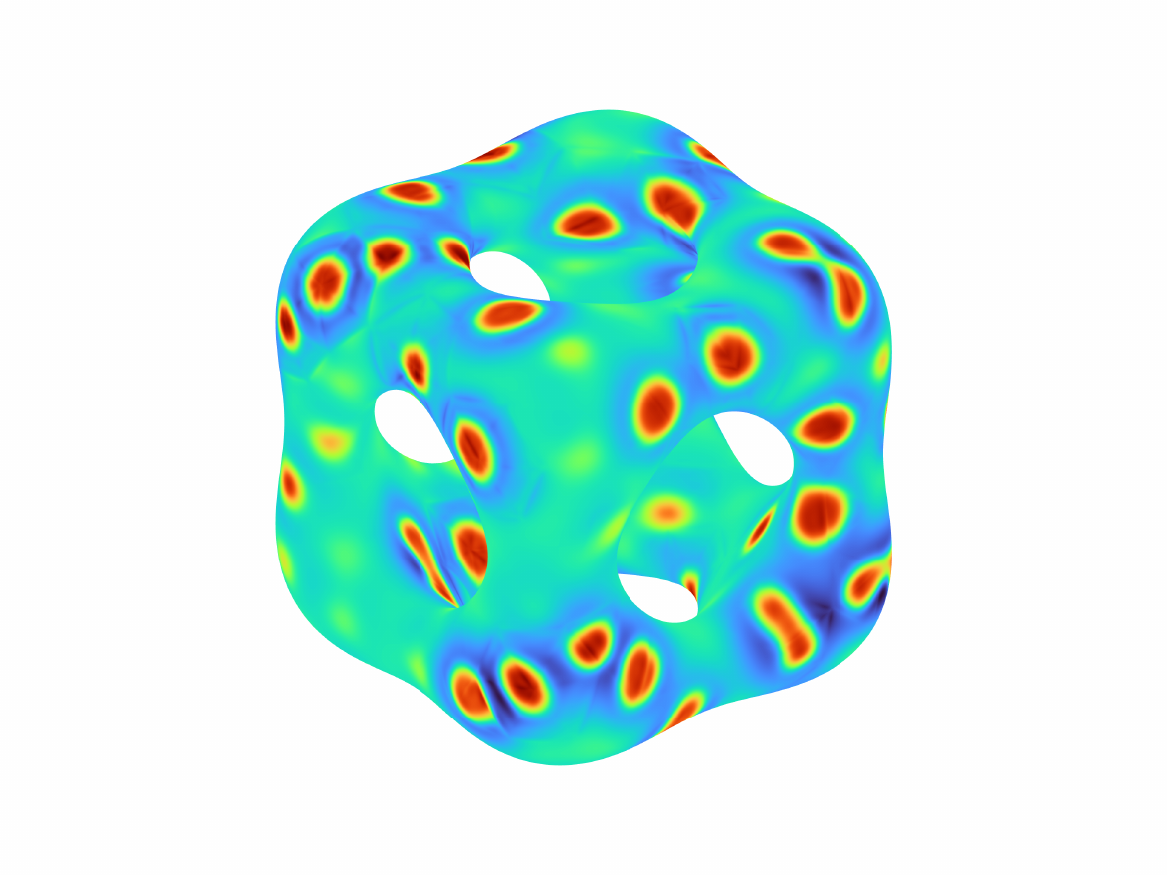} &
        \includegraphics[width=0.32\textwidth]{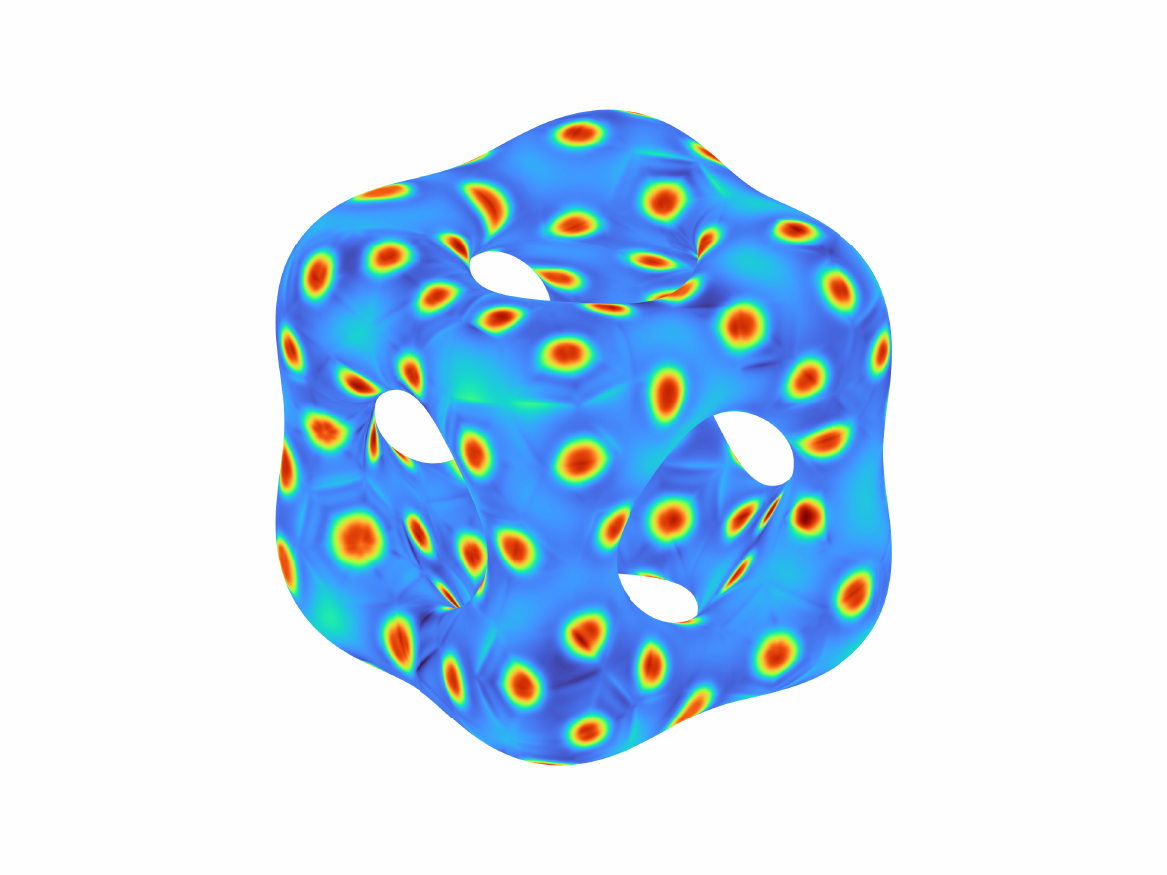}
    \end{tabular}
    \vspace{-0.8em}
    \caption{Turing model solution \( u_1 \) on a Swiss cheese surface. Top: \( r_1 = 3.5 \), \( r_2 = 0 \) at \( t = 0, 200, 600 \). Bottom: \( r_1 = 0.02 \), \( r_2 = 0.2 \) at \( t = 0, 20, 200 \). Simulated with IMEX-BDF1.}
    \label{fig:Turing_system_swiss}
\end{figure}

The resulting stable patterns exhibit either spot-like or stripe-like structures, depending on the values of the coupling parameters \( r_1 \) and \( r_2 \). The cubic coupling parameter \( r_1 \) promotes the emergence of stripes, whereas the quadratic coupling parameter \( r_2 \) tends to favor spot formation~\cite{barrio1999two}.  
Our numerical simulations indicate that, in general, spot patterns are more robust than stripes and reach a steady state significantly faster. Stripe patterns only emerge for very small values of \( r_2 \), and their orientation varies depending on the initial conditions.  
In Figure~\ref{fig:Turing_system_swiss}, we illustrate examples of these basic cases\footnote{The  Swiss cheese surface defined
implicitly by $
f(\mathbf{x}) =
(x^2 + y^2 - 4)^2
+ (z^2 - 1)^2
+ (y^2 + z^2 - 4)^2
+ (x^2 - 1)^2
+ (z^2 + x^2 - 4)^2
+ (y^2 - 1)^2
- 15$.}. For the simulations, we adopt the parameter values~\footnote{Parameters follow the non-dimensionalized Turing model in~\cite{bunow1980pattern}, where all coefficients are dimensionless. The diffusion ratio ensures scale separation needed for pattern formation.} from~\cite{bunow1980pattern}:   \( \alpha = 0.899,\; \beta=-0.91, \gamma=-\alpha \), and \( \delta_{u_1} = 0.516 \, \delta_{u_2} \), with \( \delta_{u_2} = 5\cdot10^{-3} \). The initial conditions are taken to be random: 
\( u_1(\mathbf{x}, 0) = \operatorname{rand}(\mathbf{x}) \), 
\( u_2(\mathbf{x}, 0) = \operatorname{rand}(\mathbf{x}) \).



\begin{figure}[htbp]
	\centering
        \setlength{\tabcolsep}{2pt}
    \renewcommand{\arraystretch}{0.5}
	\begin{tabular}{ccc}
        	~~\includegraphics[width=0.31\textwidth]{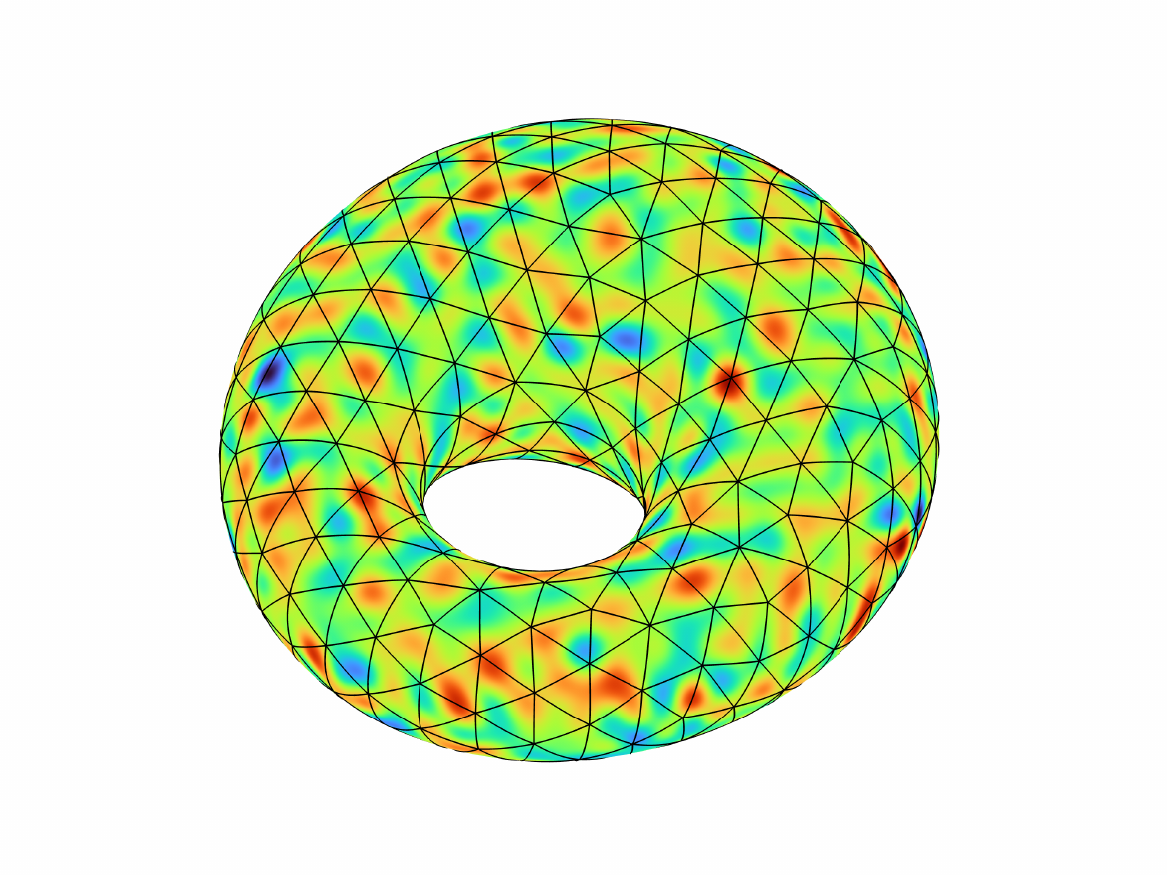} &
	~\includegraphics[width=0.31\textwidth]{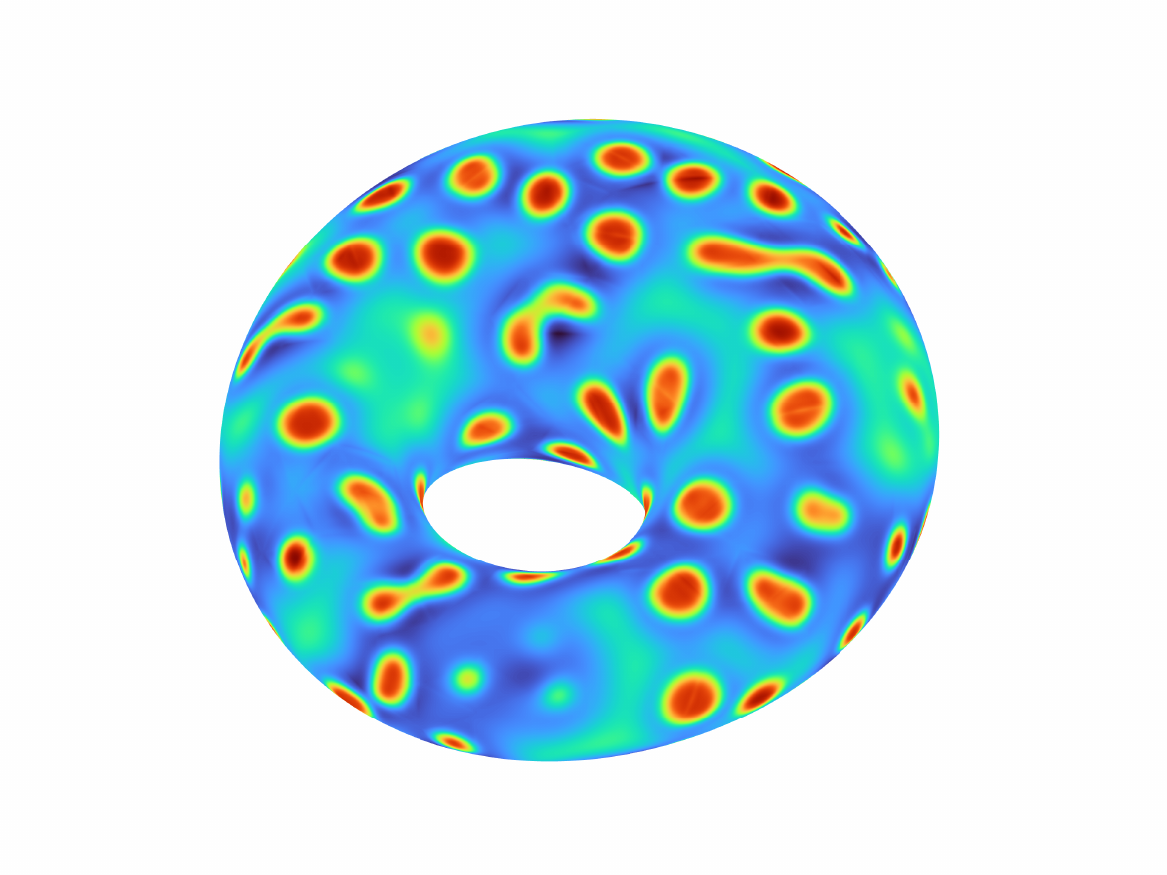} &
	~~\includegraphics[width=0.31\textwidth]{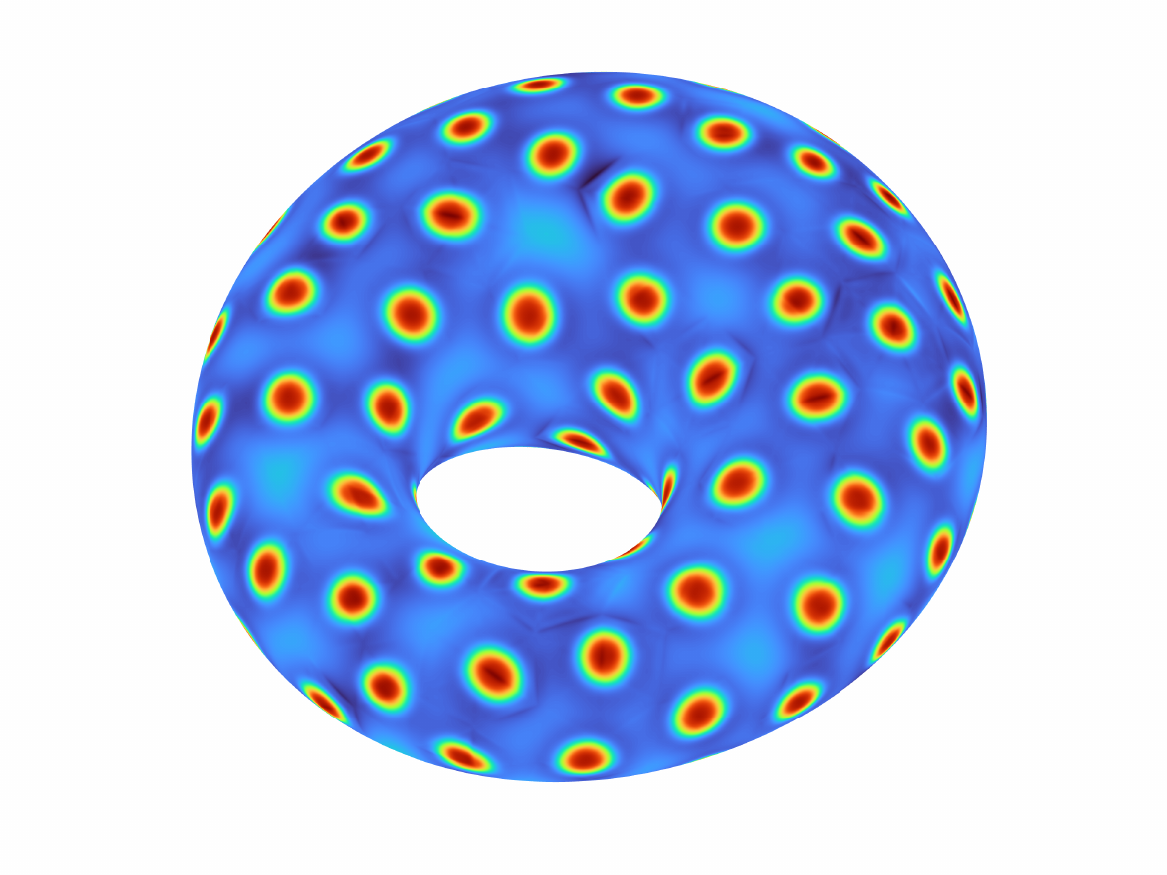}
    \\
    \includegraphics[width=0.31\textwidth]{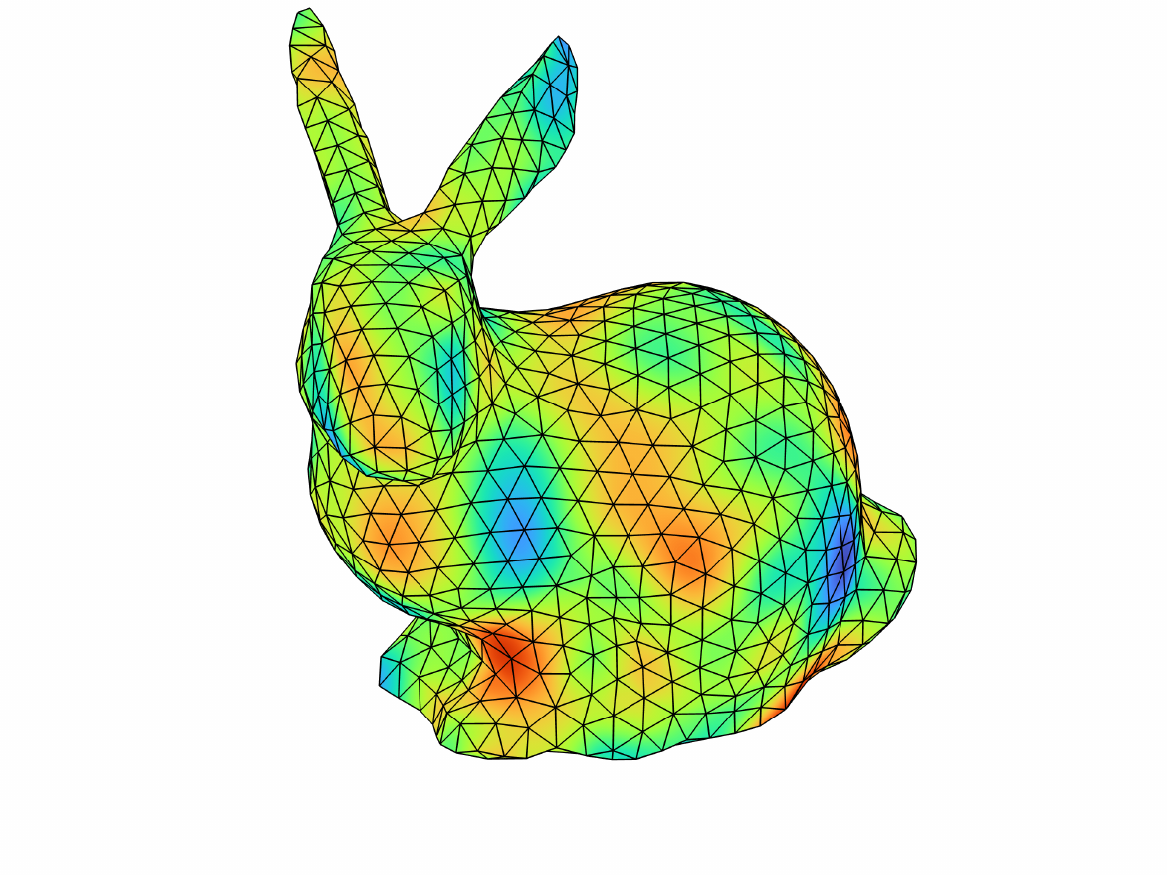} &
	\includegraphics[width=0.31\textwidth]{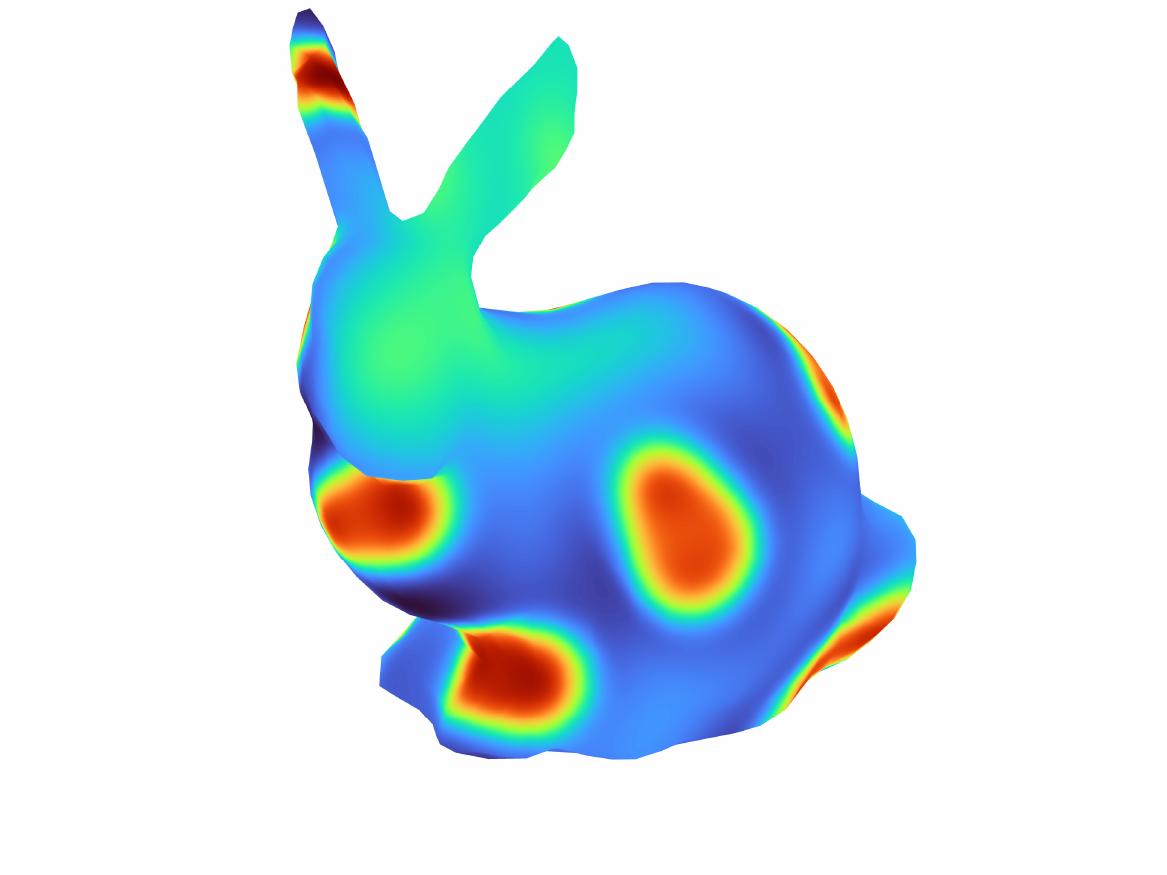} &
	\includegraphics[width=0.31\textwidth]{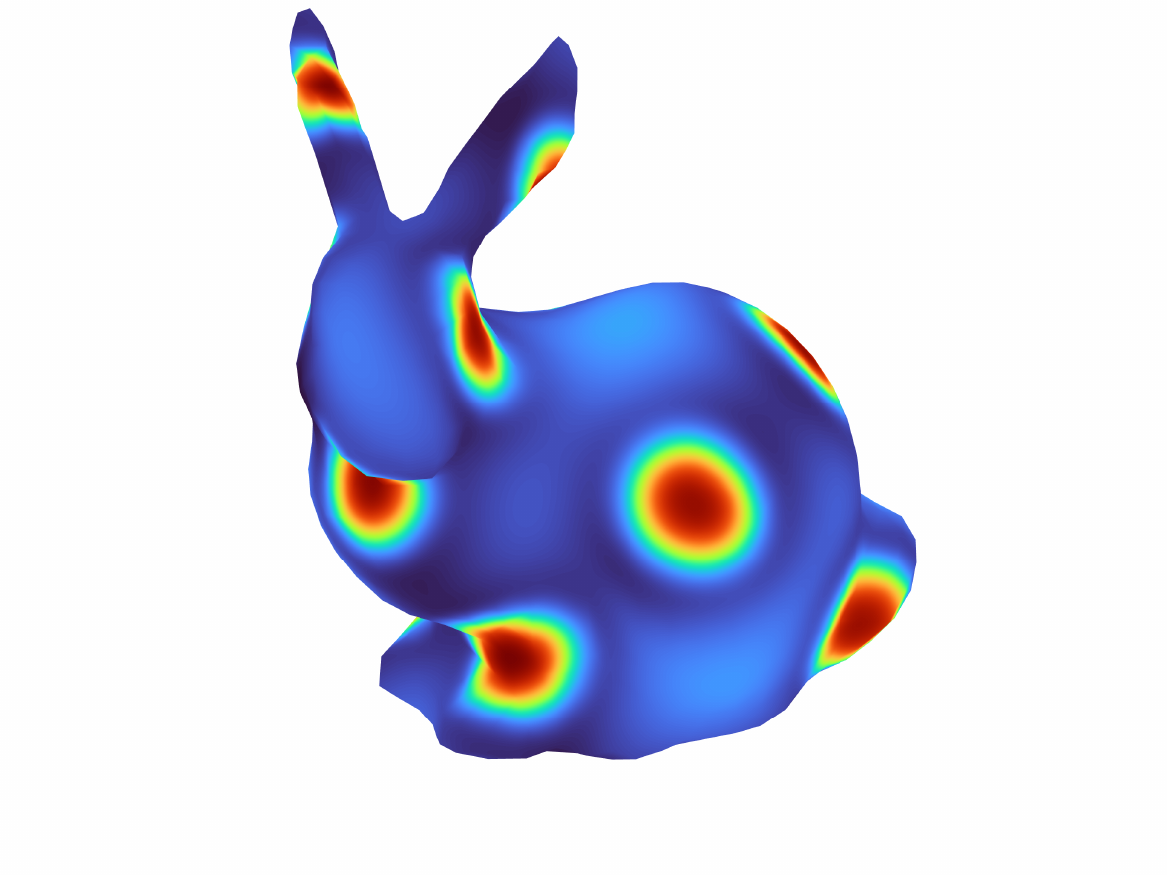} 

	\end{tabular}
\caption{Turing model solution \( u_1 \) on the asymmetric torus and Stanford Bunny at times \( t = 0 \), \( 20 \), and \( 200 \) using IMEX-BDF1 with \( r_1 = 0.02 \), \( r_2 = 0.2 \).}
\label{fig:Turing system_combined}
\end{figure}

Next, we simulate the Turing model on two representative surfaces to investigate
the influence of geometry on pattern formation. The first geometry is an
asymmetric torus,\footnote{%
The surface is given implicitly by $f(\mathbf{x}) =
\bigl(x^2 + y^2 + z^2 - d^2 + b^2\bigr)^2
- 4\bigl(a x + c^2 d\bigr)^2
- 4 b^2 y^2,$
with parameters \( a = 2 \), \( b = 1.9 \), \( d = 1 \), and
\( c^2 = a^2 - b^2 \).}
followed by a triangulated Stanford Bunny surface. In both cases, identical
reaction--diffusion parameters are used. All simulations are carried out using the IMEX-BDF1 time-stepping scheme with reaction parameters \( r_1 = 0.02 \), \( r_2 = 0.2 \). We evolve the system for 2000 time steps with a fixed time step size of \( \Delta t = 0.1 \), corresponding to a final time of \( t = 200 \). As shown in Figure~\ref{fig:Turing system_combined}, each row represents a different surface geometry (asymmetric torus, Stanford Bunny), while columns display the solution \( u \) at times \( t = 0 \), \( t = 20 \), and \( t = 200 \).  For both meshes, we employ~10\textsuperscript{th}-order elements to capture the solution. The emergence and arrangement of spots differ noticeably across geometries, emphasizing that surface shape, in addition to the system parameters, plays a critical role in the development of Turing patterns.  Unlike in one- or two-dimensional scenarios, where geometry is often simplified or ignored, solving reaction-diffusion equations on various surfaces with fixed parameters reveals distinct pattern transitions.

\begin{figure*}[htbp]
\centering

\begin{minipage}{0.06\textwidth}
\centering\textbf{\small $r=4$}
\end{minipage}%
\begin{minipage}{0.90\textwidth}
\centering
\begin{subfigure}{0.33\textwidth}
  \includegraphics[width=\linewidth]{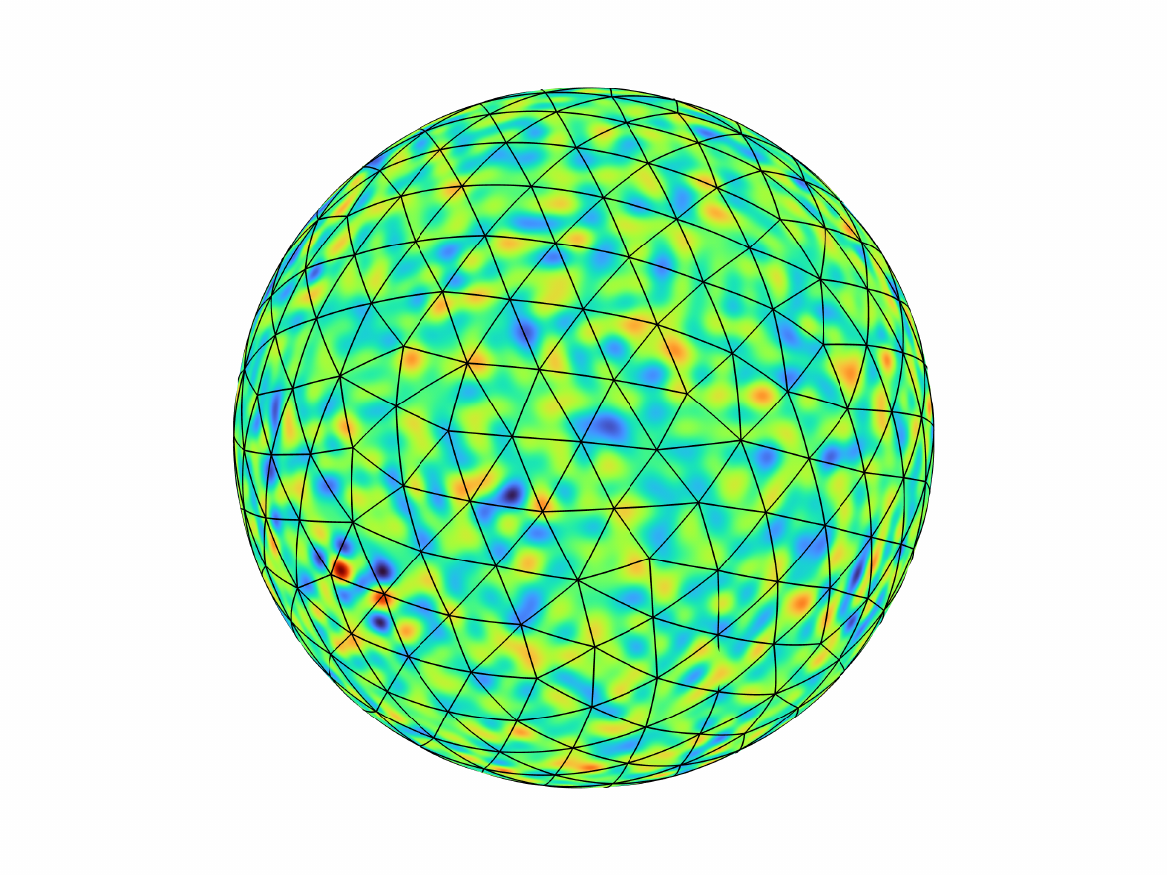}
\end{subfigure}\hfill
\begin{subfigure}{0.33\textwidth}
  \includegraphics[width=\linewidth]{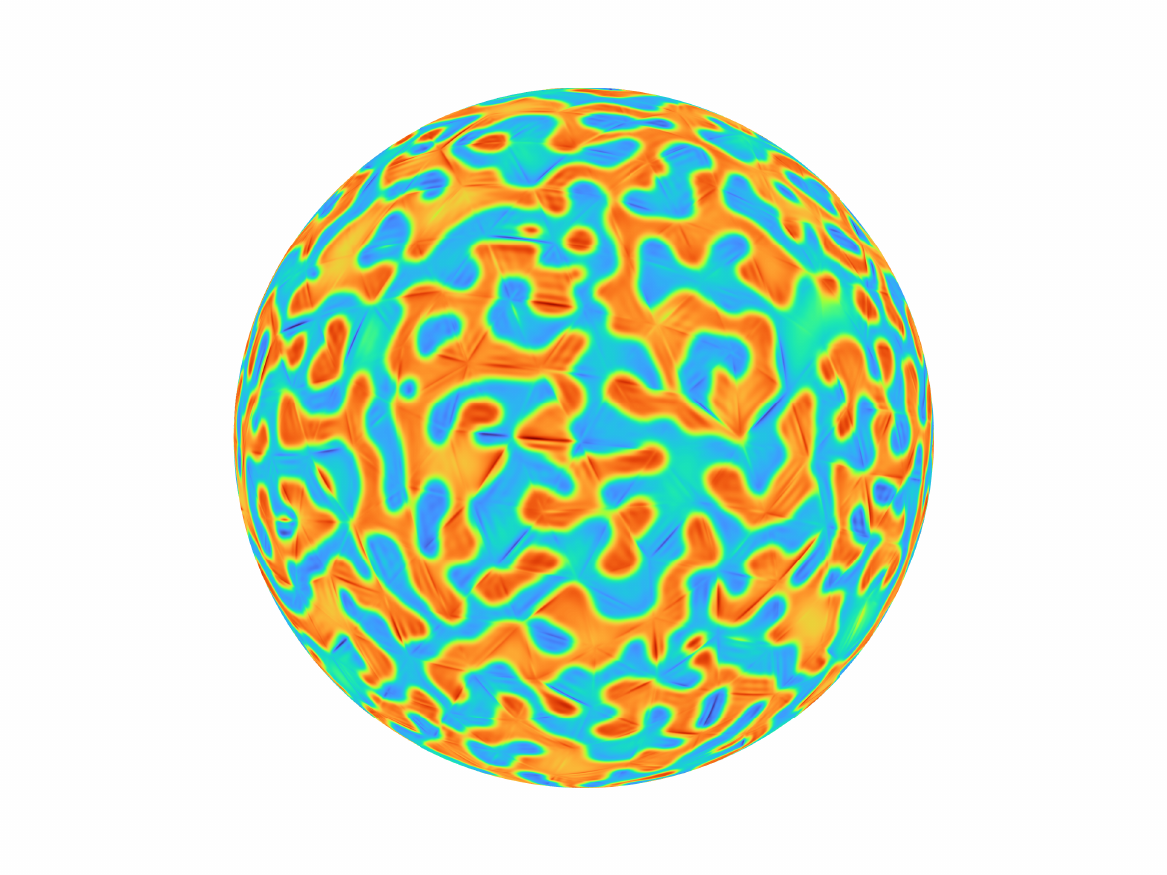}
\end{subfigure}\hfill
\begin{subfigure}{0.33\textwidth}
  \includegraphics[width=\linewidth]{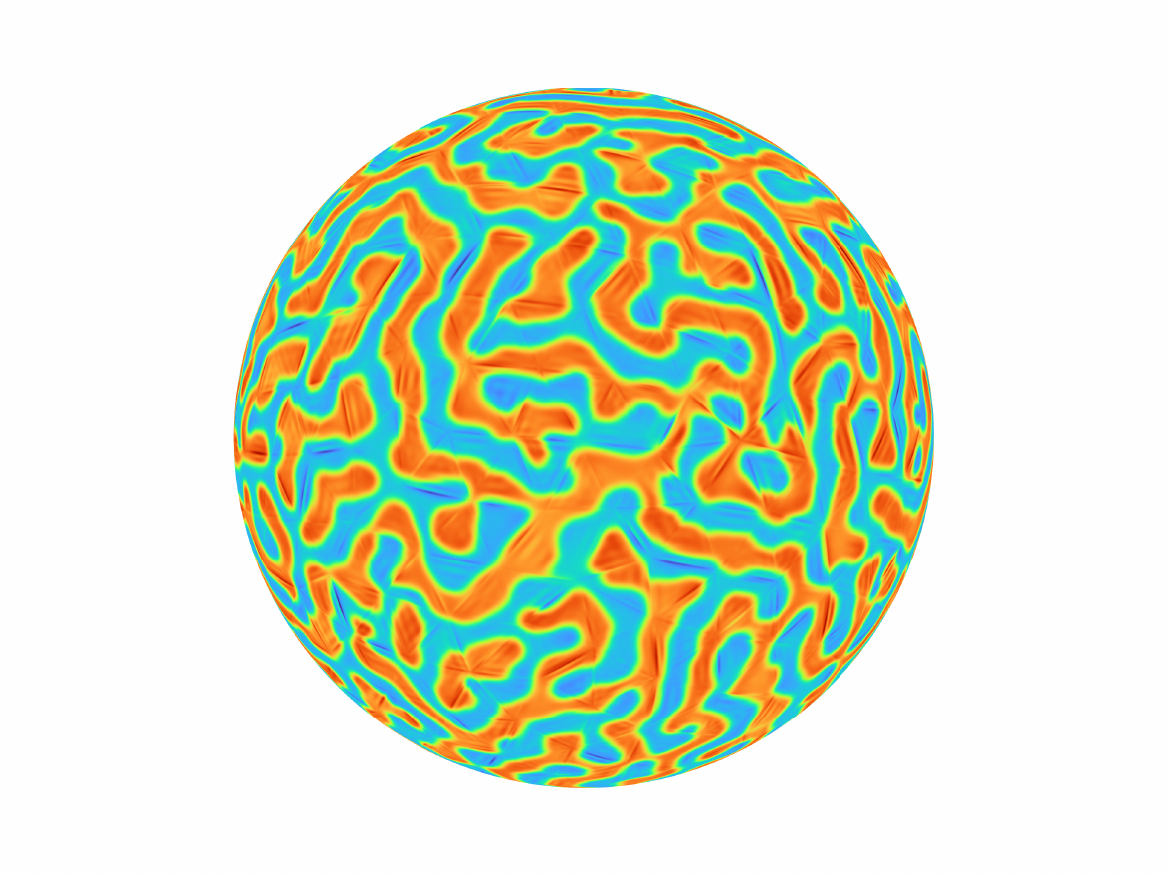}
\end{subfigure}
\end{minipage}

\vspace{1em}

\begin{minipage}{0.06\textwidth}
\centering\textbf{\small $r=2$}
\end{minipage}%
\begin{minipage}{0.90\textwidth}
\centering
\begin{subfigure}{0.33\textwidth}
  \includegraphics[width=\linewidth]{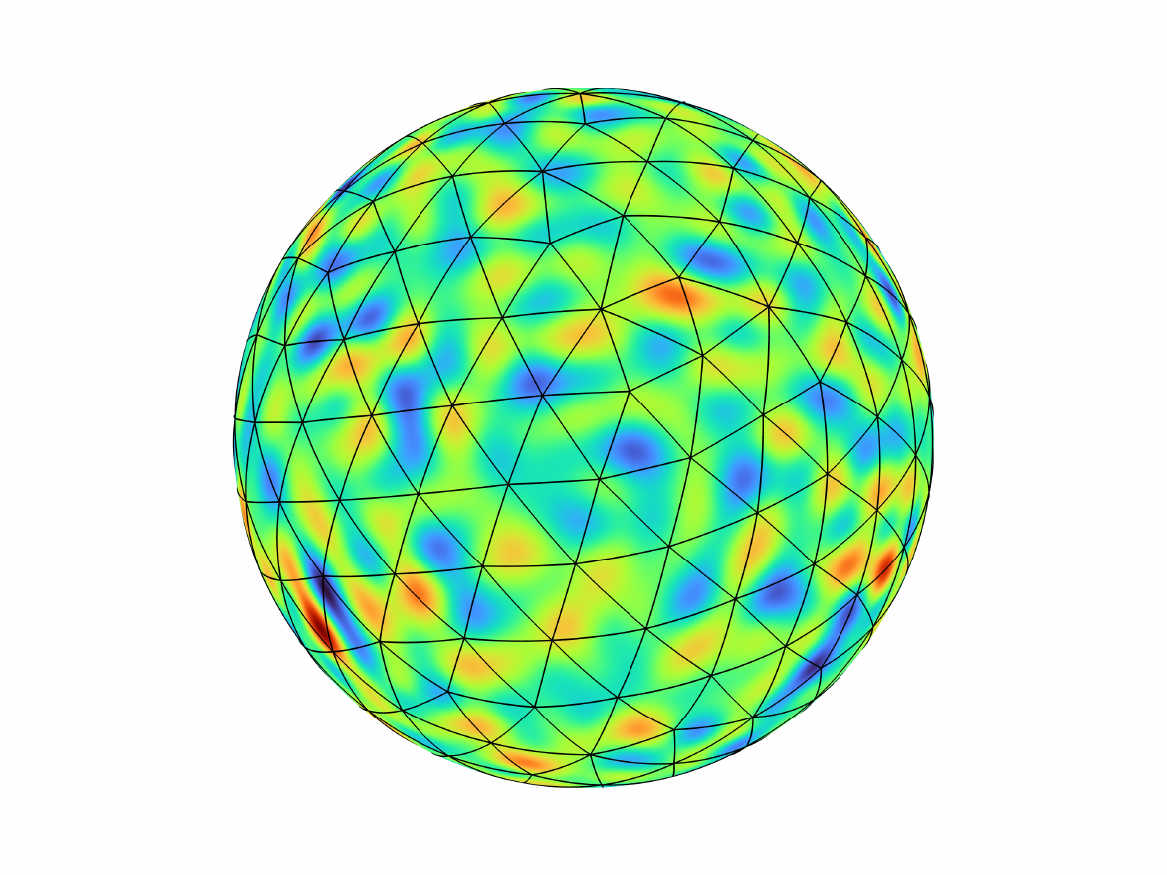}
\end{subfigure}\hfill
\begin{subfigure}{0.33\textwidth}
  \includegraphics[width=\linewidth]{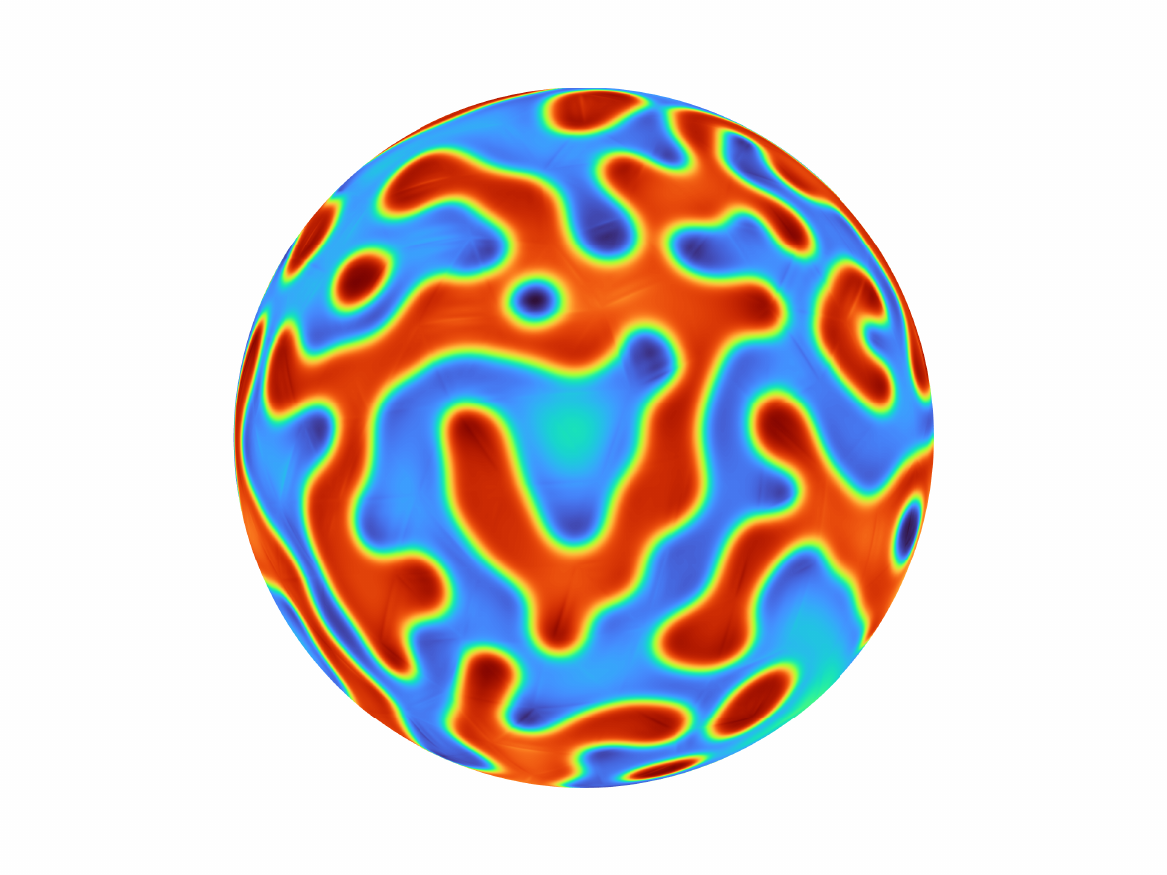}
\end{subfigure}\hfill
\begin{subfigure}{0.33\textwidth}
  \includegraphics[width=\linewidth]{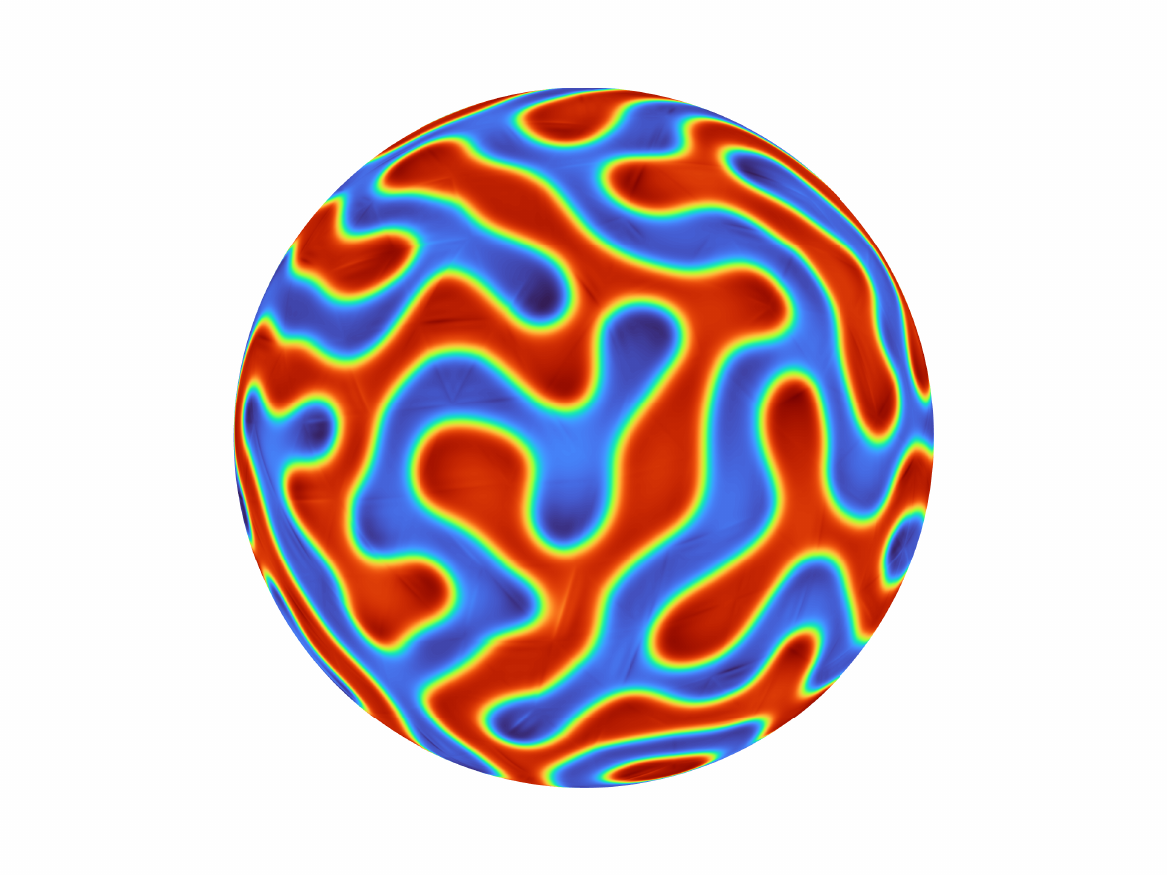}
\end{subfigure}
\end{minipage}

\vspace{1em}

\begin{minipage}{0.06\textwidth}
\centering\textbf{\small $r=1$}
\end{minipage}%
\begin{minipage}{0.90\textwidth}
\centering
\begin{subfigure}{0.33\textwidth}
  \includegraphics[width=\linewidth]{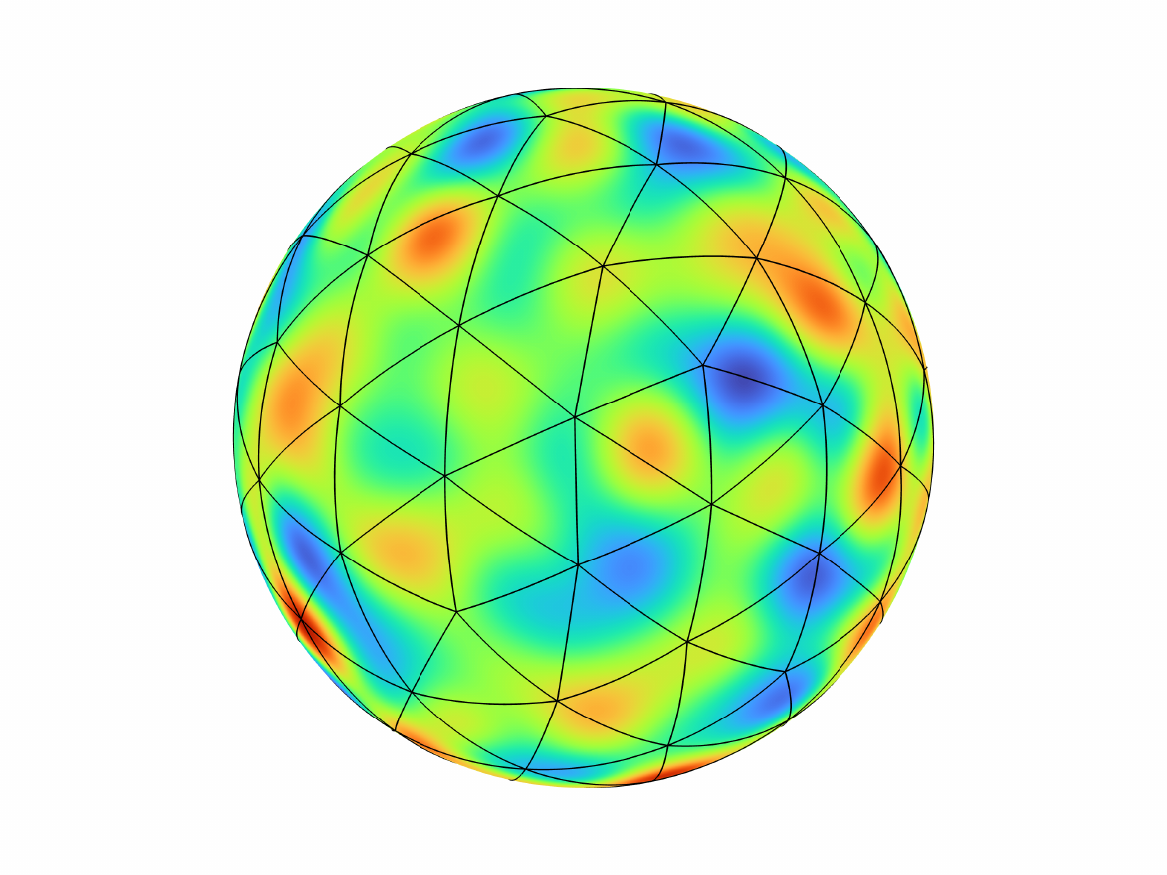}
\end{subfigure}\hfill
\begin{subfigure}{0.33\textwidth}
  \includegraphics[width=\linewidth]{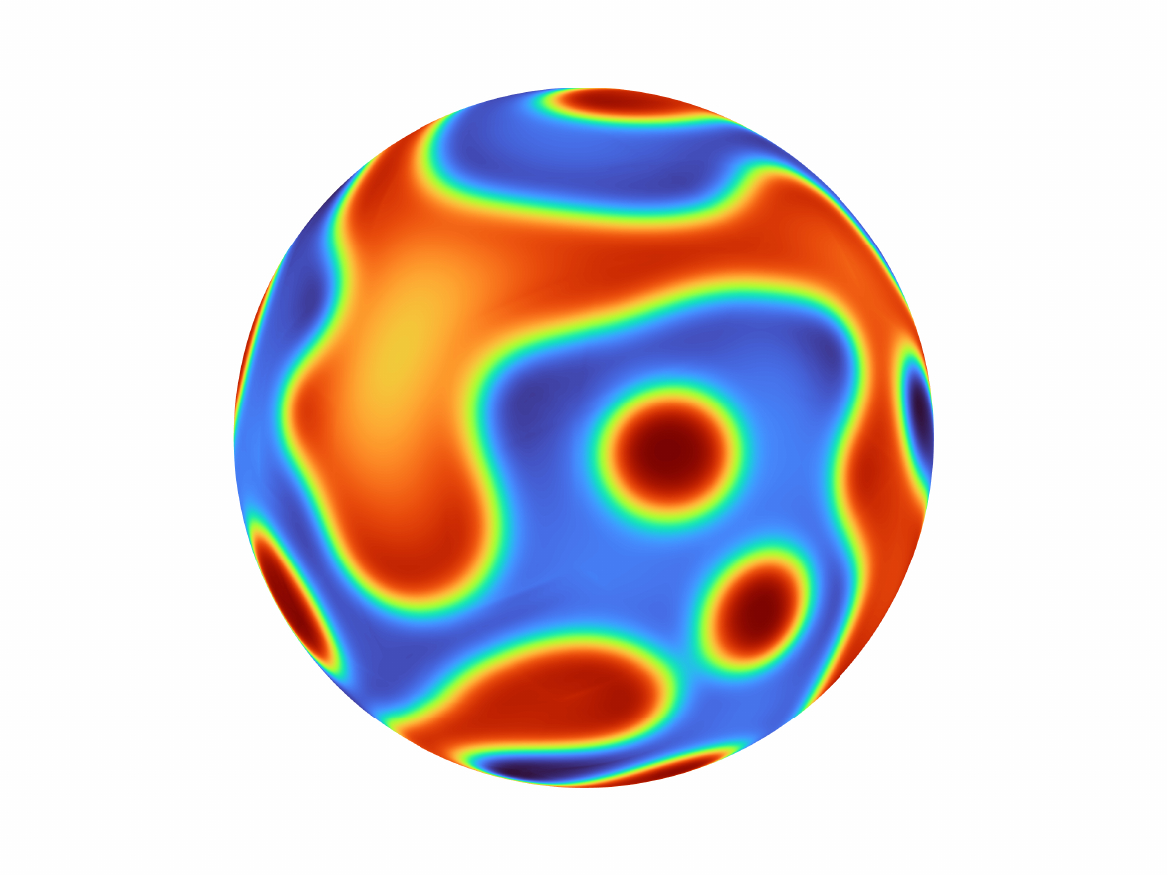}
\end{subfigure}\hfill
\begin{subfigure}{0.33\textwidth}
  \includegraphics[width=\linewidth]{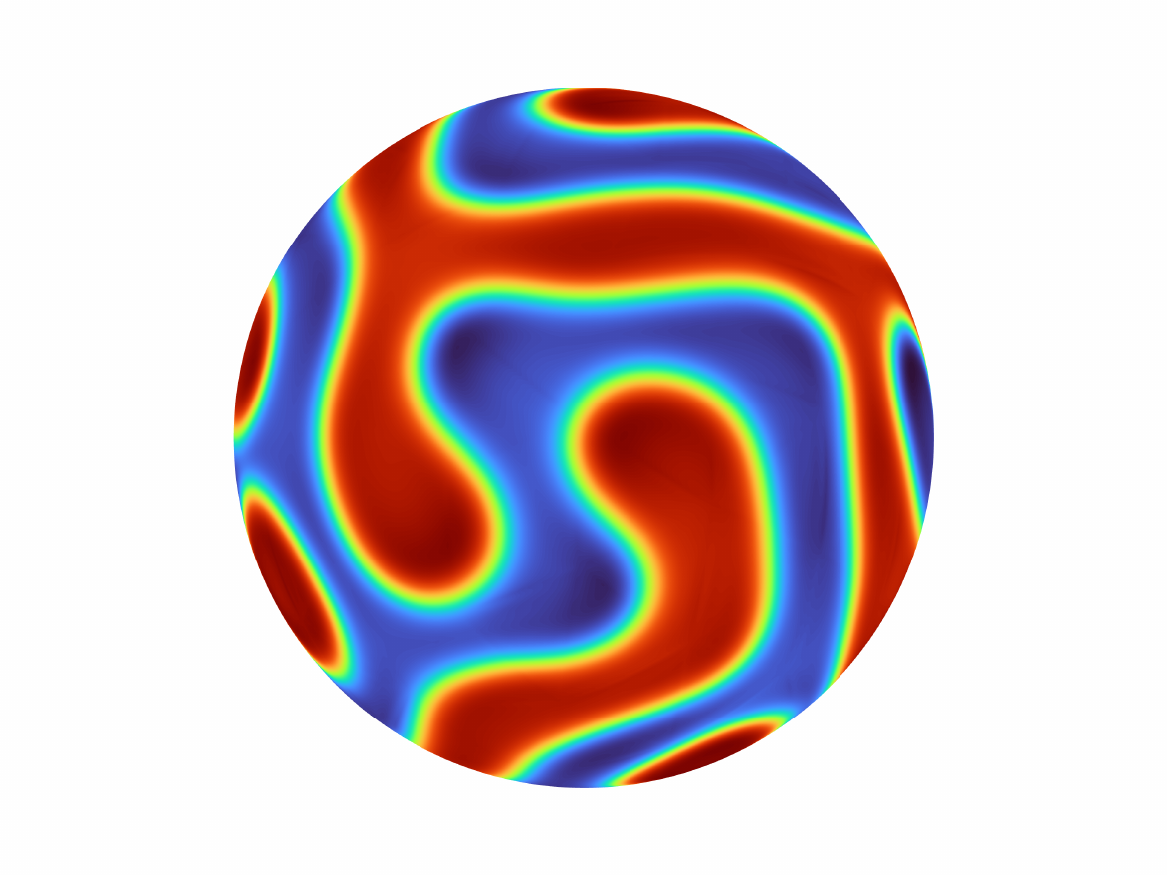}
\end{subfigure}
\end{minipage}

\caption{Pattern formation process of the reaction--diffusion model on the surface of a sphere corresponding to the concentration of the activator $u$, simulated at times \( t=0,20,200 \) using the IMEX--BDF1 scheme for different sphere radii.}
\label{fig:zebra_strip}
\end{figure*}

 In Figure~\ref{fig:zebra_strip}, we show simulations on spherical meshes with radii \( r = 1, 2, 4 \), using 10\textsuperscript{th}-order elements and random initial data. The parameters are chosen to promote stripe formation, with \( r_1 = 1.5 \), \( r_2 = 0 \), \( \alpha = 1.899 \), \( \gamma = -\alpha \), and \( \beta = -0.95 \), and diffusion coefficients set as \( \delta_{u_1} = 0.516\, \delta_{u_2} \), where \( \delta_{u_2} = 5 \cdot 10^{-3} \). Simulations are run up to \( t = 200 \) with a time step of \( \Delta t = 0.1 \). The results indicate that, as the domain expands, the number of high-activator stripes increases consistently, confirming theoretical and numerical predictions for two-dimensional growing domains ~\cite{jeong2017numerical}.
We close this section by investigating interacting Turing systems. Specifically, we consider a second Turing system in chemicals \( (u_{1}, u_{2}) \) that modulates the kinetic terms in the \( (v_{1}, v_{2}) \) system to give the model:
\begin{equation}\label{eq:interaction}
\begin{aligned}
\frac{\partial v_1}{\partial t} &= \delta_{v_1} \Delta_{\Gamma} v_1 + \alpha' v_1 (1 - r_1 v_2^2) + v_2(1 - r_2 v_1) + q_1 u_1 + q_2 u_1 v_2 + q_3 u_1 v_2^2 \\
\frac{\partial v_2}{\partial t} &= \delta_{v_2} \Delta_{\Gamma} v_2 + \beta' v_2 \left( 1 + \frac{\alpha' r_1}{\beta'} v_1 v_2 \right) + v_1(\gamma' + r_2 v_2) - q_2 u_2 v_1 - q_3 u_2^2 v_1 \\
\frac{\partial u_1}{\partial t} &= \delta_{u_1} \Delta_{\Gamma} u_1 + \alpha u_1(1 - r_1 v_2^2) + u_2(1 - r_2 u_1) \\
\frac{\partial u_2}{\partial t} &= \delta_{u_2} \Delta_{\Gamma} u_2 + \beta u_2 \left( 1 + \frac{\alpha r_1}{\beta} u_1 u_2 \right) + u_1(\gamma + r_2 u_2).
\end{aligned}
\end{equation}
This setup allows us to investigate how pattern properties change due to the interaction between the two coupled systems. For the simulations, we use the following parameter values:

\begin{itemize}
  \item For the \( (v_1, v_2) \) system:
  \[
  \alpha' = 0.398, \quad \beta' = -0.41, \quad \gamma' = -\alpha', \quad 
  \delta_{v_2} = 5 \times 10^{-3}, \quad \delta_{v_1} = 0.122\, \delta_{v_2}.
  \]
  
  \item For the \( (u_1, u_2) \) system:
  \[
  \alpha = 0.899, \quad \beta = -0.91, \quad \gamma = -\alpha, \quad 
  \delta_{u_2} =  \delta_{v_2}, \quad \delta_{u_1} = 0.516\, \delta_{u_2}.
  \]
\end{itemize}

\begin{figure*}[!t]
  \centering

  \begin{subfigure}{0.24\textwidth}
    \centering
    \includegraphics[width=\linewidth]{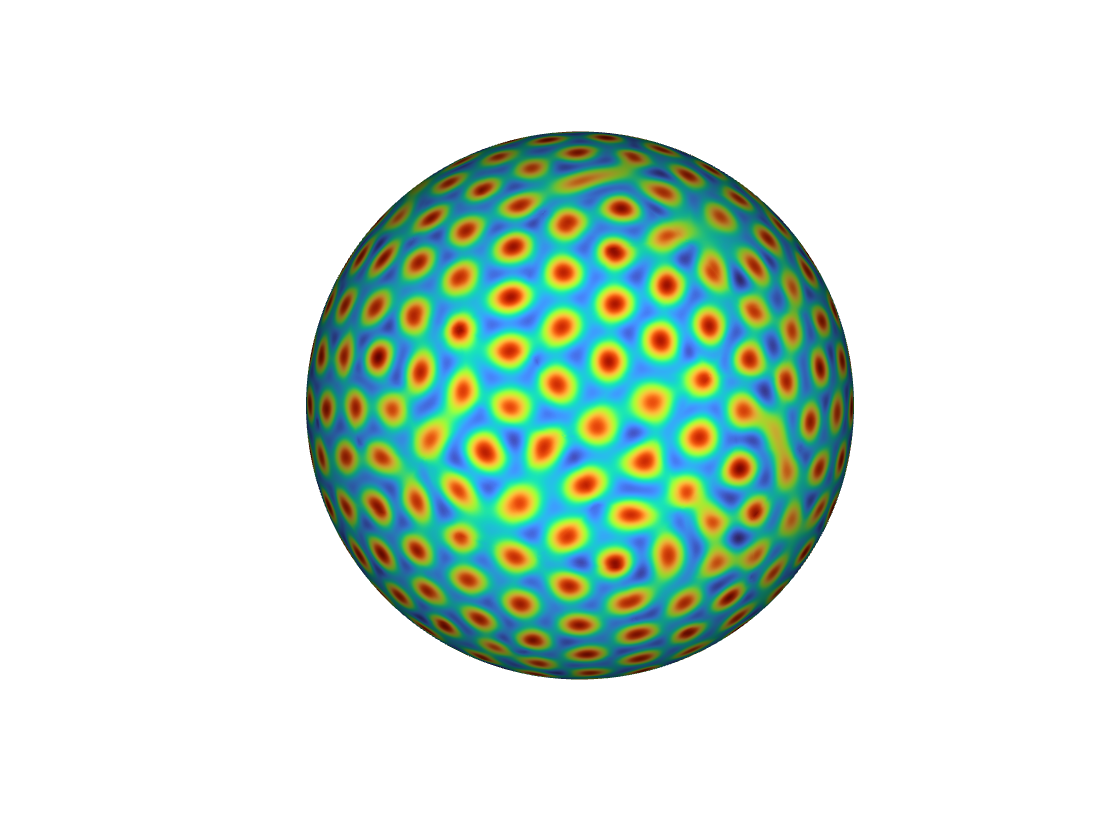}
    \caption{ $u_1$}
  \end{subfigure}\hfill
  \begin{subfigure}{0.24\textwidth}
    \centering
    \includegraphics[width=\linewidth]{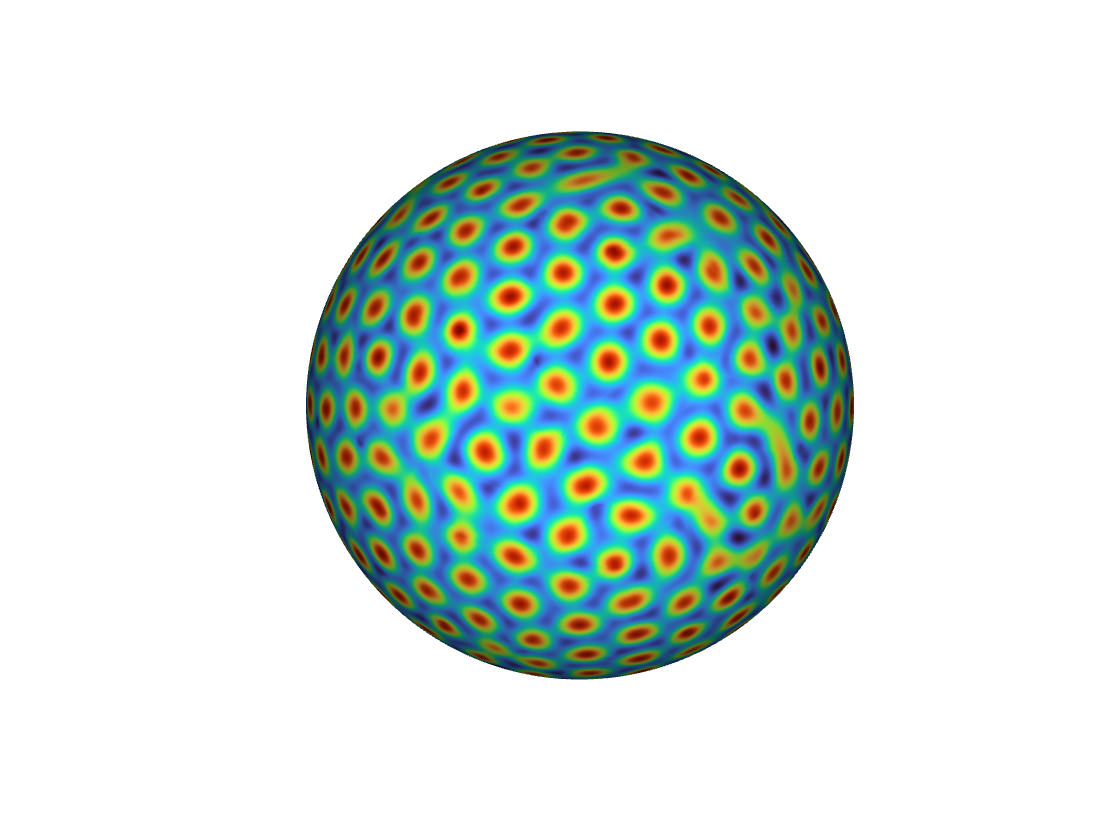}
    \caption{$v_1$, with $q_1 \neq 0$}
  \end{subfigure}\hfill
  \begin{subfigure}{0.24\textwidth}
    \centering
    \includegraphics[width=\linewidth]{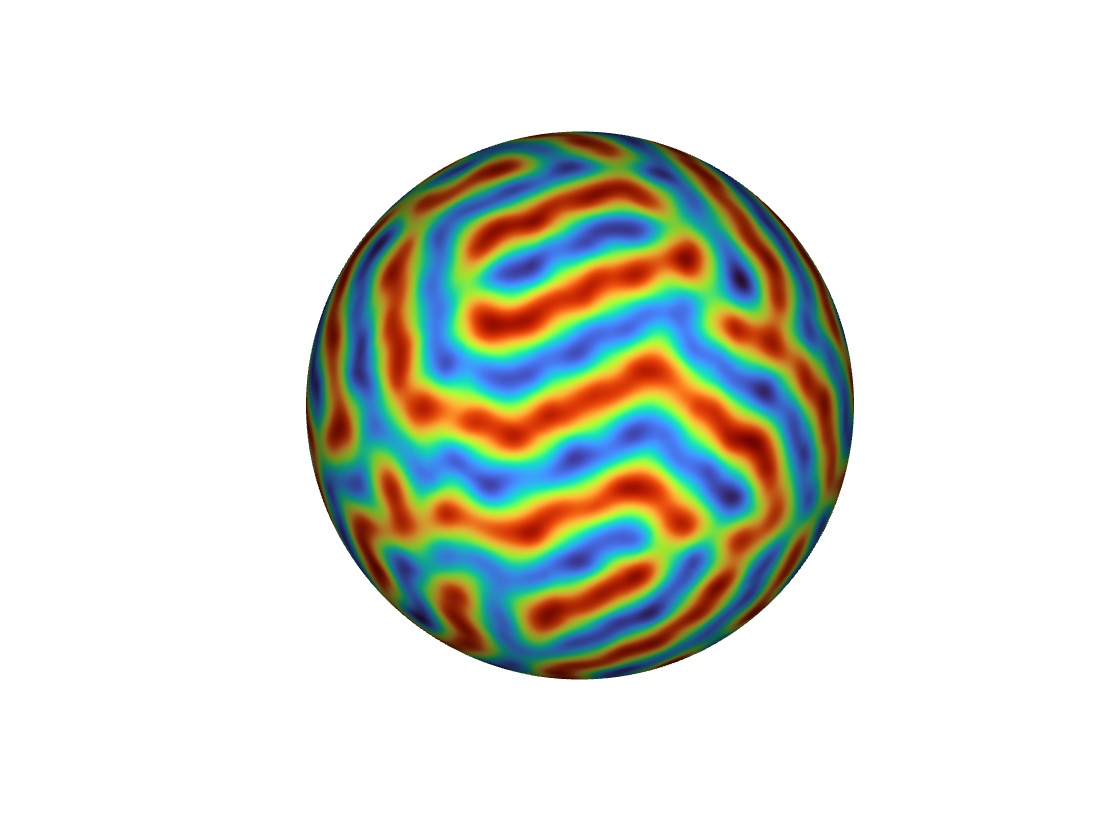}
    \caption{$v_1$, with $q_2 \neq 0$}
  \end{subfigure}\hfill
  \begin{subfigure}{0.24\textwidth}
    \centering
    \includegraphics[width=\linewidth]{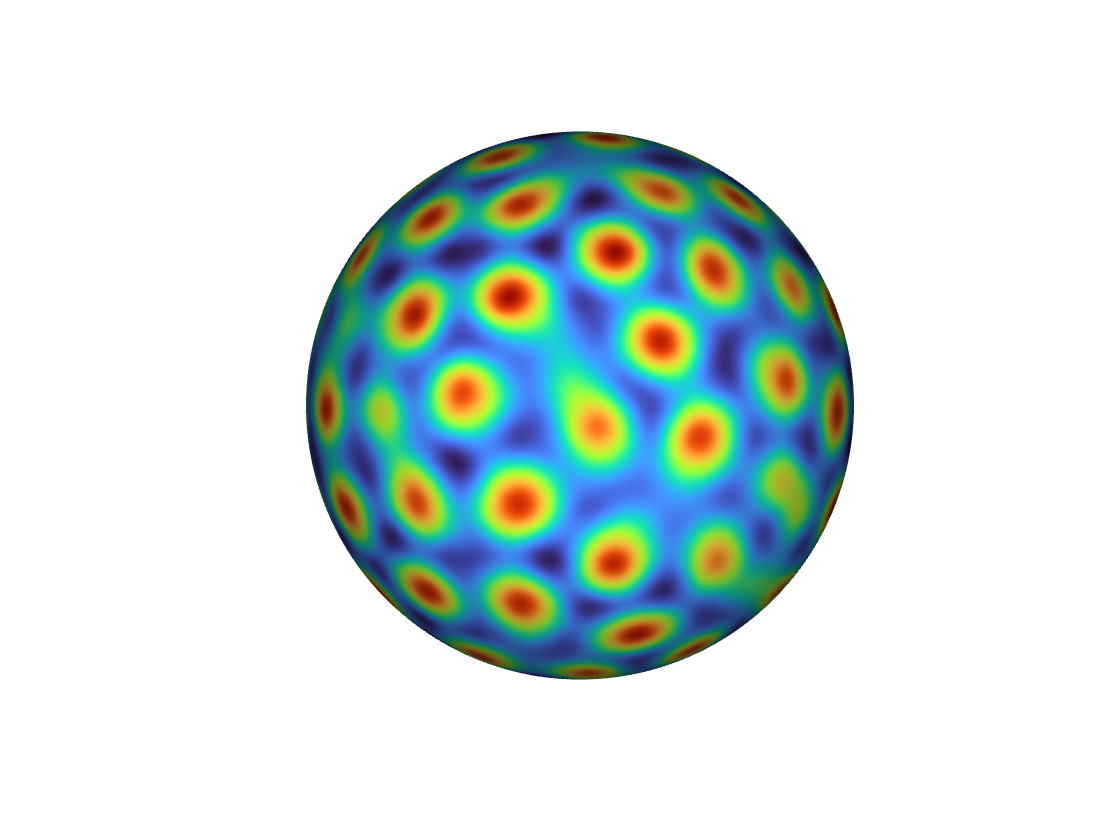}
    \caption{$v_1$, with $q_3 \neq 0$}
  \end{subfigure}

    \caption{
    Patterns obtained from the coupled system~\eqref{eq:interaction}.
    Each panel isolates a single coupling mechanism:
    linear ($q_1$), quadratic ($q_2$), or cubic ($q_3$),
    with coupling strength $0.55$.
  }
 \label{fig:coupled_patterns}
\end{figure*}

\begin{figure*}[!t]
  \centering

  \begin{subfigure}{0.24\textwidth}
    \centering
    \includegraphics[width=\linewidth]{images/up_1_q.png}
    \caption{$u_1$}
  \end{subfigure}\hfill
  \begin{subfigure}{0.24\textwidth}
    \centering
    \includegraphics[width=\linewidth]{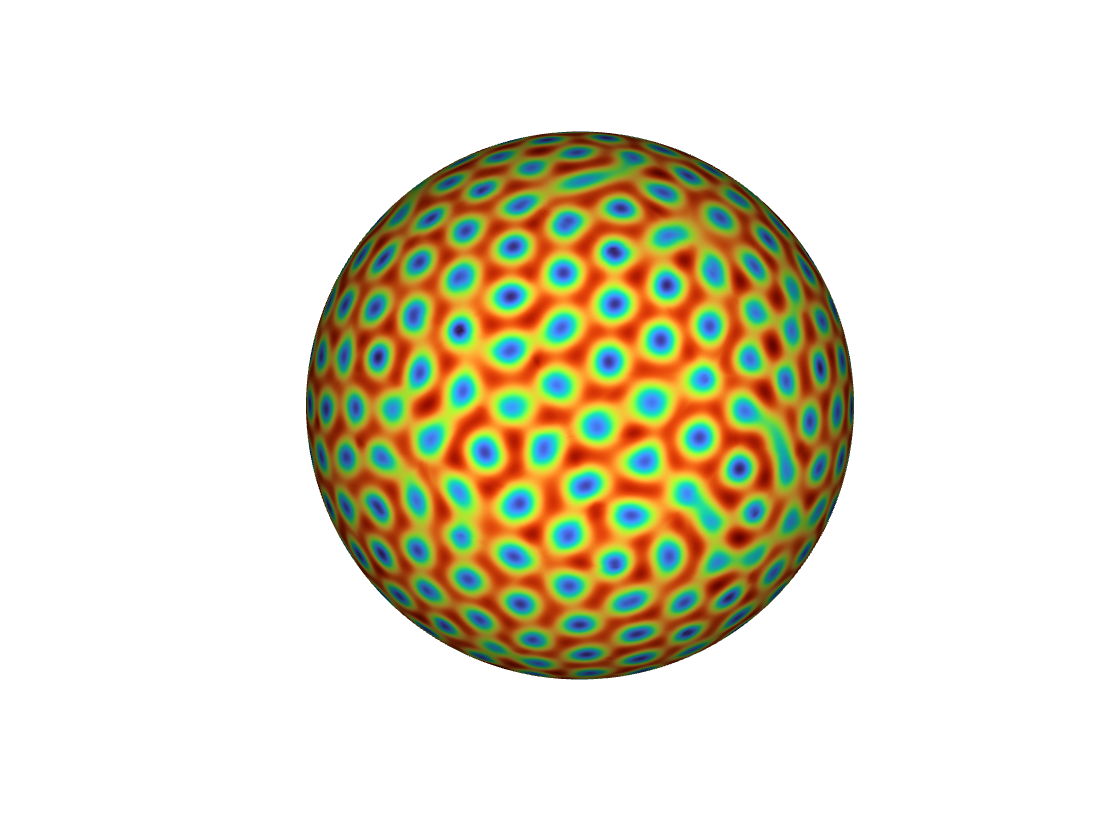}
    \caption{$v_1$, with $q_1 \neq 0$}
  \end{subfigure}\hfill
  \begin{subfigure}{0.24\textwidth}
    \centering
    \includegraphics[width=\linewidth]{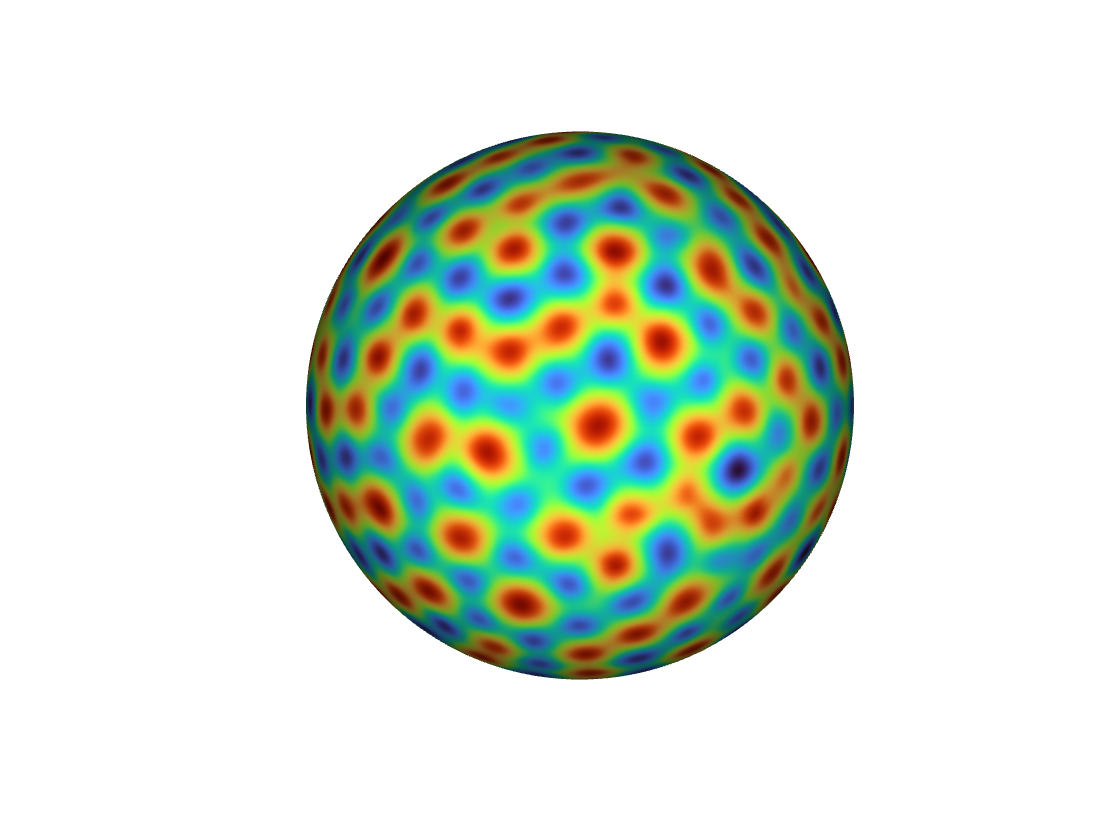}
    \caption{$v_1$, with $q_2 \neq 0$}
  \end{subfigure}\hfill
  \begin{subfigure}{0.24\textwidth}
    \centering
    \includegraphics[width=\linewidth]{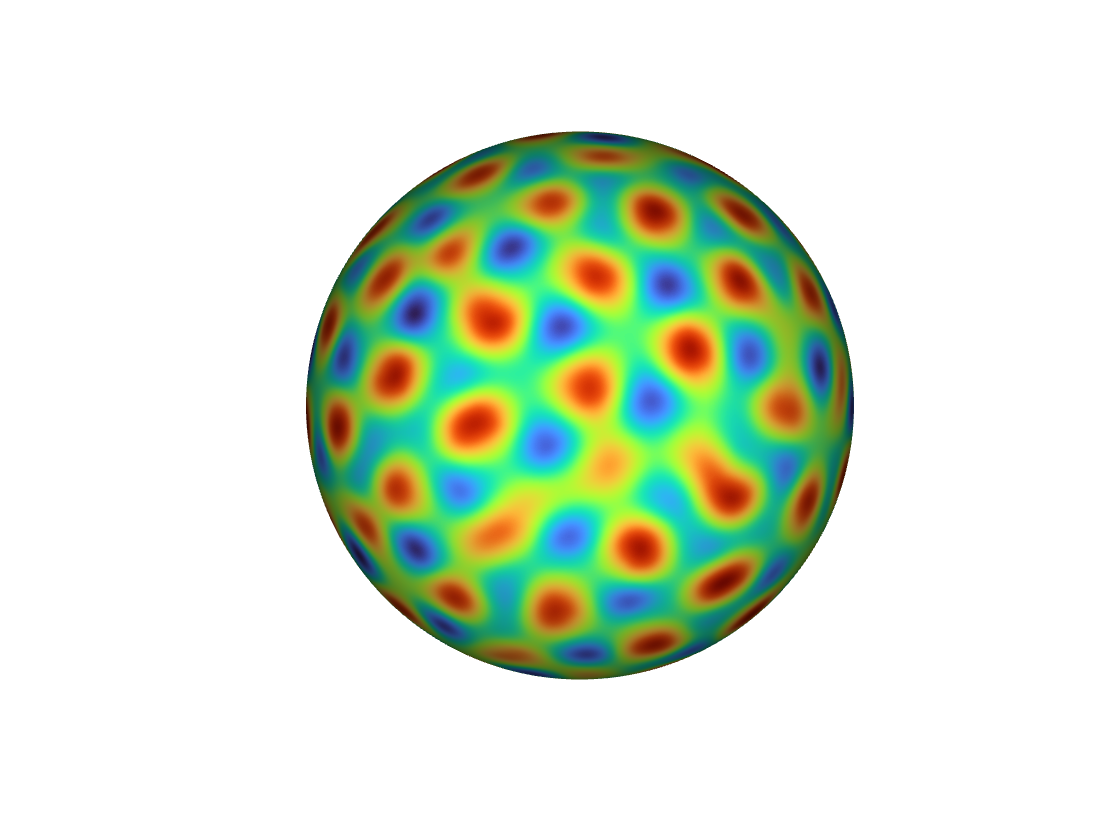}
    \caption{$v_1$, with $q_3 \neq 0$}
  \end{subfigure}

    \caption{
    Patterns obtained from the coupled system~\eqref{eq:interaction},
    using the same parameters as in Fig.~\ref{fig:coupled_patterns},
    but with negative coupling coefficients set to $-0.55$.
  }
\label{fig:-coupled_patterns}
\end{figure*}

In Figure~\ref{fig:coupled_patterns}, we show the resulting patterns for the cases of linear, quadratic, and cubic coupling. When only linear coupling is present, the pattern of \( v_1 \) becomes identical to that of \( u_1 \), indicating that the coupling completely overrides the dynamics of \( v_1 \). With cubic coupling, the solution still consists of spots, and the overall structure remains similar to the uncoupled case, suggesting that this type of interaction does not significantly alter the pattern. The most noticeable change occurs with quadratic coupling: the spot pattern is distorted and appears overlaid on a background of stripes, showing that the quadratic term introduces a strong modulation and leads to a mixed pattern with both spots and labyrinthine features. When the coupling coefficients are negative (see Figure~\ref{fig:-coupled_patterns}), the behavior changes significantly. The linear term produces a negative image of \( u_1 \) in \( v_1 \), effectively inverting the pattern. In the other two cases, the patterns in \( v_1 \) appear as a superposition of stripes and spots. Specifically, quadratic coupling now favors spot formation, while cubic coupling leads to more spatially ordered spots.

\section{Evolving surfaces}\label{sec:evol_surface}
The previous simulations on static surfaces show clearly how surface geometry affects the formation of reaction–diffusion patterns. Yet in many biological and physical systems, the surface itself changes over time, influencing the development and structure of these patterns. Accounting for surface evolution adds considerable complexity: the resulting models are non-autonomous, nonlinear parabolic PDEs, for which standard tools like linear stability analysis near bifurcations are no longer effective. This makes it essential to develop numerical methods that can handle such systems accurately and efficiently.
In response, we extend the triangle based HPS framework developed in the first
part of this work to evolving surfaces, where both the geometry and the solution
vary in time. Our approach follows a decoupled strategy:
 at each time step, we first update the surface geometry using a prescribed evolution law, then solve the PDE on the updated surface. A time-dependent surface \(\Gamma(t)\) is typically described using either an implicit or a Lagrangian representation. A widely used implicit approach defines \(\Gamma(t)\) as the zero level set of a time-dependent scalar function \(\varphi : \mathbb{R}^{d+1} \times [0, T] \to \mathbb{R}\), i.e.,
\begin{equation}\label{eq:dependet_surface} 
\Gamma(t) = \{ \mathbf{x} \in \mathbb{R}^{d+1} : \varphi(\mathbf{x},t) = 0 \}.
\end{equation}
In this formulation, the evolution of the surface is governed by the dynamics of the level set function \(\varphi\). Alternatively, a Lagrangian description tracks the motion of surface points explicitly through a velocity field \( \bm{v} : \mathbb{R}^{d+1} \times [0, T] \to \mathbb{R}^{d+1} \), which prescribes the flow driving the deformation of the surface over time.
 In this case, the surface \(\Gamma(t)\) is parameterized over its initial configuration \(\Gamma(0)\) by following the trajectories of points on the surface. For an initial point \(\mathbf{x}_0 \in \Gamma(0)\), the trajectory \(\mathbf{x}(t)\) satisfies the ordinary differential equation
\begin{equation}\label{eq:evolv_1}
    \frac{d}{dt}\mathbf{x}(t) = \bm{v}\left(\mathbf{x}(t), t \right), \quad \mathbf{x}(0) = \mathbf{x}_0.
\end{equation}

If Eq.~\eqref{eq:evolv_1} can be solved in closed form (i.e., yielding an explicit expression for \(\mathbf{x}(t)\)), then we define the parameterization of \(\Gamma(t)\) over the initial surface \(\Gamma(0)\) as:
    $\Gamma(t) := \{\mathbf{x}(t) \;|\; \mathbf{x}_0 \in \Gamma(0)\}\,$ .
If an analytical solution is not known, we employ the Explicit Euler method for simplicity to numerically solve the ODE system in Eq.~\eqref{eq:evolv_1}:
\begin{equation*}
\mathbf{x}_{m+1} = \mathbf{x}_{m} + \Delta t \, \bm{v}_{m}.
\end{equation*}

\subsection{Time dependent oriented distance function}\label{time_distance}
To study the numerical approximations of reaction-diffusion systems on evolving surfaces,  we incorporate the surface evolution law in the closest point projection, which results in a time-dependent closest point projection operator.

We assume that there exists an open bounded set \( U(t) \subset \mathbb{R}^{d+1} \) such that \( \partial U(t) = \Gamma(t) \). The \textit{time dependent oriented distance function} \( d_{\Gamma(t)}(\mathbf{x},t) \) is defined as
\[
d_{\Gamma(t)}(\mathbf{x},t):\;\mathbb{R}^{d+1} \times [0, T] \to \mathbb{R}, \quad
d_{\Gamma(t)}(\mathbf{x},t) := 
\begin{cases}
\operatorname{dist}(\mathbf{x}, \Gamma(t)) & x \in \mathbb{R}^{d+1} \setminus U(t), \\
-\operatorname{dist}(\mathbf{x}, \Gamma(t)) & \mathbf{x} \in U(t).
\end{cases}
\]

For \( \delta > 0 \) we define \( \mathscr{N}^t_\delta := \left\{ \mathbf{x} \in \mathbb{R}^{d+1} \,\middle|\, \operatorname{dist}(\mathbf{x}, \Gamma(t)) < \delta \right\} \). Clearly \( \mathscr{N}^t_\delta  \)  is an open neighbourhood of $\Gamma(t)$. For each \( \mathbf{x} \in \mathcal{N}_{\delta}^t \), there exists a unique time-dependent closest point projection operator \( \pi(\mathbf{x},t) \in \Gamma(t) \), defined by
\[
    \pi(\mathbf{x},t) = \mathbf{x} + d_{\Gamma(t)}(\mathbf{x},t) \nabla d_{\Gamma(t)}(\mathbf{x},t).
\]

By analogy with the stationary case \cite{gilbarg1977elliptic}, the existence of \( d_{\Gamma(t)} \) is readily established, and its regularity follows directly from the smoothness of the parametrization of the surface. To account for surface evolution over time, we compose the stationary closest point projection \( \pi \) with a Lipschitz continuous mapping \( \mathcal{L}_{t}: \Gamma(0) \to \Gamma(t) \), which tracks the deformation of surface points. The inverse \( \mathcal{L}_{t}^{-1} \) is also assumed to exist and to be Lipschitz continuous.

\begin{definition}
Given a Lipschitz continuous mapping \(\mathcal{L}_{t}: \Gamma(0) \rightarrow \Gamma(t), \;t\in[0,T]\), such that its inverse \(\mathcal{L}_{t}^{-1}: \Gamma(t) \rightarrow \Gamma(0)\) is Lipschitz as well. In other words, there exist constants \(c, C > 0\) such that 
\[
\|\mathcal{L}_{t}(\mathbf{x}_1) - \mathcal{L}_{t}(\mathbf{x}_2)\| \leq c\|\mathbf{x}_1 - \mathbf{x}_2\| \quad \text{and} \quad \|\mathcal{L}_{t}^{-1}(\mathbf{y}_1) - \mathcal{L}_{t}^{-1}(\mathbf{y}_2)\| \leq C\|\mathbf{y}_1 - \mathbf{y}_2\|.
\]
The time-dependent closest point projection operator \(\pi(\mathbf{x}(t),t)\) is defined by the composition of the closest point projection \(\pi\) with the Lipschitz continuous mapping \(\mathcal{L}_{t}\). That is, \(\pi: \mathcal{N}^{t}_{\delta} \rightarrow \Gamma(t)\) is
\begin{equation}\label{operator.projection}
\pi(\mathbf{x}(t),t) := (\mathcal{L}_{t} \circ \pi)(\mathbf{x}(0),t).
\end{equation}
\end{definition}

This composition ensures that the projection operator \( \pi(\mathbf{x}(t),t) \) accounts for both the geometry of the surface and its deformation in time. Note, the closest point projection on the right-hand side does not depend explicitly on time, while $\mathcal{L}_{t}$ is time-dependent. We apply the time-dependent projection operator in the numerical solution of PDEs on evolving surfaces and in the construction of triangulated surface meshes at each time step. In this context, the mapping \( \mathcal{L}_t \) is often referred to as an arbitrary Lagrangian–Eulerian (ALE) map~\cite{kovacs2019computing}. 


\section{Evolving discrete surfaces}\label{sec:surf_app}

Given a smooth initial surface $\Gamma(0)$, approximated by a triangulated surface $\Gamma^{\text{tri}}_h(0)$, i.e., a quasi-uniform family of triangulations $\widehat{\mathcal{T}}_h(0)$ of maximal element diameter $h$. Let $\mathbf{x}_k(0)$, $(k = 1, 2, \dots, N)$ denote the points of $\Gamma^{\text{tri}}_h(0)$ lying on the initial smooth surface $\Gamma(0)$. The points will be evolved in time with the given normal velocity $\bm{v}$, by solving the ODE
\begin{equation}\label{ode_surface_evolv}
\frac{d}{dt} \mathbf{x}_k(t) = \bm{v}(\mathbf{x}_k(t),t) \quad (k = 1, 2, \dots, N),
\end{equation}
\noindent
where $\bm{v}$ represents the velocity field. Obviously, the points remain on the surface $\Gamma(t)$ for all times, i.e., $d_{\Gamma(t)}(\mathbf{x}_k(t),t) = 0$ for $k = 1, 2, \dots, N$ and for all $t \in [0, T]$.


Lets \( \Delta_d \subset \mathbb{R}^{d} \) be the standard simplex of dimension $d$ and let \([0,T]\) be a time interval. A linear evolving mesh or evolving discrete surface on \([0,T]\) consists of continuously time-dependent points \(\mathbf{x}_k(t)_{k=1}^N \subset \mathbb{R}^{d+1}\) and $d$-dimensional simplex relations $\widehat{T}(t) \in \widehat{\mathcal{T}}_h(t)$, where we identify \( \widehat{T}(t)\subset \mathbb{R}^{d+1} \) with the simplex itself, which we require to satisfy:

\begin{itemize}
    \item \(\widehat{T}(t) \in \widehat{\mathcal{T}}_h(t) \) is non-degenerate, i.e., the map \( \varrho: \Delta_d \times [0,T] \to \widehat{T}(t) \),
    is a bijection.

    \item The intersection of two simplices is a common edge, a common point, or empty.

    \item There are no boundary simplices, i.e., every edge is the intersection of two different simplices.
\end{itemize}

We set the mesh width \( h \) as
\[
h(t) := \max_{\widehat{T}(t) \in \widehat{\mathcal{T}}_h(t)} \operatorname{diam}(\widehat{T}(t)), \quad h := \sup_{t \in [0,T]} h(t),
\]
where \( \operatorname{diam} \) denotes the 2-dimensional diameter, the \textit{in-ball radius} at time \( t \) is defined as
\[
    r(t) := \min_{\widehat{T}(t) \in \widehat{\mathcal{T}}_h(t)} \rho(\widehat{T}(t)),
\]
where \( \rho \) denotes the radius of the maximum inner circle. We set
\[
    \mathbf{x}_h(t) := (\mathbf{x}_k(t))_{k=1}^N \in \mathbb{R}^N \otimes \mathbb{R}^{d+1}
\]
and
\[
    \Gamma(\mathbf{x}_h) := \Gamma(\mathbf{x}_h(t)) := \Gamma^{\text{tri}}_h(t) := \bigcup_{\widehat{T}(t) \in \widehat{\mathcal{T}}_h(t)} \widehat{T}(t).
\]

\begin{figure}[ht]
    \centering
    \begin{subfigure}[b]{0.3\textwidth}
        \centering
        \includegraphics[width=1\textwidth]{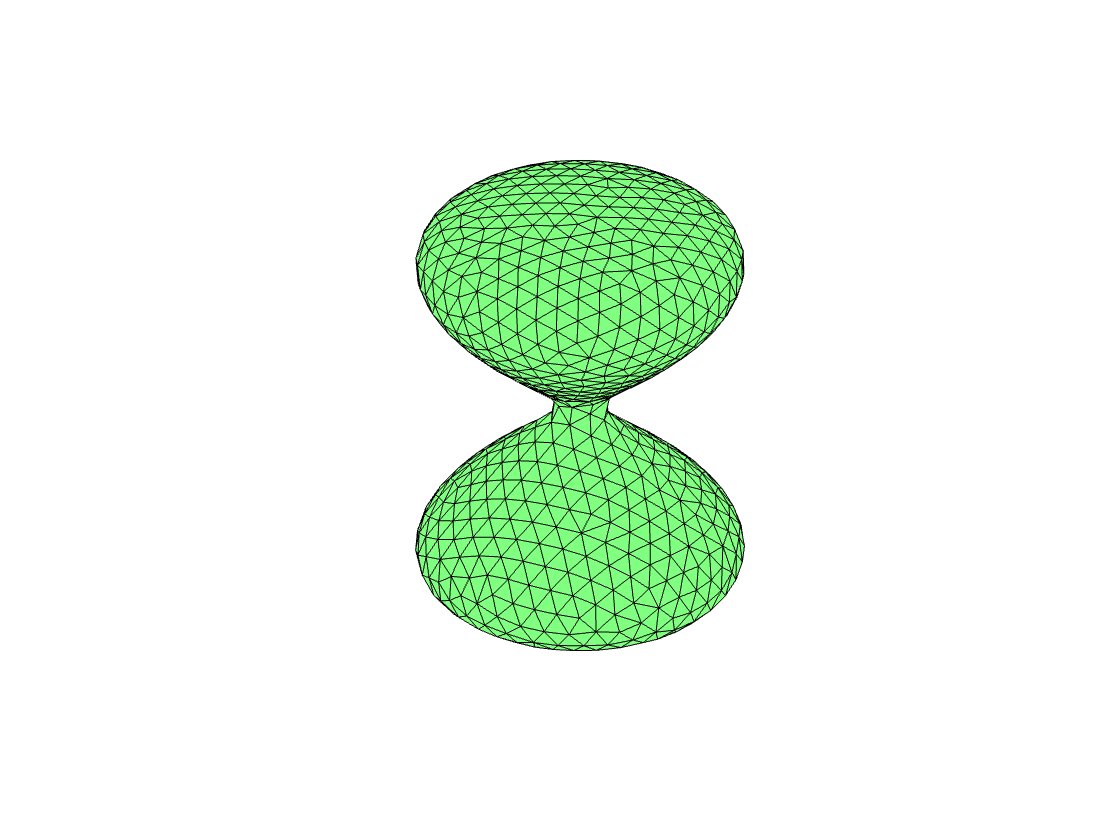}
        \label{fig:t0}
    \end{subfigure}
    \hfill
    \begin{subfigure}[b]{0.3\textwidth}
        \centering
        \includegraphics[width=1\textwidth]{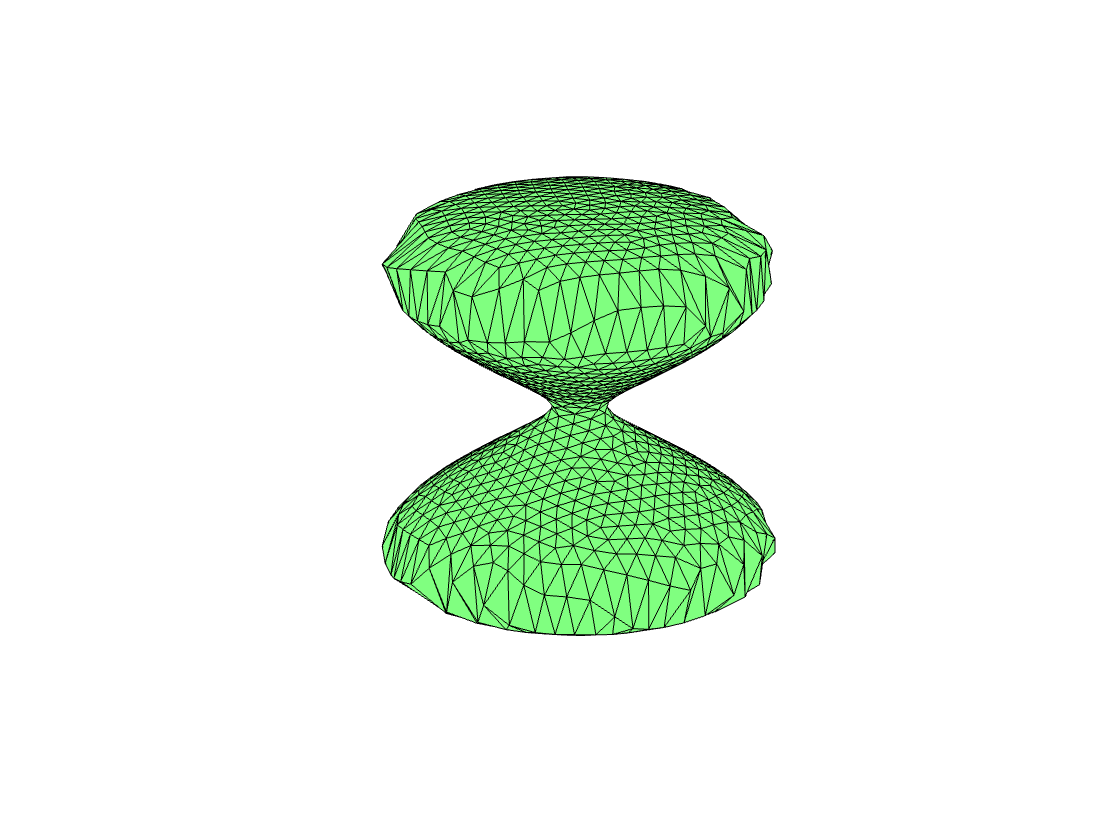}
        \label{fig:t02}
    \end{subfigure}
 \hfill
    \begin{subfigure}[b]{0.3\textwidth}
        \centering
        \includegraphics[width=1\textwidth]{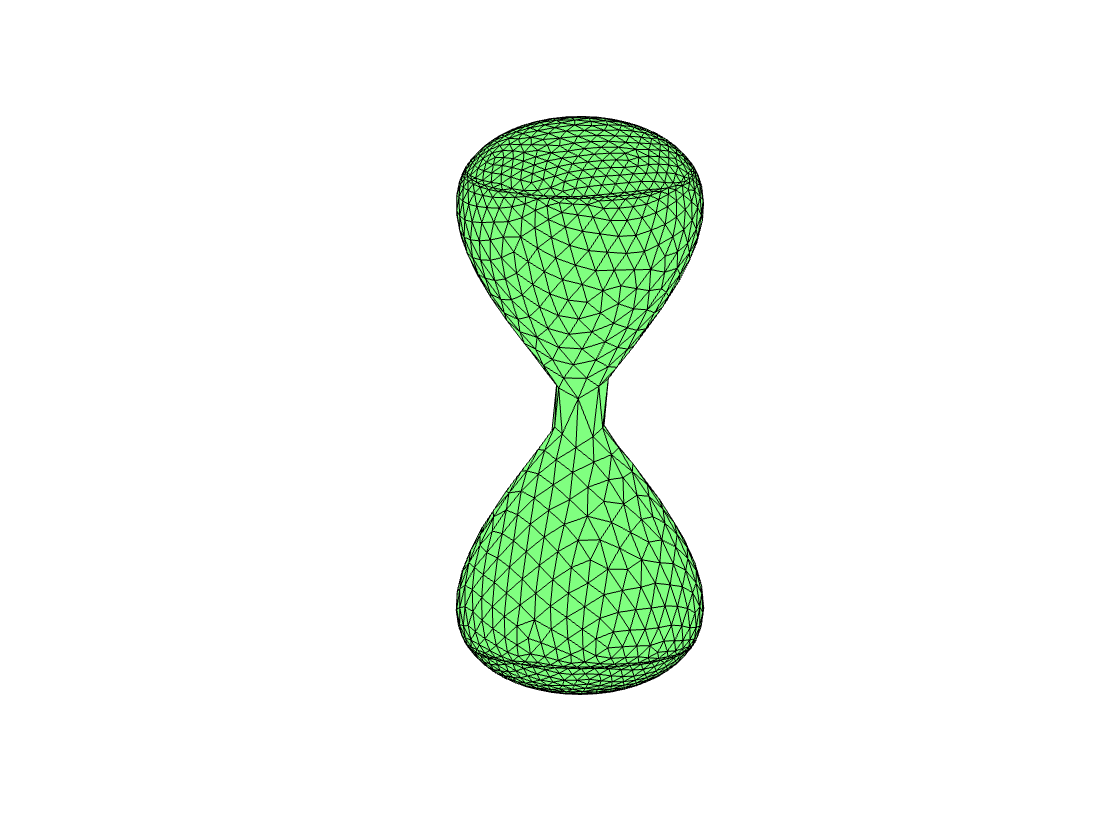}
        \label{fig:t06}
    \end{subfigure}
    \caption{Normal evolution of a closed surface at time $t = 0, 0.2, 0.9$.}
    \label{fig:surface_evolution}
\end{figure}

A family \( (\Gamma^{\text{tri}}_h(t))_{h>0} \) is called \textit{quasi-uniform}, if and only if
\[
    \sup_{t \in [0,T]}\frac{h(t)}{r(t)} 
 \leq c.
\]
Quasi-uniformity is generally not preserved during surface evolution. As shown in
Figure~\ref{fig:surface_evolution}, evolving a surface over \([0,0.9]\) leads to significant mesh degradation,
despite initially mild deformation. Normal evolution causes severe element
distortion, with triangles collapsing into near-zero angles.

\begin{figure}[ht]
    \centering
    \begin{subfigure}[b]{0.3\textwidth}
        \centering
        \includegraphics[width=1\textwidth]{images/t_0.png}
    \end{subfigure}
    \hfill
    \begin{subfigure}[b]{0.3\textwidth}
        \centering
        \includegraphics[width=1\textwidth]{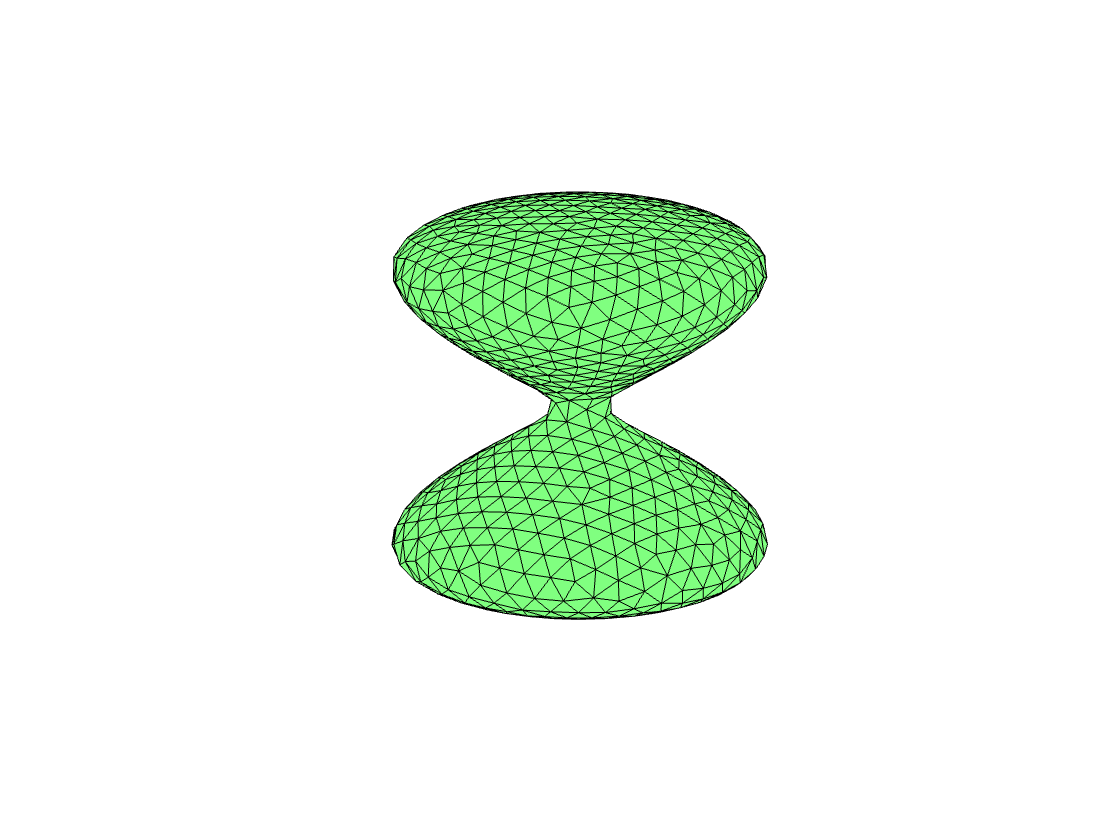}
    \end{subfigure}
    \hfill
    \begin{subfigure}[b]{0.3\textwidth}
        \centering
        \includegraphics[width=1\textwidth]{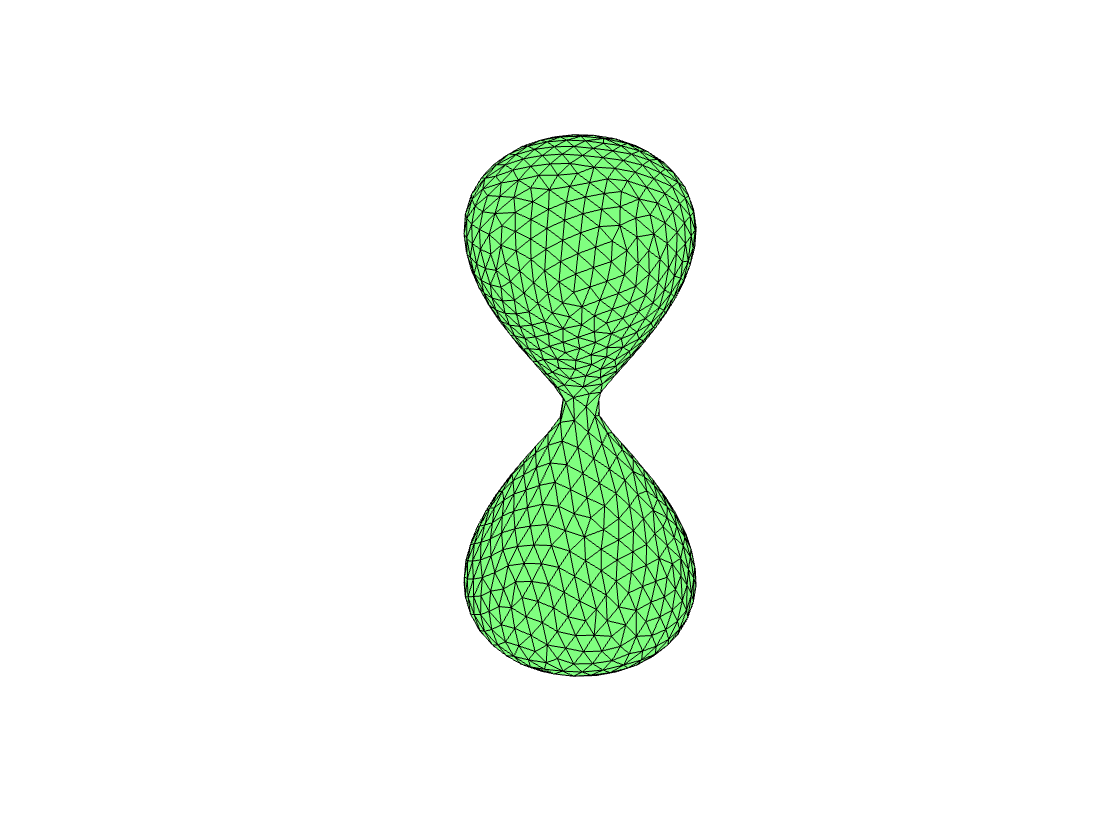}
    \end{subfigure}
    \caption{Normal evolution of a closed surface using the ALE map~\cite{kovacs2019computing} at time \( t = 0, 0.2, 0.9 \).}
    \label{fig:surface_evolution_alemap}
\end{figure}

While remeshing techniques~\cite{geuzaine2009gmsh,Persson} can restore mesh quality,
they typically require introducing or removing nodes and constructing mappings
between successive meshes, making them computationally expensive. We therefore adopt the arbitrary Lagrangian–Eulerian (ALE) approach of~\cite{kovacs2019computing}, which preserves mesh connectivity during surface evolution. As illustrated in Figure~\ref{fig:surface_evolution_alemap}, ALE
effectively prevents mesh distortion and maintains element regularity. To further
reduce computational cost, all meshes are precomputed and reused.

\section{Surface growth}\label{sec:surface_growth}

We focus on two growth models for the evolution of the surface.

\subsection{Isotropic growth}
Here we assume that each material point on the surface undergoes the same prescribed dilation. Thus, if \( ((x_1(t),x_2(t),x_3(t)) \) denotes a material point at time \( t \), then
\[
(x_1(t),x_2(t),x_3(t)) = \eta(t)(x_1(0),x_2(0),x_3(0))
\]
where the dilation factor \( \eta(t) \) is a positive growth function satisfying the property \( \eta(0) = 1 \). It's worth noting that such surface growth generally involves both normal and tangential velocities.

\subsection{Anisotropic growth}
Here we are interested in a prescribed velocity law associated with non-isotropic growth. As an illustrative example, consider the time-dependent surface \( \Gamma(t) \) given by the zero level set of the distance function
\[
d_{\Gamma(t)}(\mathbf{x}, t) = a(t)^2 b\left( \frac{x_1^2}{c(t)^2} \right) + x_2^2 + x_3^2 - a(t)^2\,.
\]
 Different surfaces may be obtained by prescribing the non-negative shape functions $a(\cdot)$, $b(\cdot)$, and $c(\cdot)$ with $b(0) = 0$ and $b(1) = 1$.

\begin{itemize}
    \item \textbf{Tangential motion} \\
For the described surface, material points can move both tangentially and normally. We represent the surface $\Gamma(t)$ using material points $(\tilde{x}_1(t), \tilde{x}_2(t), \tilde{x}_3(t))$, and parameterize it as follows: starting from $\tilde{x}_1(0) \in [-c(0),\;c(0)]$, we define $(\tilde{x}_2(0), \tilde{x}_3(0))$ to lie on the circle satisfying
\[
\tilde{x}_2(0)^2 + \tilde{x}_3(0)^2 = a(0)^2\left(1 - b\left(\frac{\tilde{x}_1(0)^2}{c(0)^2}\right)\right),
\]
with the time evolution of the coordinates given by
\[
\tilde{x}_1(t) = \tilde{x}_1(0)\frac{c(t)}{c(0)}, \quad \tilde{x}_2(t) = \tilde{x}_2(0)\frac{a(t)}{a(0)}, \quad \tilde{x}_3(t) = \tilde{x}_3(0)\frac{a(t)}{a(0)}.
\]
Thus, the condition $d_{\Gamma(t)}(\mathbf{\tilde{x}}(t), t) = 0$ holds for the surface at any time $t$.
    \item \textbf{Normal motion} \\
Alternatively, we may choose material points on the same surface such that the velocity of the material points, $\bm{v}\left( \tilde{\mathbf{x}}(t),t\right)$, is in the normal direction. The material velocity is described by
\begin{equation}\label{normal_velocity}
\bm{v}(\mathbf{x},t) = V(\mathbf{x},t)\mathbf{n}(\mathbf{x},t), \quad \mathbf{n}(\mathbf{x},t) = \frac{\nabla d_{\Gamma(t)}(\mathbf{x},t)}{|\nabla d_{\Gamma(t)}(\mathbf{x},t)|}, \quad \text{and} \quad V(\mathbf{x},t) = -\frac{\partial_t d_{\Gamma(t)}(\mathbf{x},t)}{|\nabla d_{\Gamma(t)}(\mathbf{x},t)|},
\end{equation}
where $\mathbf{n}(\mathbf{x},t)$ represents the unit normal to the surface at each point.

\end{itemize}





\section{Numerical study of pattern formation on evolving surfaces}\label{patter_evolve}
In this section, we investigate the impact of surface evolution on Turing pattern formation by solving reaction–diffusion systems on time-dependent geometries. Through a series of examples, we illustrate how both isotropic and anisotropic surface deformations influence the emergence, spatial distribution, and symmetry of patterns over time, and how effectively the solver is able to capture these features. The reaction–diffusion dynamics are integrated in time using a first-order implicit–explicit backward differentiation formula (IMEX-BDF1) scheme.

\subsection{Isotropic growth and decay}
\paragraph{Logistic growth}

In this example, we consider an evolving sphere with initial radius \( r = 1 \), whose evolution is driven by a logistic growth function \( \eta(t) \), defined as follows:
\begin{equation}\label{eq:logistic}
\eta(t) = \frac{e^{g_{\text{rate}} t}}{1 + \frac{1}{K} \left( e^{g_{\text{rate}} t} - 1 \right)},
\quad \text{and} \quad 
\mathcal{L}_t(\mathbf{x},t) := \eta(t) (x_1(0),x_2(0),x_3(0)),
\end{equation}
where \(g_{\text{rate}} =0.1\) and \(K=1.5\). Here \(g_{\text{rate}}\) is the growth rate of the surface and \(K\) is the limiting final fixed size of the radius.  To project points onto the initial sphere \( \Gamma(0) \), we use the closest point projection operator $\pi(\mathbf{x}) := \frac{\mathbf{x}}{\|\mathbf{x}\|}$.
The time-dependent closest point projection \( \pi(\mathbf{x}, t) : \mathcal{N}^{t}_\delta \to \Gamma(t) \) is then given by the composition
\begin{equation}\label{eq:projection}
\pi(\mathbf{x}(t), t) := (\mathcal{L}_t \circ \pi)(\mathbf{x}) = \eta(t) \frac{\mathbf{x}(0)}{\|\mathbf{x}(0)\|},
\end{equation}
which ensures that \( \pi(\mathbf{x}(t), t) \) lies on the evolving surface \( \Gamma(t) \) with radius \( \eta(t) \). This composition accounts for both the geometry of the initial surface and its deformation over time.

We consider reaction-diffusion systems defined on evolving surfaces $\surf(t)$,
\begin{subequations}\label{eq:reaction_diffusion_2}
\begin{align}
    \frac{\partial u_1}{\partial t} &= \delta_{u_1} \Delta_{\Gamma(t)} u_1 + \alpha u_1 \left(1 - r_1 u_2^2\right) + u_2 \left(1 - r_2 u_1\right), \label{eq:reaction_diffusion_u1}  \\
    \frac{\partial u_2}{\partial t} &= \delta_{u_2} \Delta_{\Gamma(t)} u_2 + \beta u_2 \left(1 + \frac{\alpha r_1}{\beta} u_1 u_2\right) + u_1 \left(\gamma + r_2 u_2\right).\label{eq:reaction_diffusion_u2} 
\end{align}
\end{subequations}
For the simulation, we take $\alpha = 0.899$, $\beta = -0.91$, $\gamma = -\alpha$, $r_1 = 0.02$, $r_2 = 0.15$, \( \delta_{u_1} = 0.516\, \delta_{u_2} \), where \( \delta_{u_2} = 5 \cdot 10^{-3} \). Whenever a parameter differs from these default values in a specific experiment, it is explicitly stated.

\begin{figure}[!t]
    \centering
    \begin{subfigure}[b]{0.32\textwidth}

        \includegraphics[width=0.8\linewidth]{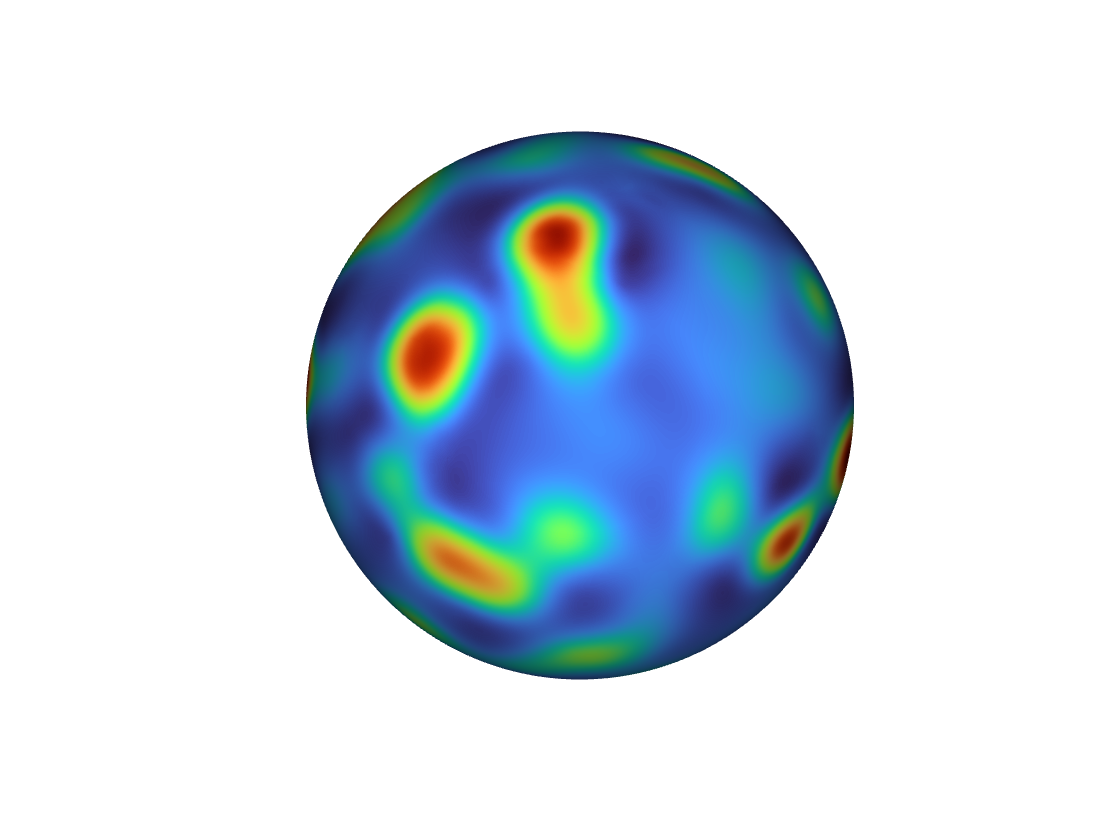}
        \caption{$t=10$}
    \end{subfigure}
    \hfill
    \begin{subfigure}[b]{0.32\textwidth}

        \includegraphics[width=0.9\linewidth]{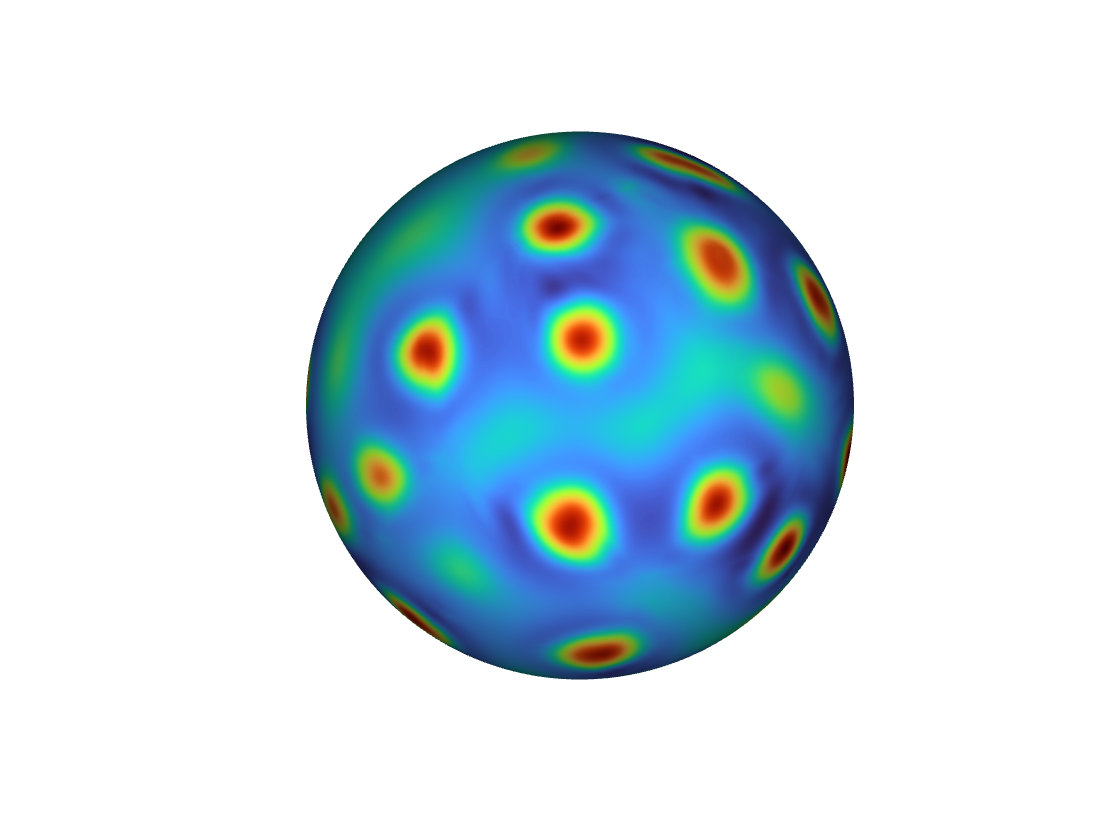}
        \caption{$t=30$}
    \end{subfigure}
    \hfill
    \begin{subfigure}[b]{0.32\textwidth}

     \includegraphics[
  width=1.00\linewidth
]{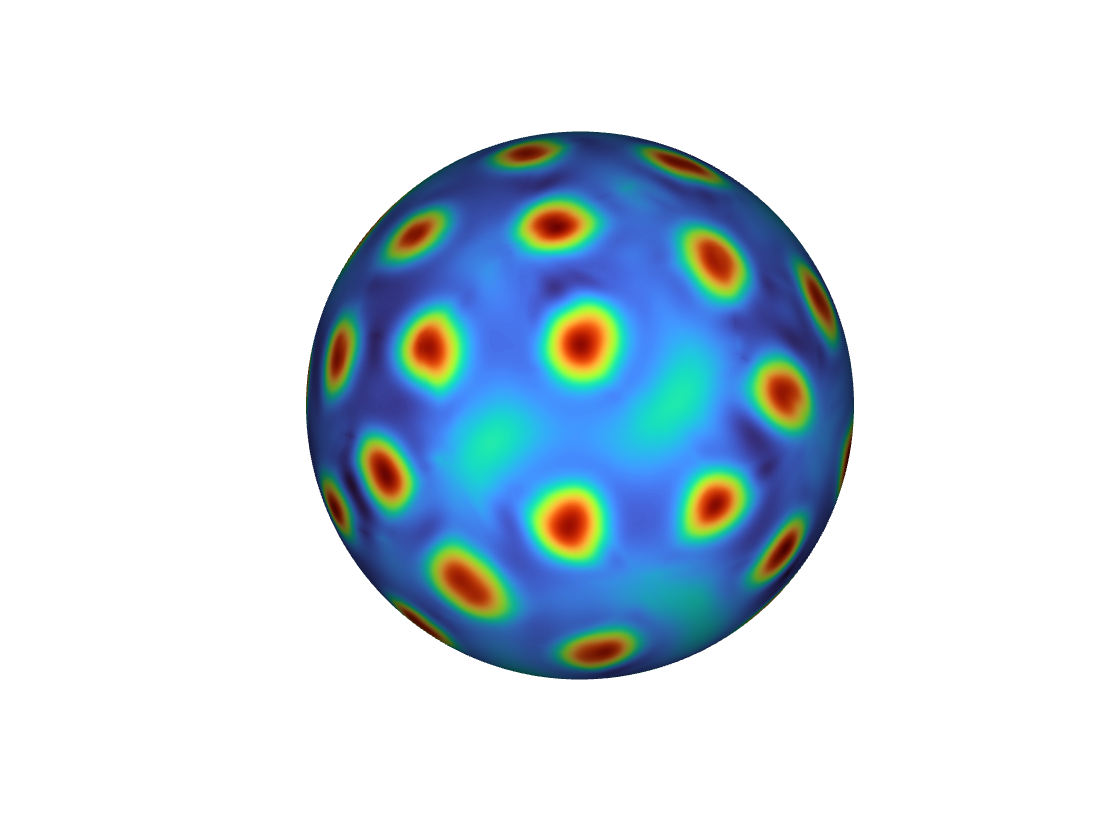}
        \caption{$t=50$}
    \end{subfigure}

    \caption{Concentrations for the Turing model on the sphere. 
    The sphere evolves according to the logistic growth function \eqref{eq:logistic} 
    with $\Delta t = 10^{-2}$.}
    \label{fig:stability_s}
\end{figure}

Figure~\ref{fig:stability_s} shows the evolution of the \( u_1 \) concentration in the Turing model on a sphere that expands isotropically according to the logistic growth law~\eqref{eq:logistic}. The plots correspond to times \( t = 10, 30 \), and \( 50 \). Initially, the pattern begins to emerge from small perturbations around the homogeneous state. As the surface grows, more spots form and gradually stabilize. The number and spacing of the spots appear to be influenced by the increasing surface area, which allows the activation of higher spatial modes over time.

\begin{figure*}[!t]
    \centering
    \begin{subfigure}[b]{0.45\textwidth}
        \centering
        \includegraphics[width=\linewidth]{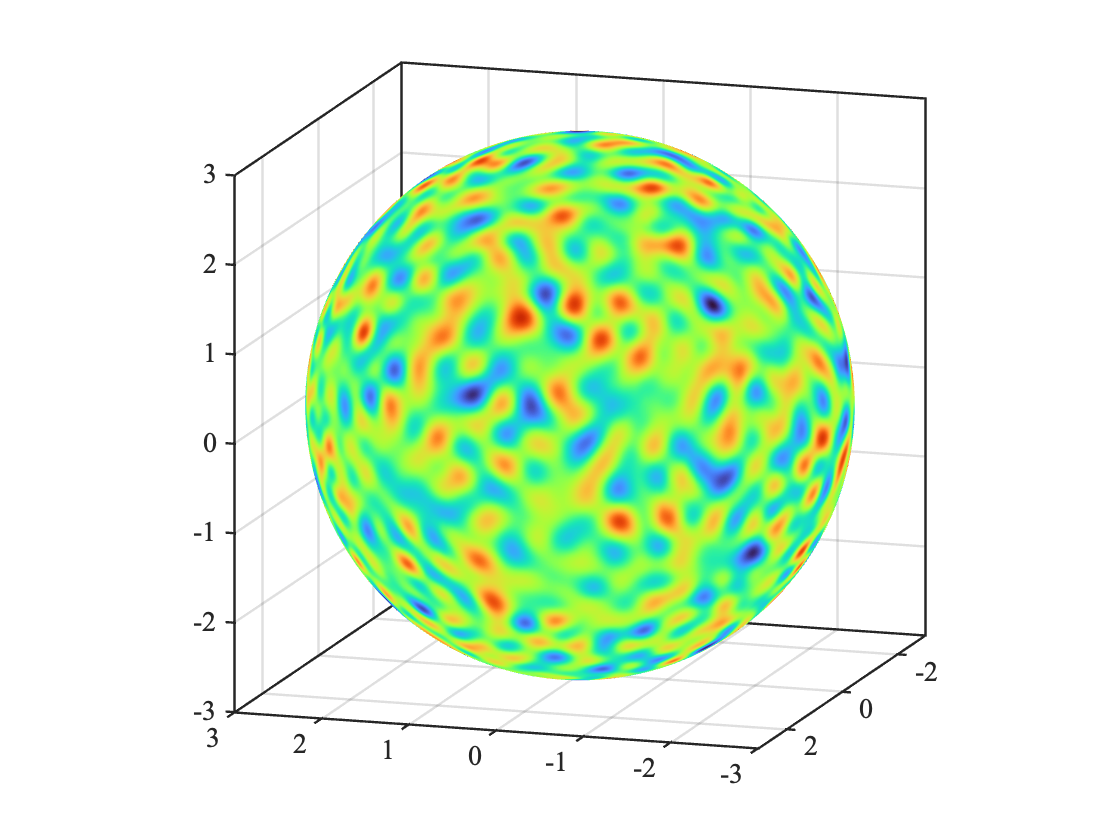}
        \caption{$t=0$}
    \end{subfigure}
    \hfill
    \begin{subfigure}[b]{0.45\textwidth}
        \centering
        \includegraphics[width=\linewidth]{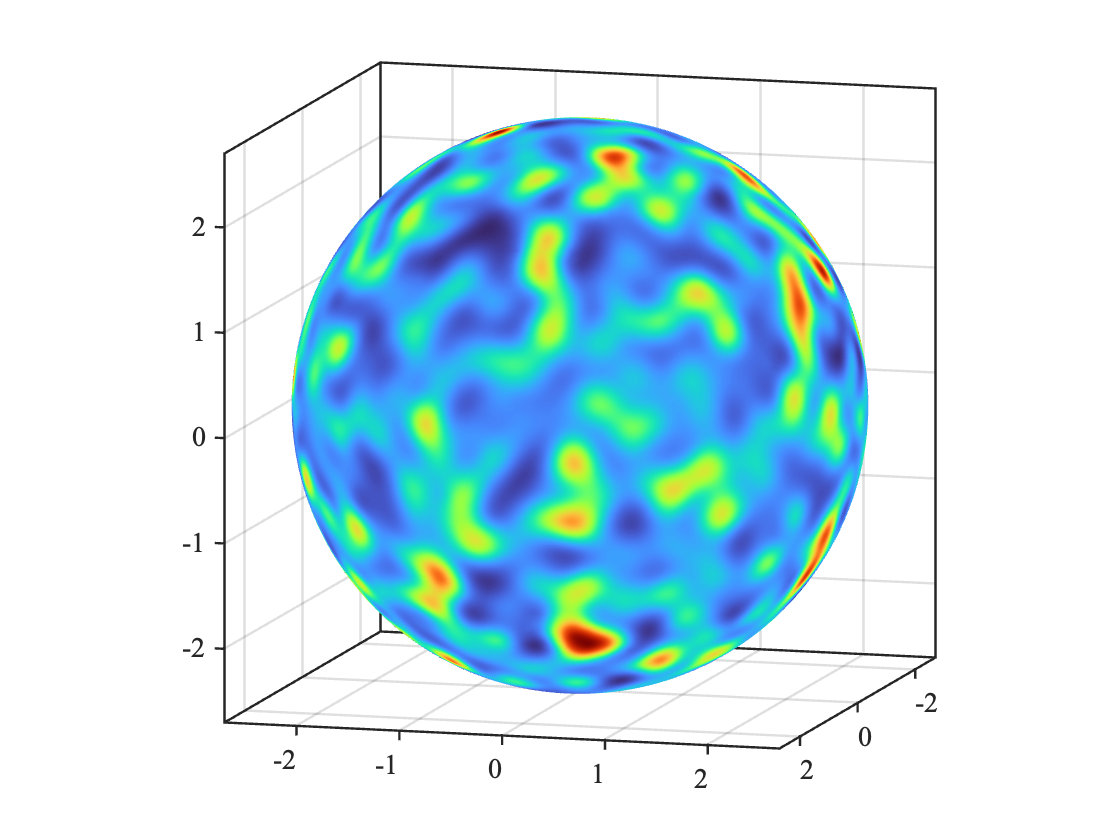}
        \caption{$t=5$}
    \end{subfigure}
    \vfill
    \begin{subfigure}[b]{0.45\textwidth}
        \centering
        \includegraphics[width=\linewidth]{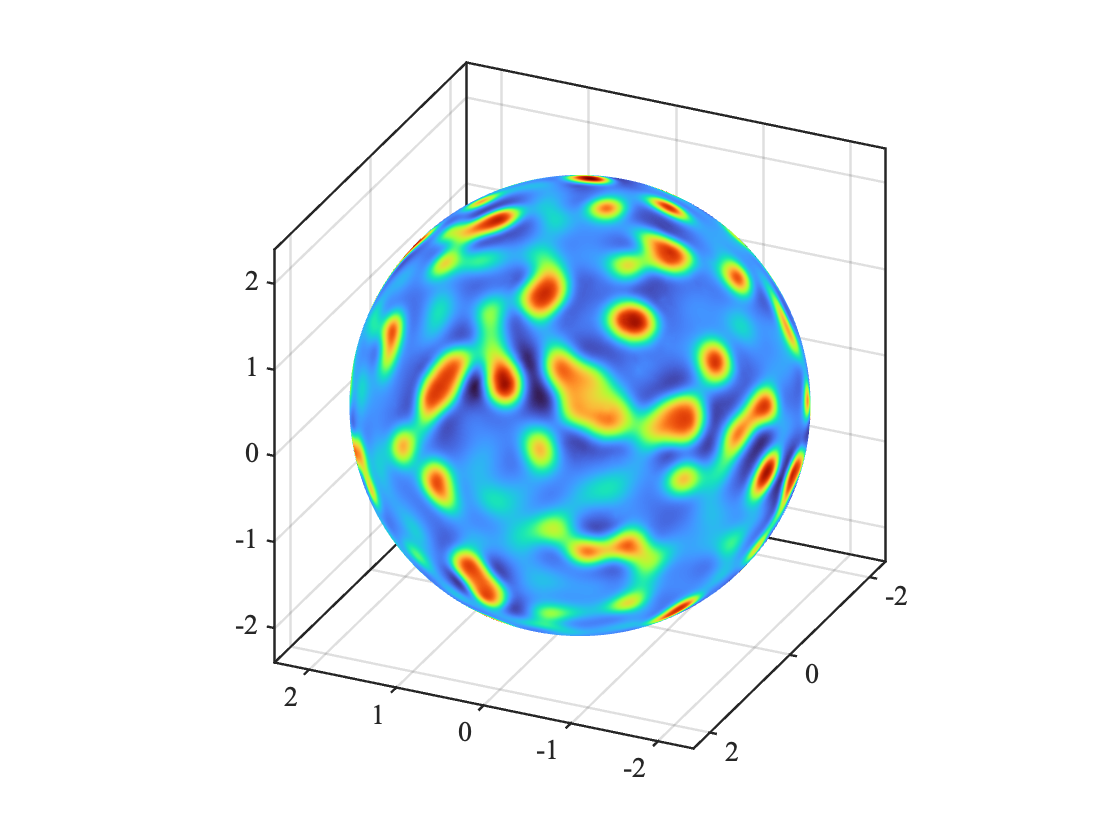}
        \caption{$t=10$}
    \end{subfigure}
        \hfill
    \begin{subfigure}[b]{0.45\textwidth}
        \centering
        \includegraphics[width=\linewidth]{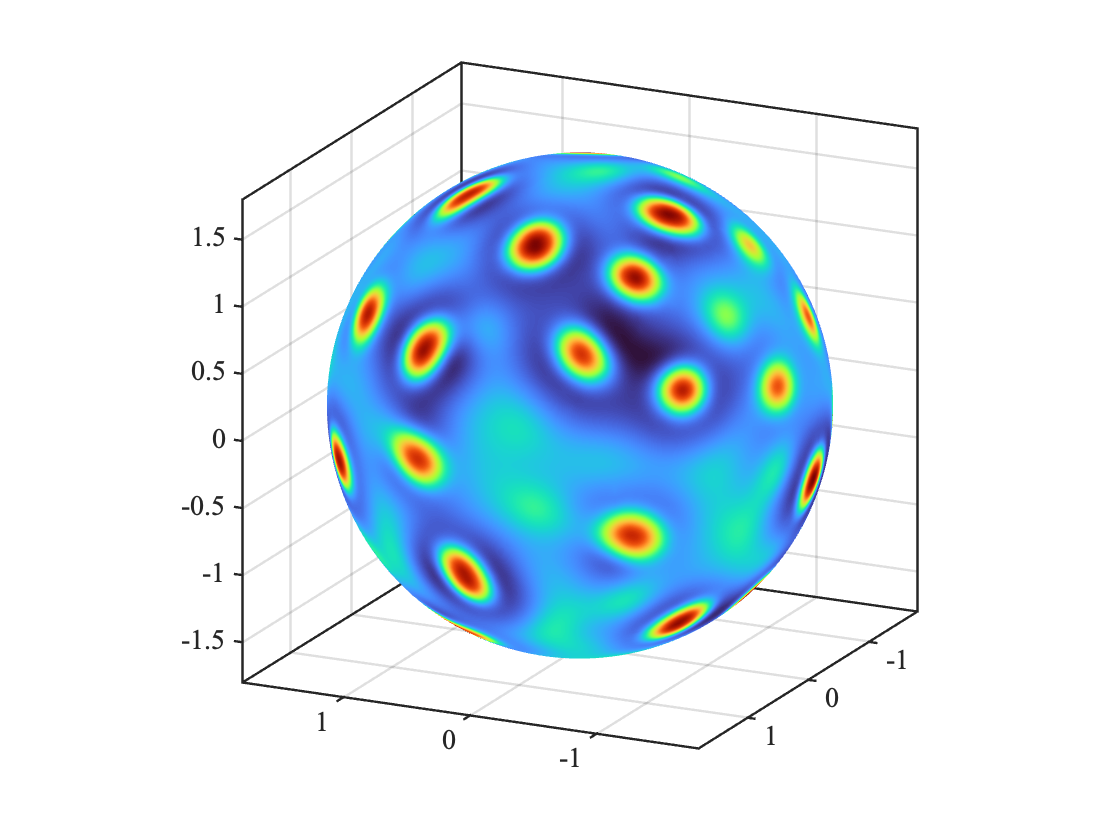}
        \caption{$t=20$}
    \end{subfigure}
    \vfill
    \begin{subfigure}[b]{0.45\textwidth}
        \centering
        \includegraphics[width=\linewidth]{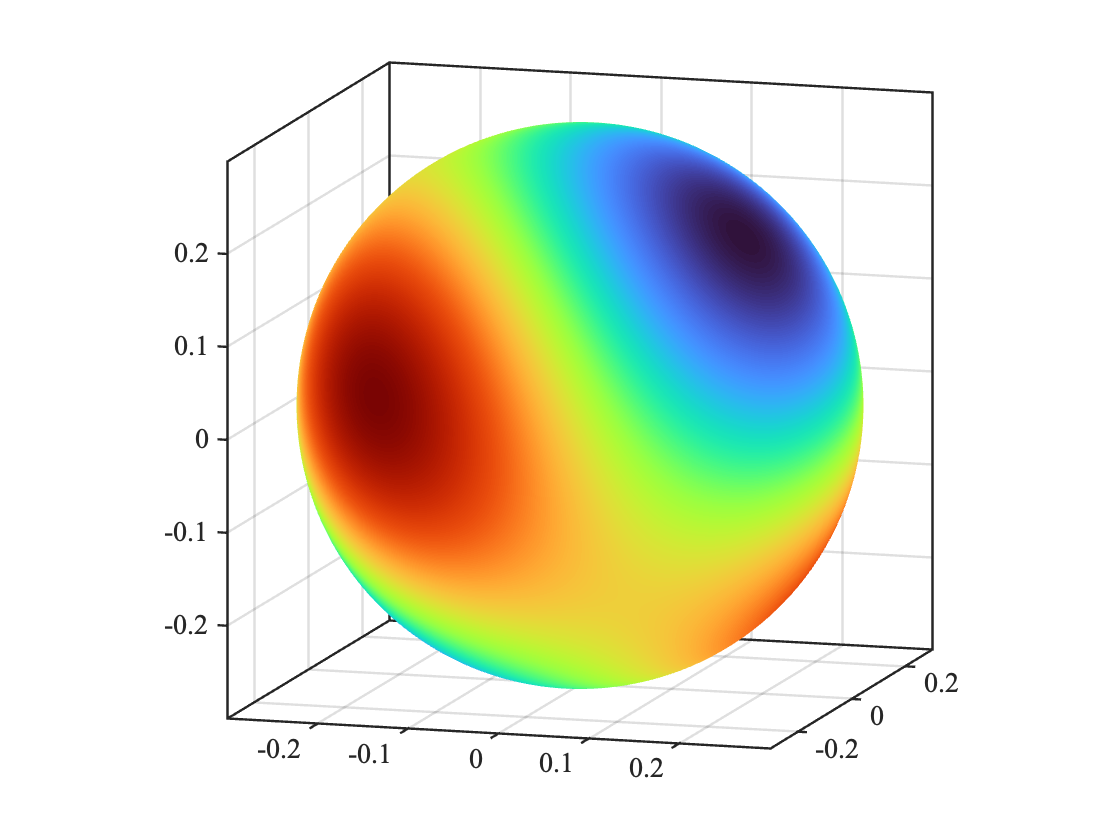}
        \caption{$t=45$}
    \end{subfigure}
    \hfill
     \begin{subfigure}[b]{0.45\textwidth}
        \centering
        \includegraphics[width=\linewidth]{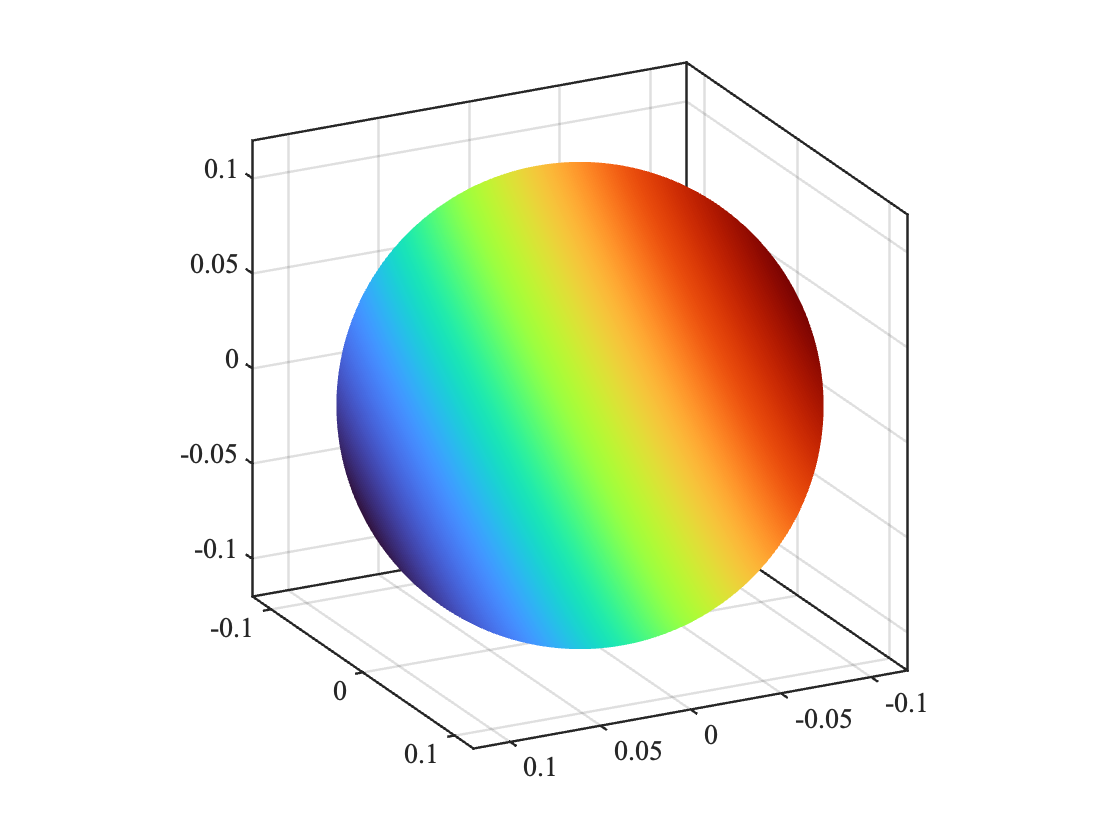}
        \caption{$t=48$}
    \end{subfigure}

    \caption{Patterns arising from the reaction-diffusion system \eqref{eq:reaction_diffusion} on a contracting sphere.}
    \label{fig:Turing_system_contract_sphere}
\end{figure*}

\begin{figure}[htbp]
    \centering
    \begin{tabular}{cccc}
        \includegraphics[width=0.3\textwidth]{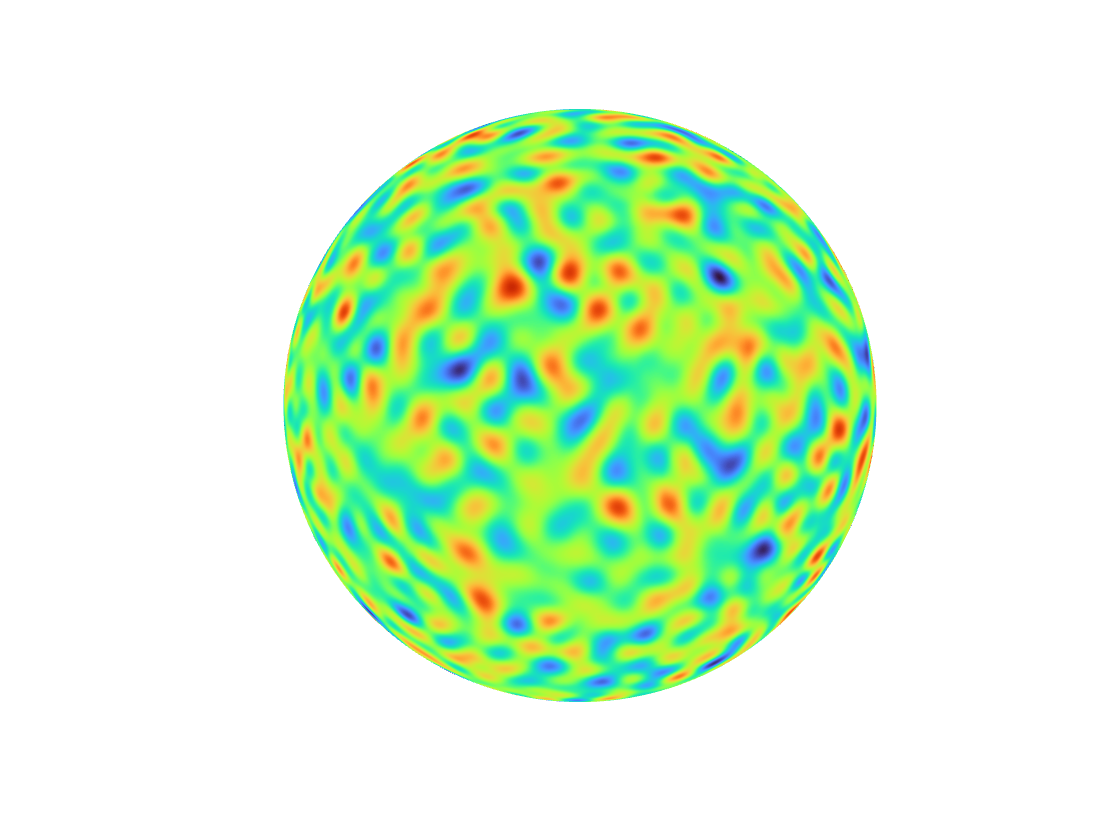} &
        \includegraphics[width=0.3\textwidth]{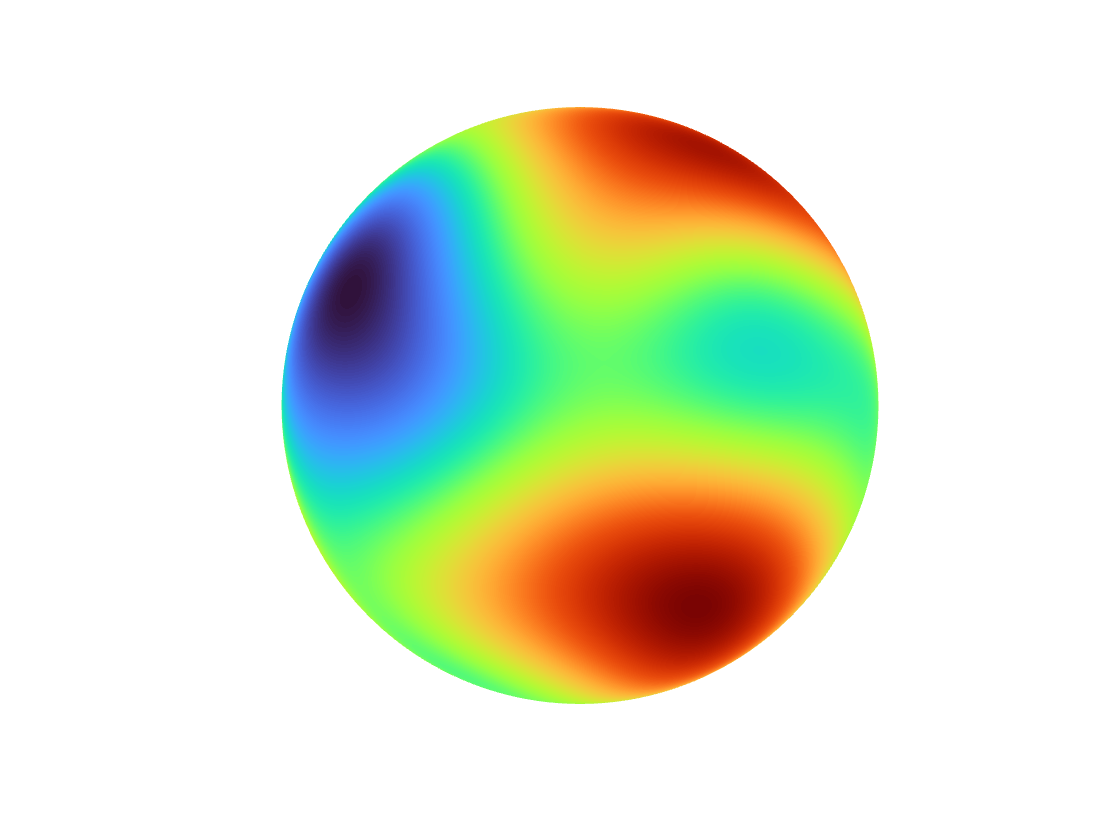} &
        \includegraphics[width=0.3\textwidth]{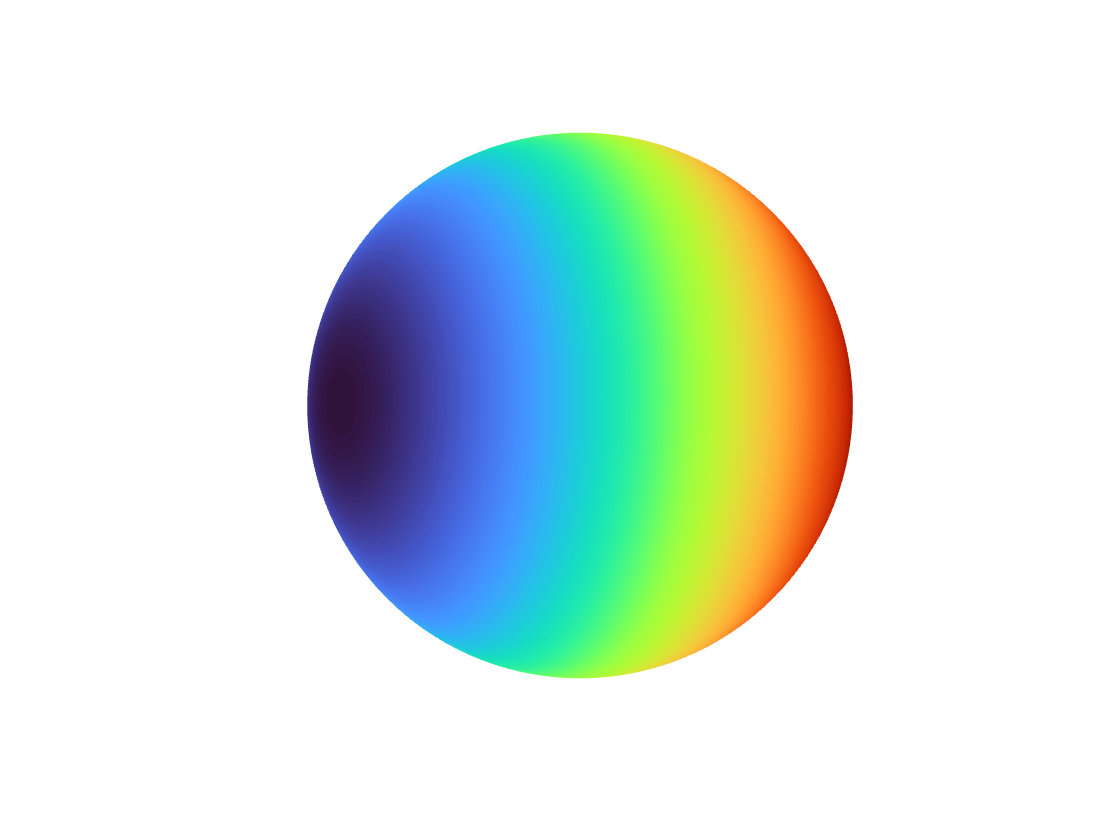} 
        
    \end{tabular}
   \caption{Patterns arising from the reaction--diffusion system~\eqref{eq:reaction_diffusion} on a stationary sphere with \( r = 3 \). The plots show the chemical concentration corresponding to increasing values of \( \delta_{u_1} \), with \( \delta_{u_1} = 0.0026 \), \( \delta_{u_1} = 0.0176 \), and \( \delta_{u_1} = 0.2758 \). The configurations are shown from left to right.}
    \label{fig:Turing_system_diff_parra}
\end{figure}

\paragraph{Linear contracting sphere}
In the second example, we repeat the experiment using the same reaction-diffusion system \eqref{eq:reaction_diffusion_2} and parameter values as in the expanding case. However, instead of surface expansion, we consider surface contraction with $\eta(t) = 1 - g_{rate} t$ with $g_{rate}=0.02$.  The initial surface $\Gamma(0)$ is a sphere of radius $r = 3$, which contracts over time to a sphere of radius $0.1$.  As it is shown in Figure~\ref{fig:Turing_system_contract_sphere} ,the final pattern on the sphere with radius equal to $0.1$ is the
same as the pattern forming  on the sphere with radius $r=3$ when $\delta_{u_1} = 0.2758$ (see Figure~\ref{fig:Turing_system_diff_parra}). Hence, this experiment demonstrates that contracting the surface can be thought of as increasing the values of $\delta_{u_1}$ without surface evolution. It is clear that less and less patterns are forming as the surface continues to contract, in contrast to the formation of more patterns when the surface is expanding as shown in Figure~\ref{fig:stability_s}.

\subsection{Anisotropic growth}
\begin{figure*}[!t]
    \centering
    \begin{subfigure}[b]{0.4\textwidth}
        \centering
        \includegraphics[width=\linewidth]{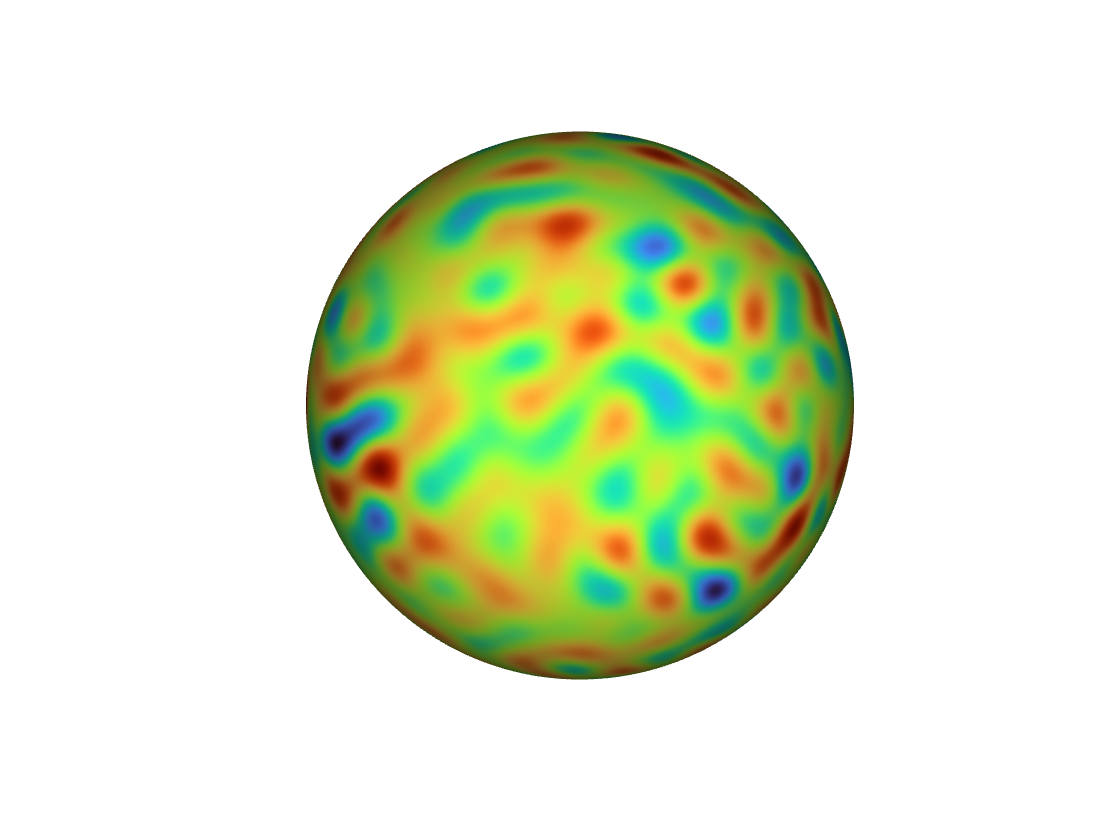}
        \caption{$t=0$}
    \end{subfigure}
    \hfill
    \begin{subfigure}[b]{0.4\textwidth}
        \centering
        \includegraphics[width=\linewidth]{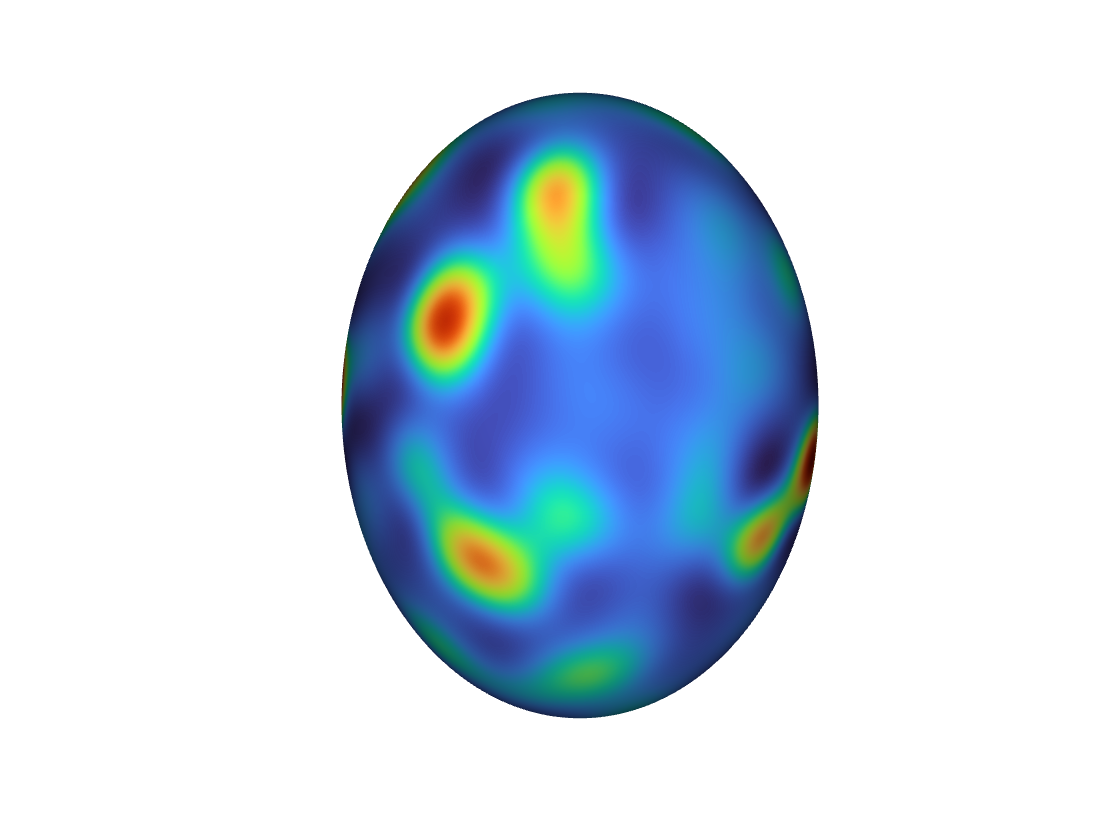}
        \caption{$t=10$}
    \end{subfigure}
    \vfill
    \begin{subfigure}[b]{0.4\textwidth}
        \centering
        \includegraphics[width=\linewidth]{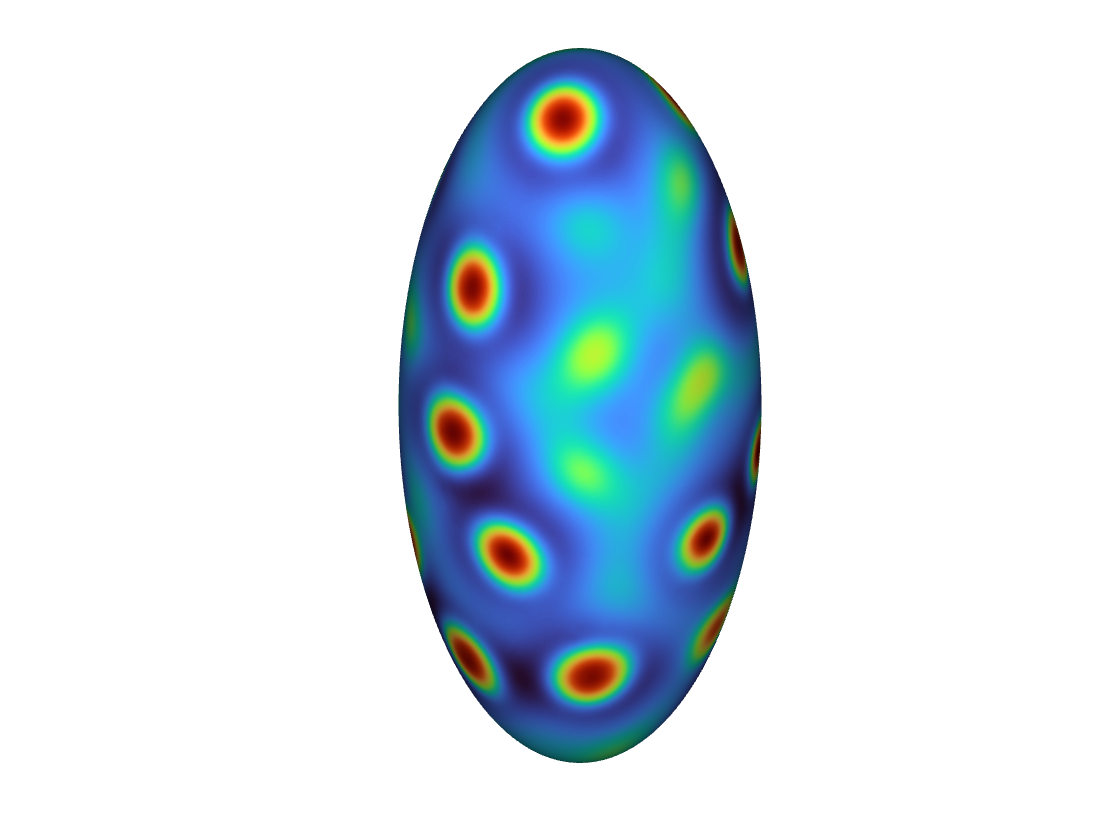}
        \caption{$t=30$}
    \end{subfigure}
    \hfill
    \begin{subfigure}[b]{0.4\textwidth}
        \centering
        \includegraphics[width=\linewidth]{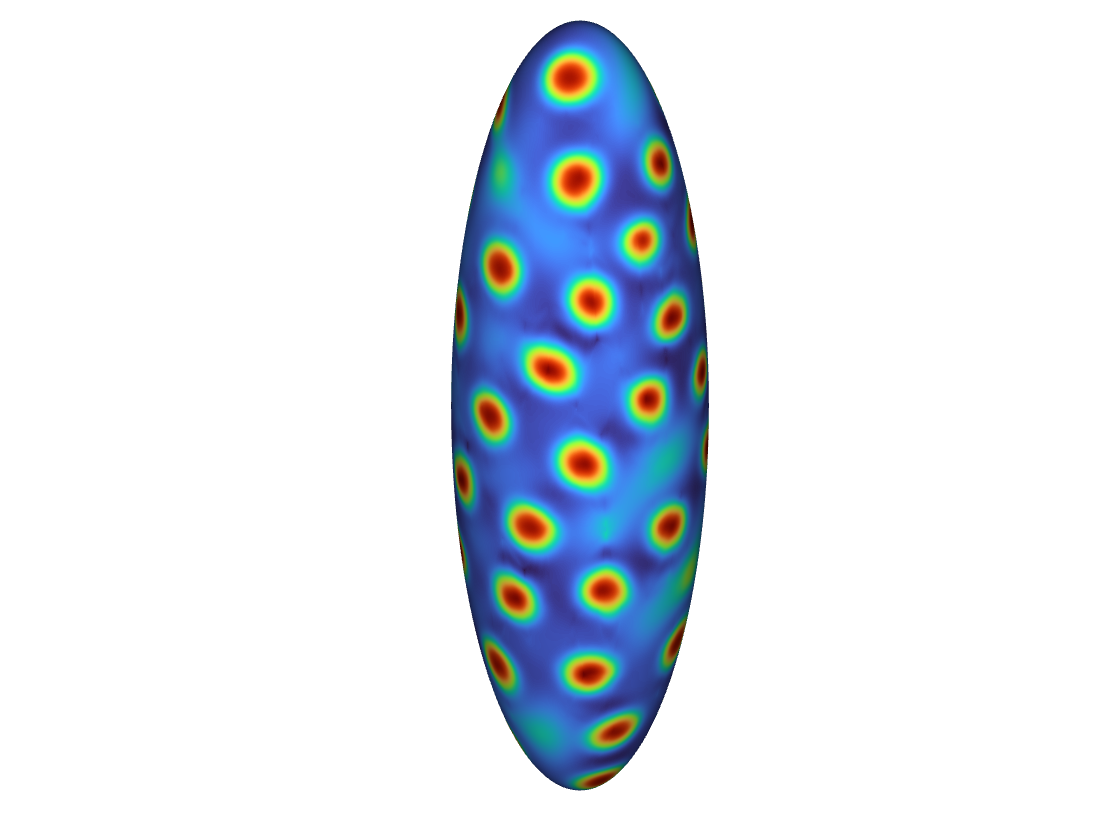}
        \caption{$t=60$}
    \end{subfigure}
    \caption{Patterns arising from the reaction-diffusion system \eqref{eq:reaction_diffusion} on an anisotropic growth of the unit sphere, where the sphere evolves into an ellipsoid. The simulations use the implicit–explicit backward differentiation formula (IMEX-BDF1) scheme with with $\dt=0.1$.}
    \label{fig:TS_sphere_into_ellipsoid}
\end{figure*}

\paragraph{Anisotropic evolution of a sphere into an ellipsoid.}
We begin with the initial surface $\surf(0)$ and assume that it moves with a velocity field given by $\bm{v}(\mathbf{x}(t), t) := \left(0, 0, g_{\text{rate}} x_3(0)\right),$
where $g_{\text{rate}}$ is a fixed growth factor. This implies that the anisotropic deformation of the surface is governed by the following time-dependent closest point projection
\begin{equation}\label{eq:projection_1}
\pi(\mathbf{x}(t), t) := (\mathcal{L}_t \circ \pi)(\mathbf{x}(0),t) 
= \left( \frac{x_1(0)}{\|\mathbf{x}(0)\|}, \; \frac{x_2(0)}{\|\mathbf{x}(0)\|}, \; \frac{x_3(0)}{\|\mathbf{x}(0)\|}(1 + g_{\text{rate}} t) \right)^{\!\top}.
\end{equation}
Accordingly, the surface at time $t$ can be described as the zero level set of a time-dependent distance function:
\[
\surf(t) := \left\{ \mathbf{x} \in \mathbb{R}^3 \;:\; d_{\Gamma(t)}(\mathbf{x}, t) = 0 \right\},\quad\text{where}\quad d_{\Gamma(t)}(\mathbf{x}(t), t) := x_1^2 + x_2^2 + \frac{x_3^2}{(1 + g_{\text{rate}} t)^2} - 1.
\]
The final time is set to $T = 60$ and the growth rate is $g_{\text{rate}} = 0.04$, so that the final surface $\surf(T)$ is an ellipsoid whose equation is $x_1^2+x_2^2+\frac{x_3^2}{(\frac{17}{5})^2}=1.$ In Figure~\ref{fig:TS_sphere_into_ellipsoid}, we observe the continuous formation of spots as the surface evolves. This demonstrates how anisotropic deformation into a prolate ellipsoid enhances directional effects, leading to elongated and spatially organized patterns aligned with the principal axis of growth.

\paragraph{An anisotropic evolution of the dumbbell}
We define the evolving closed surface \( \Gamma(t) \) as the zero level set of a time-dependent distance function:
\begin{equation}
d_{\Gamma(t)}\left(\mathbf{x}(t),t\right) = x_1^2 + x_2^2 + a(t)^2 b\left(\frac{x_3^2}{c(t)^2}\right) - a(t)^2, \quad \text{i.e., } \surf(t):=\{\mathbf{x}\in \mathbb{R}^{3}\; : \; d_{\Gamma(t)}\left(\mathbf{x}(t),t\right)=0\}.
\end{equation}
Here the functions $a$, $b$, and $c$ are given by
\[
a(t) = 0.1 + 0.05 \sin(2\pi t),\quad
b(s) = 200s\left(s - \frac{199}{200}\right),\quad 
c(t) = 1 + 0.2 \sin(4\pi t).
\]
 The surface velocity $\bm{v}$ in $\mathbf{x}_k(t)$ is given by Eq.~\eqref{normal_velocity}, and the coordinates evolve in time according to
\[
(\mathbf{x}_k(t))_1 = (\mathbf{x}_k(0))_1 \frac{a(t)}{a(0)}, \quad (\mathbf{x}_k(t))_2 = (\mathbf{x}_k(0))_2 \frac{a(t)}{a(0)}, \quad (\mathbf{x}_k(t))_3 = (\mathbf{x}_k(0))_3 \frac{c(t)}{c(0)},
\]

\begin{figure}[htbp]
    \centering
    \begin{tabular}{cccc}
        \includegraphics[width=0.3\textwidth]{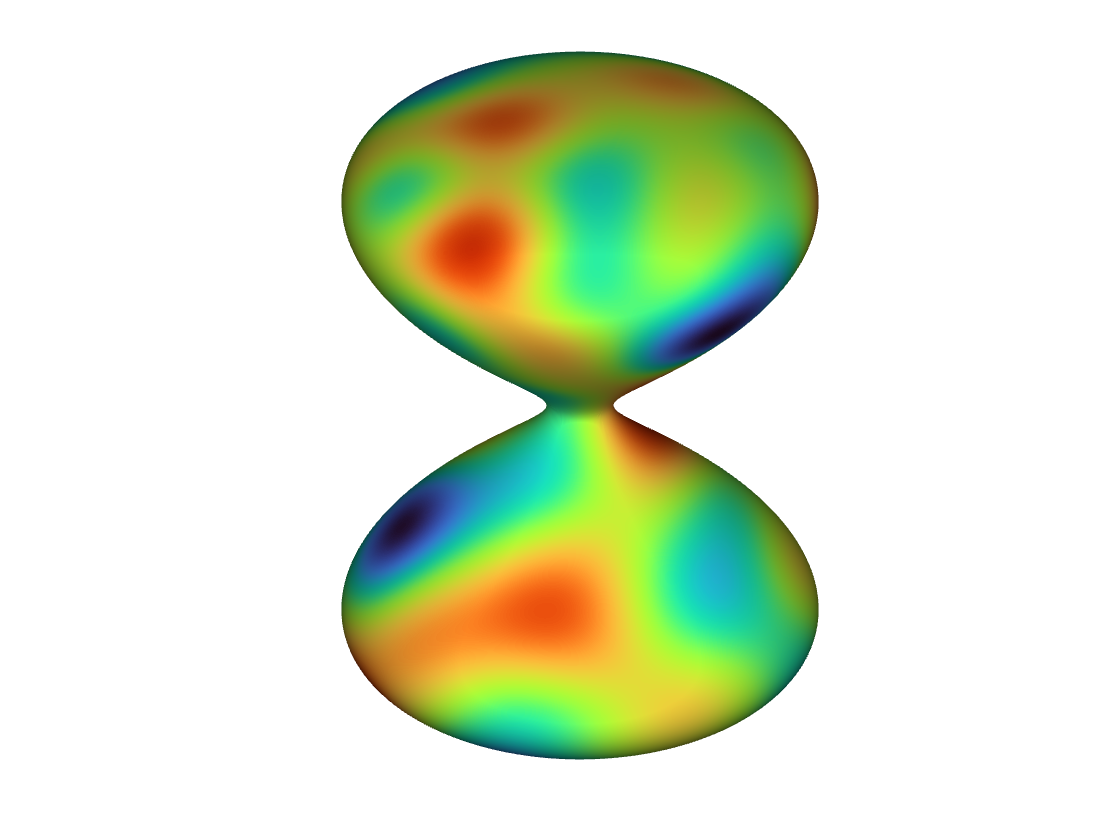} &
        \includegraphics[width=0.3\textwidth]{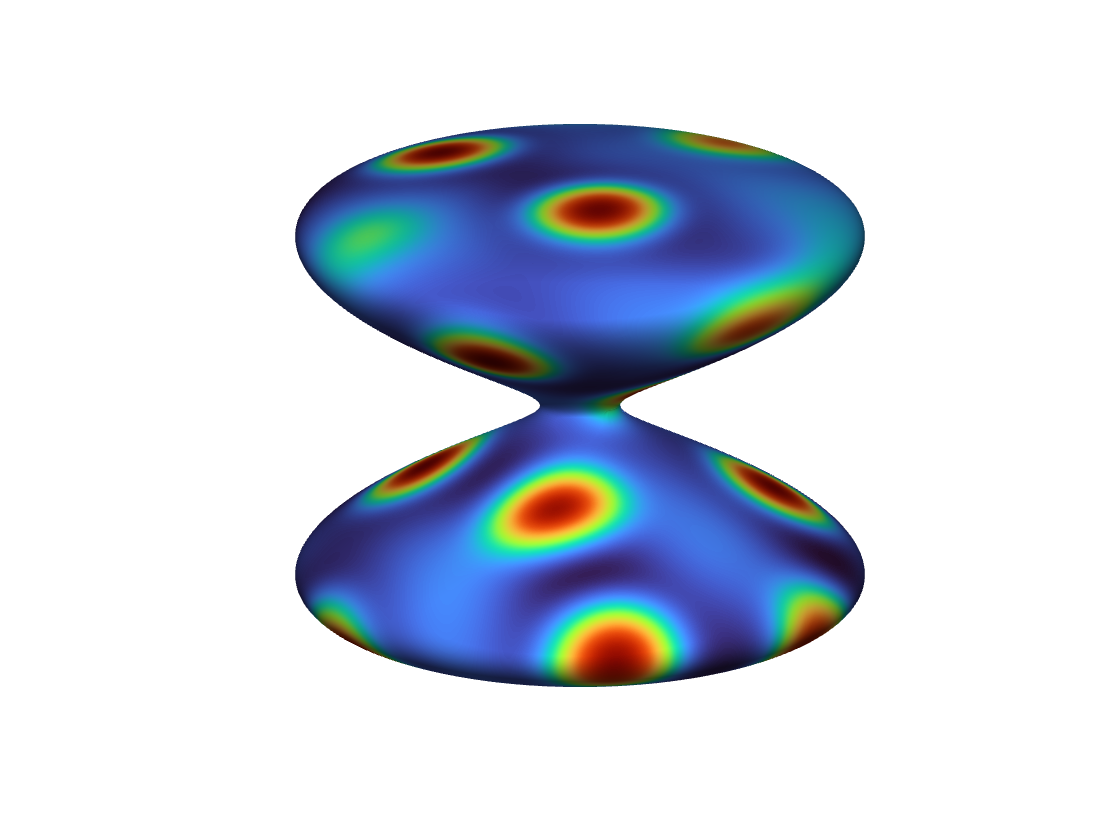} &
        \includegraphics[width=0.3\textwidth]{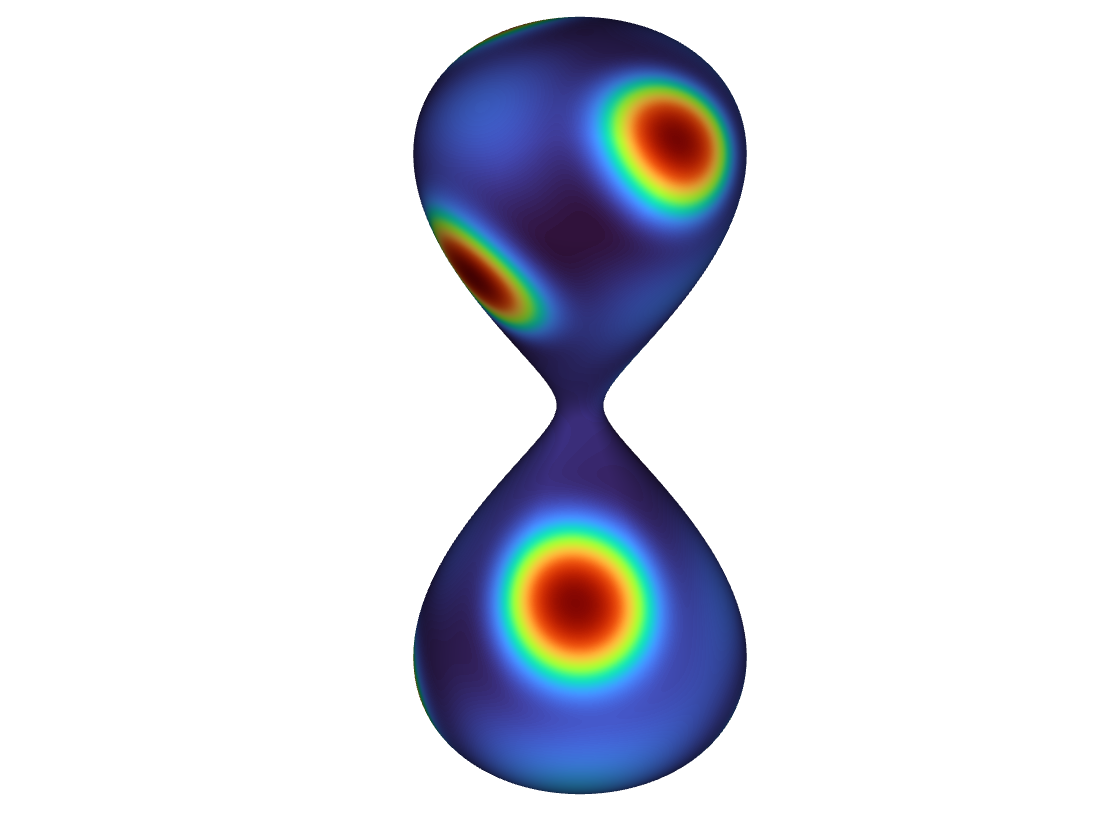} 
        
    \end{tabular}

    \caption{Patterns arising from the reaction-diffusion system \eqref{eq:reaction_diffusion} on an anisotropic growth of the dumbbell shape.
    The simulations use the implicit–explicit backward differentiation formula (IMEX-BDF1) scheme with with $\dt=0.1$.}
    \label{fig:Turing_system_combined_dumbell}
\end{figure}

\begin{figure}[htbp]
    \centering
    \begin{tabular}{cccc}
        \includegraphics[width=0.3\textwidth]{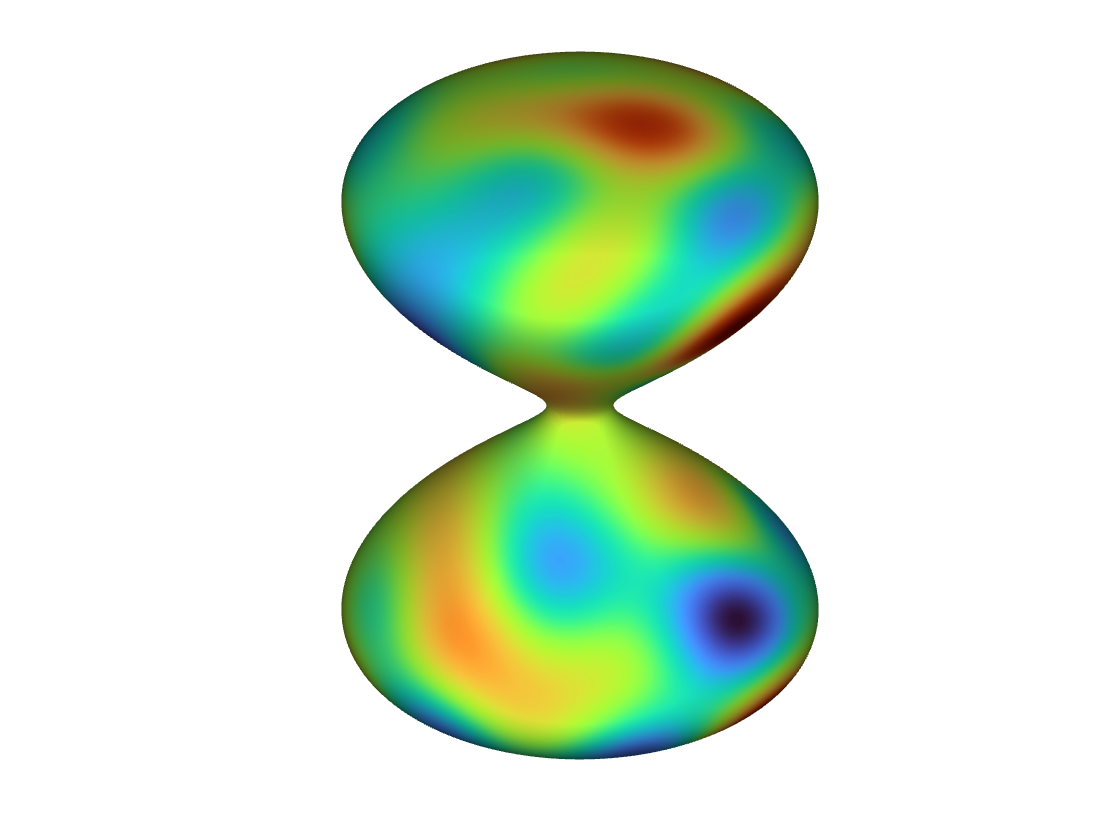} &
        \includegraphics[width=0.3\textwidth]{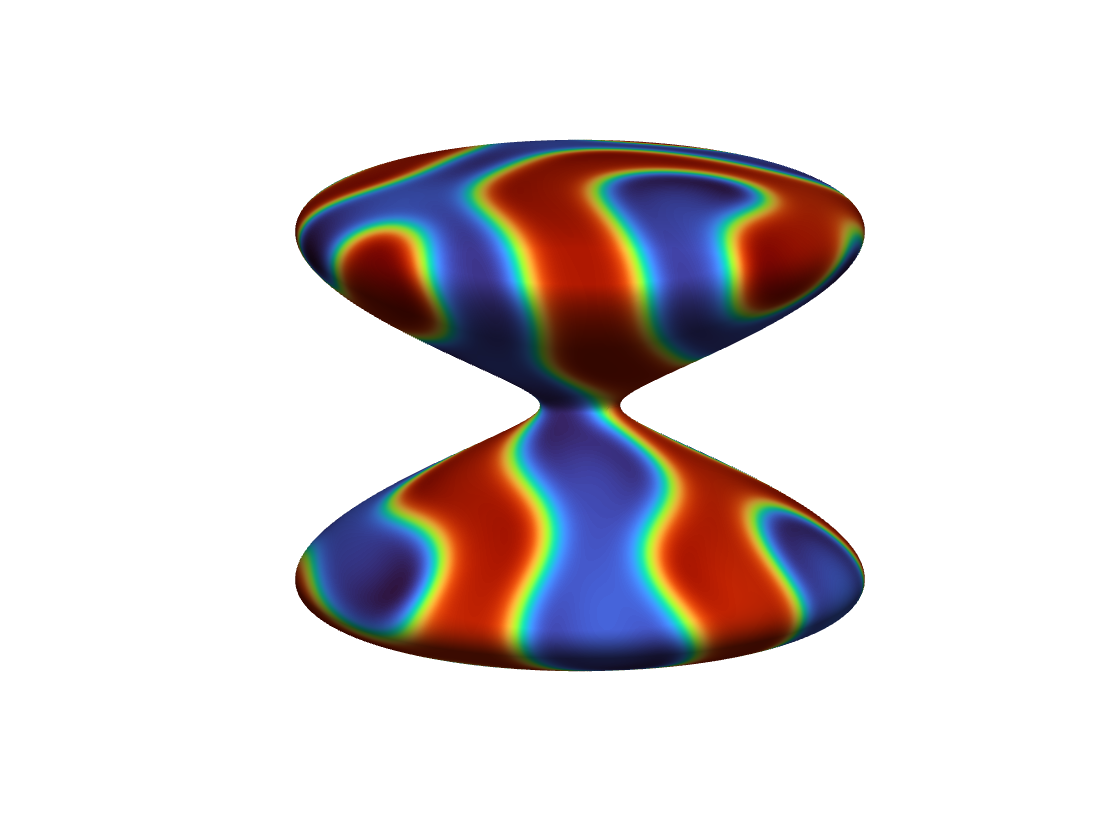} &
        \includegraphics[width=0.3\textwidth]{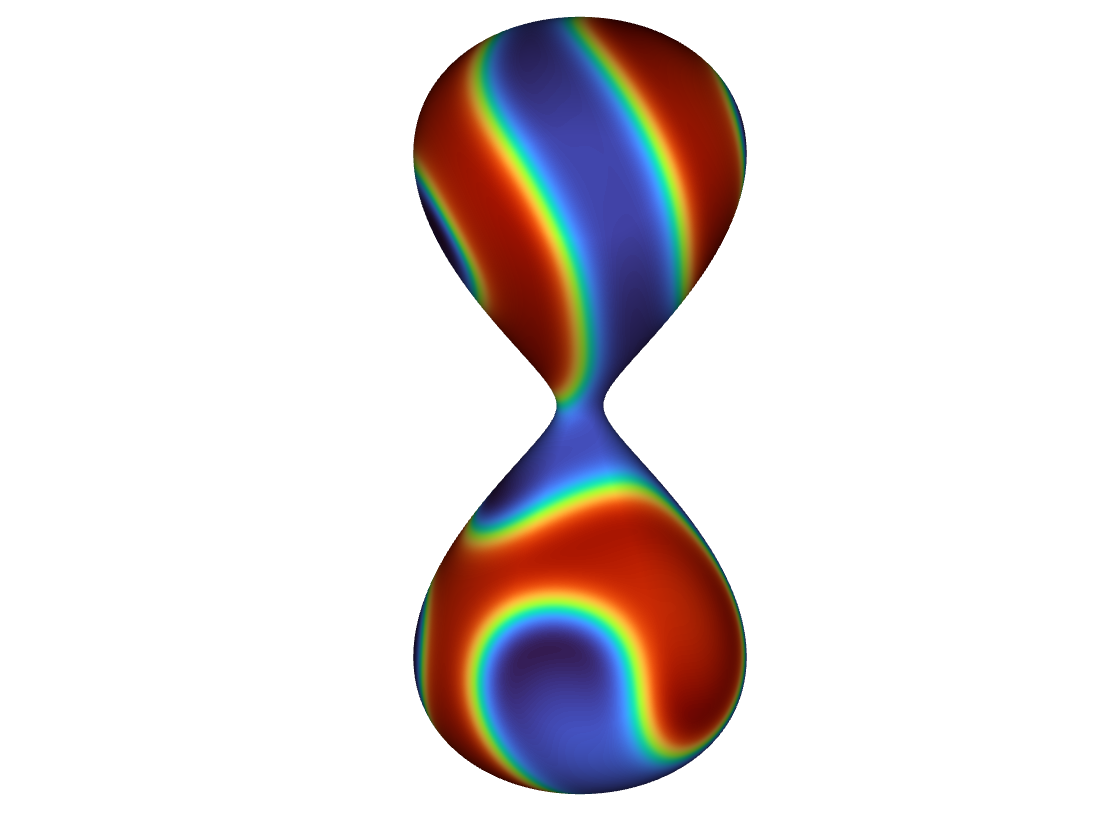} 
        
    \end{tabular}

    \caption{Patterns arising from the reaction-diffusion system \eqref{eq:reaction_diffusion} on an anisotropic growth of the dumbbell shape. The red regions represent the activation part, and the blue regions
represent the deactivation patterns, respectively.}
    \label{fig:Turing_system_strip_dumbell}
\end{figure}

Figure~\ref{fig:Turing_system_combined_dumbell} shows the time evolution of the \( u_1 \) solution in the Turing model, obtained with the spot-forming parameter set listed in Table~\ref{tab_par}. The simulation runs until final time \( T = 60 \), with snapshots shown at \( t = 0, 10, 40 \), and \( 60 \). As the surface deforms, localized spot patterns emerge and stabilize, influenced by the curvature variations near the neck and lobes. These geometric features appear to guide the spatial distribution and persistence of the spots, indicating a strong coupling between surface shape and pattern dynamics. In contrast, Figure~\ref{fig:Turing_system_strip_dumbell} presents results obtained with a different parameter set that favors stripe formation. While the underlying surface evolution remains the same, the altered kinetics lead to the emergence of stripe-like structures aligned along the axial direction. This highlights the role of both geometry and reaction parameters in determining the symmetry and mode of pattern formation.
\begin{table}[h!]
\centering
\caption{The table shows the values of the parameters of equations~\eqref{eq:reaction_diffusion_u1} and \eqref{eq:reaction_diffusion_u2} used in the numerical experiments shown in Figure~ \eqref{fig:Turing_system_combined_dumbell} and~\eqref{fig:Turing_system_strip_dumbell}. We set $\delta_{u_1} = 0.5166\delta_{u_2}$. }
\vspace{0.4cm} 
\begin{tabular}{|c|c|c|c|c|c|c|c|}
\hline
Pattern & $\delta_{u_2}$ & $\alpha$ & $\beta$ & $\gamma$ & $r_1$ & $r_2$ & Final time \\ \hline
Spots   & $4.5 \times 10^{-3}$ & 0.899 & -0.91 & -0.899 & 0.02 & 0.15 & 60 \\ \hline
Stripes & $5 \times 10^{-3}$ & 1.899 & -0.95 & -1.899 & 1.5  & 0   & 60 \\ \hline
\end{tabular}
\label{tab_par}
\end{table}

\section*{Acknowledgment}
The author gratefully acknowledges Balázs Kovács for providing the MATLAB
implementation of the arbitrary Lagrangian--Eulerian (ALE) mapping algorithms
used in this work.

\bibliographystyle{siamplain}
\bibliography{Ref.bib}

@article{bunow1980pattern,
  title={Pattern formation by reaction-diffusion instabilities: Application to morphogenesis in {D}rosophila},
  author={Bunow, Barry and Kernevez, Jean-Pierre and Joly, Gislaine and Thomas, Daniel},
  journal={Journal of theoretical biology},
  volume={84},
  number={4},
  pages={629--649},
  year={1980},
  publisher={Elsevier}
}

@book{fornberg1998practical,
  author    = {Bengt Fornberg},
  title     = {A Practical Guide to Pseudospectral Methods},
  year      = {1998},
  publisher = {Cambridge University Press}
}

@article{canuto1988schwarz,
  title={The {Schwarz} algorithm for spectral methods},
  author={Canuto, Claudio and Funaro, Daniele},
  journal={SIAM journal on numerical analysis},
  volume={25},
  number={1},
  pages={24--40},
  year={1988},
  publisher={SIAM}
}

@article{gillman2014direct,
  title={A direct solver with $\mathcal{O}(N)$ complexity for variable coefficient elliptic PDEs discretized via a high-order composite spectral collocation method},
  author={Gillman, Adrianna and Martinsson, Per-Gunnar},
  journal={SIAM Journal on Scientific Computing},
  volume={36},
  number={4},
  pages={A2023--A2046},
  year={2014},
  publisher={SIAM}
}

@article{dziuk2013finite,
  title={Finite element methods for surface {PDEs}},
  author={Dziuk, Gerhard and Elliott, Charles M},
  journal={Acta Numerica},
  volume={22},
  pages={289--396},
  year={2013},
  publisher={Cambridge University Press}
}

@article{Persson,
author = {Persson, Per-Olof and Strang, Gilbert},
title = {A Simple Mesh Generator in {MATLAB}},
journal = {SIAM Review},
volume = {46},
number = {2},
pages = {329-345},
year = {2004},
doi = {10.1137/S0036144503429121},

}

@article{maini1997spatial,
  title={Spatial pattern formation in chemical and biological systems},
  author={Maini, Philip K and Painter, Kevin J and Chau, Helene Nguyen Phong},
  journal={Journal of the Chemical Society, Faraday Transactions},
  volume={93},
  number={20},
  pages={3601--3610},
  year={1997},
  publisher={Royal Society of Chemistry}
}

@article{fortunato2022highorder,
  title={A high-order fast direct solver for surface PDEs},
  author={Fortunato, Daniel},
  journal={SIAM Journal on Scientific Computing},
  volume={46},
  number={4},
  pages={A2582--A2606},
  year={2024},
  publisher={SIAM}
}

@article{martinsson2013direct,
  title={A direct solver for variable coefficient elliptic PDEs discretized via a composite spectral collocation method},
  author={Martinsson, Per-Gunnar},
  journal={Journal of Computational Physics},
  volume={242},
  pages={460--479},
  year={2013},
  publisher={Elsevier}
}

@article{jeong2017numerical,
  title={Numerical simulation of the zebra pattern formation on a three-dimensional model},
  author={Jeong, Darae and Li, Yibao and Choi, Yongho and Yoo, Minhyun and Kang, Dooyoung and Park, Junyoung and Choi, Jaewon and Kim, Junseok},
  journal={Physica A: Statistical Mechanics and its Applications},
  volume={475},
  pages={106--116},
  year={2017},
  publisher={Elsevier}
}

@article{ciarlet2002finite,
  title={The Finite Element Method for Elliptic Problems (SIAM, Philadelphia)},
  author={Ciarlet, PG},
  year={2002}
}

@book{mathew2008domain,
  title={Domain decomposition methods for the numerical solution of partial differential equations},
  author={Mathew, Tarek Poonithara Abraham},
  year={2008},
  publisher={Springer}
}

@article{jonsson2017cut,
  title={Cut finite element methods for elliptic problems on multipatch parametric surfaces},
  author={Jonsson, Tobias and Larson, Mats G and Larsson, Karl},
  journal={Computer Methods in Applied Mechanics and Engineering},
  volume={324},
  pages={366--394},
  year={2017},
  publisher={Elsevier}
}

@article{davis2004algorithm,
  author  = {Timothy A. Davis},
  title   = {Algorithm 832: \textsc{{UMFPACK}} V4.3---An Unsymmetric-Pattern Multifrontal Method},
  journal = {ACM Transactions on Mathematical Software (TOMS)},
  volume  = {30},
  number  = {2},
  pages   = {196--199},
  year    = {2004}
}

@book{briggs2000multigrid,
  title={A multigrid tutorial},
  author={Briggs, William L and Henson, Van Emden and McCormick, Steve F},
  year={2000},
  publisher={SIAM}
}

@article{turing1990chemical,
  title={The chemical basis of morphogenesis},
  author={Turing, Alan Mathison},
  journal={Bulletin of mathematical biology},
  volume={52},
  pages={153--197},
  year={1990},
  publisher={Springer}
}

@article{ascher1995implicit,
  title={Implicit-explicit methods for time-dependent partial differential equations},
  author={Ascher, Uri M and Ruuth, Steven J and Wetton, Brian TR},
  journal={SIAM Journal on Numerical Analysis},
  volume={32},
  number={3},
  pages={797--823},
  year={1995},
  publisher={SIAM}
}

@book{strikwerda2004finite,
  title={Finite difference schemes and partial differential equations},
  author={Strikwerda, John C},
  year={2004},
  publisher={SIAM}
}

@book{varga1962matrix,
  title={Matrix iterative methods},
  author={Varga, Richard S},
  year={1962},
  publisher={Prentice Hall Incorporated}
}

@book{martinsson2019fast,
  title={Fast direct solvers for elliptic PDEs},
  author={Martinsson, Per-Gunnar},
  year={2019},
  publisher={SIAM}
}

@article{george1973nested,
  title={Nested dissection of a regular finite element mesh},
  author={George, Alan},
  journal={SIAM journal on numerical analysis},
  volume={10},
  number={2},
  pages={345--363},
  year={1973},
  publisher={SIAM}
}

@article{fortunato2021ultraspherical,
  title={The ultraspherical spectral element method},
  author={Fortunato, Daniel and Hale, Nicholas and Townsend, Alex},
  journal={Journal of Computational Physics},
  volume={436},
  pages={110087},
  year={2021},
  publisher={Elsevier}
}

@book{M,
  address = {Cambridge},
  author = {Biggs, Norman},
  edition = {2nd},
  publisher = {Cambridge University Press},
  title = {Algebraic Graph Theory},
  year = 1993
}

@article{kovacs2019computing,
  title={Computing arbitrary {L}agrangian {E}ulerian maps for evolving surfaces},
  author={Kov{\'a}cs, Bal{\'a}zs},
  journal={Numerical Methods for Partial Differential Equations},
  volume={35},
  number={3},
  pages={1093--1112},
  year={2019},
  publisher={Wiley Online Library}
}

@book{gilbarg1977elliptic,
  title={Elliptic partial differential equations of second order},
  author={Gilbarg, David and Trudinger, Neil S and Gilbarg, David and Trudinger, NS},
  volume={224},
  number={2},
  year={1977},
  publisher={Springer}
}

@article{barrio1999two,
  title={A two-dimensional numerical study of spatial pattern formation in interacting {T}uring systems},
  author={Barrio, RA and Varea, C and Arag{\'o}n, JL and Maini, PK},
  journal={Bulletin of mathematical biology},
  volume={61},
  number={3},
  pages={483--505},
  year={1999},
  publisher={Elsevier}
}

@book{trefethen2000spectral,
author = {Trefethen, Lloyd N.},
title = {Spectral Methods in MATLAB},
publisher = {Society for Industrial and Applied Mathematics},
year = {2000},
doi = {10.1137/1.9780898719598},
address = {},
edition   = {}
}

@book{quarteroni2008numerical,
  title={Numerical approximation of partial differential equations},
  author={Quarteroni, Alfio and Valli, Alberto},
  volume={23},
  year={2008},
  publisher={Springer Science \& Business Media}
}

@article{geuzaine2009gmsh,
  title={Gmsh: A 3-D finite element mesh generator with built-in pre-and post-processing facilities},
  author={Geuzaine, Christophe and Remacle, Jean-Fran{\c{c}}ois},
  journal={International journal for numerical methods in engineering},
  volume={79},
  number={11},
  pages={1309--1331},
  year={2009},
  publisher={Wiley Online Library}
}

@article{isaac2020recursive,
  title={Recursive, parameter-free, explicitly defined interpolation nodes for simplices},
  author={Isaac, Tobin},
  journal={SIAM Journal on Scientific Computing},
  volume={42},
  number={6},
  pages={A4046--A4062},
  year={2020},
  publisher={SIAM}
}

@article{eskilsson2004triangular,
  title={A triangular spectral/hp discontinuous Galerkin method for modelling 2D shallow water equations},
  author={Eskilsson, Claes and Sherwin, Spencer J},
  journal={International Journal for Numerical Methods in Fluids},
  volume={45},
  number={6},
  pages={605--623},
  year={2004},
  publisher={Wiley Online Library}
}

@article{dubiner1991spectral,
  title={Spectral methods on triangles and other domains},
  author={Dubiner, Moshe},
  journal={Journal of Scientific Computing},
  volume={6},
  number={4},
  pages={345--390},
  year={1991},
  publisher={Springer}
}
\end{document}